\documentclass[11pt]{article}
\usepackage{amsmath,amssymb,amsthm,amscd}
\numberwithin{equation}{section} \textwidth 160mm \oddsidemargin
-1pt 

\def\p{\partial}
\def\b{\bar}

\def\cC{{\mathcal C}}
\def\cB{{\mathcal B}}
\def\cD{{\mathcal D}}

\def\cF{{\mathcal F}}

\def\cH{{\mathcal H}}
\def\cI{{\mathcal I}}

\def\cM{{\mathcal M}}
\def\cN{{\mathcal N}}
\def\cO{{\mathcal O}}

\def\cR{{\mathcal R}}
\def\cS{{\mathcal S}}

\def\cU{{\mathcal U}}
\def\cV{{\mathcal V}}
\def\cW{{\mathcal W}}

\def\bD{{\mathbf D}}
\def\bE{{\mathbf E}}

\def\NN{{\mathbb N}}

\def\RR{{\mathbb R}}
\def\CC{{\mathbb C}}

\def\dbd{\overline{\partial}\partial}

\newtheorem{prop}{Proposition}[subsection]
\newtheorem{theo}[prop]{Theorem}
\newtheorem{lem}[prop]{Lemma}
\newtheorem{cor}[prop]{Corollary}
\newtheorem{rem}[prop]{Remark}

\newtheorem{defi}[prop]{Definition}

\def\begeq{\begin{equation}}
\def\endeq{\end{equation}}

\def\and{\quad{\rm and}\quad}

\let\lra=\longrightarrow

\def\mapright\#1{\,\smash{\mathop{\lra}\limits^{\#1}}\,}

\newcommand{\integer}{\ensuremath{{\mathbf Z}}}
\newcommand{\naturals}{\ensuremath{{\mathbf N}}}

\newcommand{\complex}{\ensuremath{{\mathbf C}}}

\newcommand{\GL}[1]{\ensuremath{{\mathrm {GL}_{ #1 }}}}
\newcommand{\GLC}[1]{\GL{#1}(\complex)}

\newcommand{\BB}{{\mathcal B}}

\newcommand{\II}{\mathcal{I}}
\newcommand{\LL}{\mathcal{L}}

\newcommand{\kk}{{\mathbf{k}}}

\newcommand{\Grr}{\ensuremath{{\mathrm{Gr}}}}

\newcommand{\Grrinf}{\Grr(H^{(n)})}
\newcommand{\Grrr}[1]{\Grr^{(#1)}}

\newcommand{\LG}{{\LL} \GLC{n}}
\newcommand{\LGp}{{\LL}^+ \GLC{n}}
\newcommand{\LGm}{{\LL}^- \GLC{n}}

\title{Geometry of K\"ahler Metrics and Foliations by Holomorphic Discs}
\author{X. X. Chen and G. Tian\thanks{Both authors are supported by NSF research grants and the second
author was also partly supported by a J. Simons fund}}
\begin{document}
\ifx\href\undefined\else\hypersetup{linktocpage=true}\fi 
\bibliographystyle{plain}
\date{} 
\maketitle
\quad\quad\quad\quad\quad\quad{ Dedicated to Professor E. Calabi for his 80th birthday}

\tableofcontents

\section{Introduction and Main Results}
The purpose of this paper is to establish a completely new partial
regularity theory on certain homogeneous complex Monge-Ampere
equations. Our partial regularity theory will be obtained by
studying foliations by holomorphic curves and and their relations
to homogeneous complex Monge-Ampere equations. As applications, we
will prove the uniqueness of extremal K\"ahler metrics and give an
necessary condition for existence of extremal K\"ahler metrics.
Further applications will be discussed in our forthcoming papers.

\subsection{A brief Tour of Extremal K\"ahler Metrics}

Following \cite{calabi82}, we call a K\"ahler metric extremal if
the complex gradient of its scalar curvature is a holomorphic
vector field. In particular, any K\"ahler metric of constant
scalar curvature is extremal, conversely, if the underlying
K\"ahler manifold has no holomorphic vector fields, then an
extremal K\"ahler metric is of constant scalar curvature. It
follows from the standard Hodge theory that a K\"ahler metric of
constant scalar curvature becomes K\"ahler-Einstein if its
K\"ahler class is propotional to the first Chern class of the
underlying manifold.

In the 50's, E. Calabi proposed the problem of studying existence
of K\"ahler-Einstein metrics on compact K\"ahler manifolds with
definite first Chern class. In 1976, S. T. Yau solved the famous
Calabi conjecture. A corollary of this implies that any K\"ahler
manifold with vanishing first Chern class has a Calabi-Yau metric,
that is, a Ricci-flat K\"ahler metric. Around the same time, T.
Aubin and S. T. Yau independently proved existence of
K\"ahler-Einstein metrics on compact K\"ahler manifolds with
negative first Chern class. The remaining case is technically more
involved. In \cite{Tian90}, the second named author proved that a
complex surface with positive first Chern clas admits a
K\"ahler-Einstein metric if and only if its automorphism group is
reductive. For higher dimensions, he proved in \cite{Tian97} that
the existence of K\"ahler-Einstein metrics with positive scalar
curvature is equivalent to an analytic stability of the underlying
K\"ahler manifold. It remains open how this analytic stability is
related to certain algebraic stability from geometric invariant
theory (cf. \cite{Tian97}, \cite{Tian00},
\cite{Dona01},\cite{PhongSturm02}, \cite{Paultian04}). There has
not been much progress made on the existence of general extremal
metrics, even in the case of complex surfaces. One possible reason
is that the corresponding equation is highly nonlinear and of 4th
order. We will give an necessary condition for existence of
extremal metrics.

On the other hand, there have been many results on uniqueness of
extremal metrics. Using the maximum principle, E. Calabi observed
in 50's that K\"ahler-Einstein metrics of non-positive scalar
curvature are unique in their K\"ahler class. In \cite{Bando87},
Bando and Mabuchi proved uniqueness of K\"ahler-Einstein metric of
positive scalar curvature modulo holomorphic automorphisms. In
\cite{Tianzhu00}, X.H. Zhu and the second named author proved
uniqueness of K\"{a}hler-Ricci Solitons on any K\"{a}hler
manifolds with positive first Chern class. Following a suggestion
of Donaldson, the first named author proved in \cite{chen991}
uniqueness for K\"ahler metrics of constant scalar curvature in
any K\"ahler class when the underlying manifold has non-positive
first Chern class. In \cite{Dona001}, S. Donaldson proved that
K\"ahler metrics of constant scalar curvature are unique in any
rational K\"ahler class on any projective manifolds without
non-trivial holomorphic vector fields.\footnote{During preparation
of this paper, we learned from T. Mabuchi that he was able to
remove the assumption on non-existence of holomorphic vector
fields in the special case of projective manifolds.}

In this paper, we will prove
\begin{theo}
\label{th:uniqueness} Let $(M, \Omega)$ be a compact K\"ahler
manifold with a K\"ahler class $\Omega\in H^2(M,\RR)\cap
H^{1,1}(M,\CC)$. Then there is at most one extremal K\"ahler
metric with K\"ahler class $\Omega$ modulo holomorphic
transformations, that is, if $\omega_1$ and $\omega_2$ are
extremal K\"ahler metrics with the same K\"ahler class, then there
is a holomorphic transformation $\sigma$ such that
$\sigma^*\omega_1 =\omega_2$.
\end{theo}

In order to give a necessary condition for the existence of
K\"ahler metrics with constant scalar curvature, we recall the
K-energy introduced by T. Mabuchi in \cite{Ma87}: For any $\phi$
with $\omega+ \sqrt{-1}\partial\overline{\partial}\phi > 0$, which
will be denoted by $\omega_\phi$, define
$$ \bE_\omega(\phi)=-\int_0^1\int_M
\dot\phi(s(\omega_{\phi_t})-\mu)\omega_{\phi_t}^n\wedge dt,$$
where $\omega_{\phi_t}$ is any path of K\"ahler metrics joining
$\omega$ and $\omega_\phi$, $s(\omega_{\phi_t})$ denotes the
scalar curvature and $\mu $ is its average. It turns out that
$\mu$ is determined by the first Chern class $c_1(M)$ and the
K\"ahler class $[\omega]$.

\begin{theo}
\label{intro:th:minimum} Let $M$ be a compact K\"ahler manifold
with a constant scalar curvature K\"ahler metric $\omega$. Then
$\bE_\omega(\phi)\ge 0$ for any $\phi$ with $\omega_\phi
>0$.
\end{theo}

This theorem was proved for K\"ahler-Einstein metrics in
\cite{Bando87} (also see \cite{Tian97}) and in \cite{chen991} for
K\"ahler manifolds with non-positive first Chern
class.\footnote{After we finished the first draft of this paper,
we learned that S. Donaldson proved this theorem in the case of
projective manifolds and without holomorphic vector fields
\cite{Donaldson04}. His method is completely different from ours.}
This theorem can be also generalized to arbitrary extremal
K\"ahler metrics by using the modified K-energy. This theorem
gives a partial answer to a conjecture of the second author: $M$
has a constant scalar curvature K\"ahler metric in a given
K\"ahler class $[\omega]$ if and only if the K-energy is proper in
a suitable sense on the space of K\"ahler metrics with the fixed
K\"ahler class $[\omega]$. Combining Theorem
\ref{intro:th:minimum} with results in \cite{Tian00} and
\cite{Paultian04}, we can deduce

\begin{cor} Let $(M,L)$ be a polarized algebraic manifold, that is,
$M$ is algebraic and $L$ is a positive line bundle. If there is a
constant scalar curvature K\"ahler metric with K\"ahler class equal
to $c_1(L),$ then $(M,L)$ is asymptotically K-semistable or
CM-semistable in the sense of \cite{Tian97} (also see \cite{Tian00})
\footnote{According to \cite{Paultian04}, the CM-stability
(semistability) is equivalent to the K-stability (semistability).}.
\end{cor}

\subsection{Space of K\"ahler Metrics.}

Now let us explain how the above theorems can be proved. First let
us give a direct approach suggested by S. Donaldson and developed in
\cite{chen991}. Let $\omega$ be a fixed K\"ahler metric on $M$. It
follows from Hodge theory that the space of K\"ahler metrics with
K\"ahler class $[\omega]$ can be identified with the space of
K\"ahler potentials
\[
{\cH}_{\omega} = \{ \phi \mid \omega_{\phi} = \omega +
 \sqrt{-1} \p \b\p \phi > 0 \;{\rm on} \; M\}/ \sim,
\]
where $\phi_1 \sim \phi_2$ if and only if $\phi_1=\phi_2+ c$ for
some constant $c$. We will drop the subscript $\omega$ if no
possible confusion may occur. A tangent vector in $T_\phi\cH$ is
just a function ${\phi_0}$ such that
$$\int_M {\phi_0}\omega_\phi^n =0.$$ Its norm in the $L^2$-metric
on $\cH$ is given by (cf. \cite{Ma87})
\[
\|{\phi_0}\|^2_{\phi} =\int_{M}{\phi_0}^2\;\omega^n_\phi.
\]
A straightforward computation shows that the geodesic equation of
this $L^2$ metric is
\[
\phi''(t) - {1\over 2} \langle d \phi', d\phi'\rangle_\phi \; = \;0.
\]
where $\langle \cdot, \cdot\rangle_\phi$ denotes the natural inner
product on $T^*M$ induced by the K\"ahler metric $\omega_\phi$,
$\phi(t)\in \cH$ for any $t\in [0,1]$ and $\phi', \phi''$ denote
the derivatives of $\phi$ on $t$. Set $\phi (t,\theta,
x)=\phi(t)(x)$ for $t\in [0,1]$, $\theta\in S^1$ and $x\in M$,
then the path $\{\phi(t)\}$ satisfies the geodesic equation if and
only if the function $\phi$ on $[0,1]\times S^1\times M$ satisfies
the homogeneous complex Monge-Ampere equation
\begin{equation}
\label{eq:hcma0} (\pi_2^* \omega + \sqrt{-1}\partial
\overline{\partial }\phi)^{n+1} \;\; =\;\; 0, \qquad {\rm on} \;
\Sigma \times M,
\end{equation}
where $\Sigma = [0,1] \times S^1$ and $\pi_2: \Sigma\times M\mapsto
M$ and $\pi_1: \Sigma\times M\mapsto \Sigma$ are two orthogonal
projections. Using the convexity of the K-energy along geodesics of
the $L^2$-metric, S. Donaldson observed in \cite{Dona96} that both
Theorem 1.1.1 and 1.1.2 follow if one can show: For any
$\phi_0,\phi_1 \in \cH$, (\ref{eq:hcma0}) has a smooth solution
$\phi$ such that $\phi(t,\theta,\cdot)\in \cH$,
$\phi(0,\theta,\cdot)= \phi_0$ and $\phi(1,\theta,\cdot)=\phi_1$.
However, this turns out to be a very difficult problem and remains
open until now. In fact, one can consider (\ref{eq:hcma0}) over a
general Riemann surface $\Sigma$ with boundary condition $\phi
=\phi_0$ along $\partial \Sigma$, where $\phi_0$ is a smooth
function on $\partial\Sigma \times M$ such that $\phi_0(z,\cdot)\in
\cH$ for each $z\in
\partial \Sigma$.\footnote{we often regard $\phi_0$ as a smooth
map from $\p \Sigma$ into $\cH$.} It also has geometric meaning.
The equation (\ref{eq:hcma0}) can be regarded as the infinite
dimensional version of the WZW equation for maps from $\Sigma$
into $\cH$ (cf. \cite{Dona96}).\footnote{Original WZW equation is
for maps from a Riemann surface into a Lie group.} The following
theorem was proved by the first named author in \cite{chen991} and
will play a fundamental role in this paper.
\begin{theo} (\cite{chen991})
\label{th:chen991} For any smooth map $\phi_0: \p \Sigma
\rightarrow \cH$, there exists a unique $C^{1,1}$ solution $\phi$
of (\ref{eq:hcma0}) such that $\phi = \phi_0$ on $\p \Sigma$ and
$\phi(z,\cdot)\in \overline{\cH}$ for each $z\in
\Sigma$.\footnote{Here $\overline{\cH}$ denotes the closure of
$\cH$ in any $C^{1,\alpha}(\Sigma\times M)$-topology ($\forall
\;\alpha\in (0,1)$). A function $\varphi \in \b \cH$ if there exists a sequence of potential functions
$\{\varphi_m \in \cH, m \in \NN\}$ with uniform $C^{1,1}$ bound such that $\varphi_m \rightarrow \varphi$
in weak $C^{1,1}$ topology. }
\end{theo}
The lack of sufficient regularity for above solution is the
obstruction to proving Theorem 1.1.1 and 1.1.2 in general cases by
using geodesics. We should point out that complex Monge-Ampere
equations have been studied extensively (cf. \cite{CNS84},
\cite{CKNS85}, \cite{Bedford76} etc.). However, the regularity for
homogeneous complex Monge-Ampere equations beyond $C^{1,1}$ has
been missing. Indeed, there are examples in which some solutions
are only $C^{1,1}$, here is a simple example: Let $\Omega$ be the
unit ball in $\CC^2$ and define
\[ u = \left\{
\begin{array}{lcl} 0 & {\rm if} &\; |z_1|^2, |z_2|^2 \leq {1\over 2};
\\ ({1\over 2}-|z_1|^2)^2 &{\rm if} & \;\; |z_1|^2 \geq {1\over 2};
\\ ({1\over 2}-|z_2|^2)^2 &  {\rm if} &\;\; |z_2|^2 \geq {1\over 2};
\end{array}\right.
\]
then $(\partial\overline{\partial} u)^2=0$ on $\Omega$ and
$u|_{\p\Omega}$ is smooth, but $u$ is only $C^{1,1}$.

Therefore, it may not be possible to have smooth solutions for
(\ref{eq:hcma0}) with general boundary values. Of course, our
boundary values are special and they satisfy certain positivity
condition, that is, any boundary function $\phi_0$ takes values in
$\cH$. This fact will play a crucial role in our approach.

\subsection{An Approach to Uniqueness and lower bound for the K-energy}

Our proof of Theorem 1.1.1 and 1.1.2 starts with the following
observations: Given two functions $\phi_0$ and $\phi_1$ in $\cH$, if
there is a sufficiently nice solution $\phi$, which is not
necessarily smooth, of (\ref{eq:hcma0}) on $[0,1]\times\RR \times
M$, then the evaluation function $f=\bE (\phi(z,\cdot))$ is bounded,
subharmonic\footnote{This follows from sub-harmonic property of the
K-energy with respect to almost smooth solution as shown in the
subsequent computation.} and constant along each boundary component
of $[0,1]\times \RR$. Furthermore, if $\phi_0$ is a critical metric
of $\bE$, then it follows from the Maximum principle that
$\bE(\phi_1)\ge \bE(\phi_0)$ and equality holds if and only if each
$\phi(z,\cdot)$ is a critical metric of $\bE$. The infinite strip
$[0,1]\times \RR$ can be approximated by bounded disks. Hence, if
(\ref{eq:hcma0}) has a uniformly bounded and sufficiently smooth
solution when $\Sigma$ is a disk, then Theorem 1.1.1 and 1.1.2
follow.

In this paper, we will establish a partial regularity theory for
(\ref{eq:hcma0}) in the case that $\Sigma$ is a unit disk. This
seems to be the first partial regularity theory for homogeneous
complex Monge-Ampere equations.

First let us introduce some notations. Suppose that $\phi$ is a
$C^{1,1}$ solution of (\ref{eq:hcma0}), we denote by $\cR_\phi$
the set of all $(z,x)\in \Sigma\times M$ near which $\phi$ is
smooth and $\omega_{\phi(z',\cdot)}=\omega
+\sqrt{-1}\partial\overline{\partial} \phi(z',\cdot)$ is a
K\"ahler metric. We may regard $\cR_\phi$ as the regular set of
$\phi$. It is open, but {\it aprior}, it may be empty. We have a
distribution $\cD_\phi\subset T(\Sigma\times M)$ over $\cR_\phi$:
\begin{eqnarray}
\label{eq:cDovercR} \cD_\phi|_{(z,x)} = \{ v \in T_z\Sigma\times
T_xM~|~i_v \left (\pi^*_2\omega + \sqrt{-1} \p\overline{\p}
\phi\right ) =0\}, ~~~(z,x)\in \cR_\phi.
\end{eqnarray}
Here $i_v$ denotes the interior product. Since the form is closed,
$\cD_\phi$ is integrable. We say that $\cR_\phi$ is saturated in
$\cV\subset \Sigma\times M$ if every maximal integral sub-manifold
of $\cD_\phi$ in $\cR_\phi\cap \cV$ is a disk and closed in $\cV$.
By nature of product manifold, we may write any vector in $\cD_\phi$
as (after re-scaling)
\begin{equation}
{\p \over {\p z}} + X \in \cD_\phi\mid_{(z,x)},\qquad {\rm
where}\;\; X \in T^{1,0}_x M. \label{def:leafvector}
\end{equation}
\begin{defi}
\label{def:partiallysmooth} A solution $\phi$ of (\ref{eq:hcma0})
is called partially smooth if it is $C^{1,1}$-bounded on
$\Sigma\times M$ and $\cR_\phi$ is open and saturated in
$\Sigma\times M$, but dense in $\p \Sigma \times M$, such that the
varying volume form $\omega_{\phi(z,\cdot)}^n$ extends to a
continuous $(n,n)$ form on $\Sigma^0\times M$, where
$\Sigma^0=(\Sigma\backslash \partial \Sigma)$.
\end{defi}

Clearly, if $\phi$ is a partially smooth solution, then its
regular set $\cR_\phi$ consists of all points where the vertical
volume form $\omega_{\phi(z,\cdot)}^n$ is positive in
$\Sigma\times M$.

\begin{theo} \label{th:partiallysmooth}
Suppose that $\Sigma$ is a unit disk. For every smooth map
$\phi_0: \p \Sigma \rightarrow \cH$, there exists a unique
partially smooth solution to (\ref{eq:hcma0}).
\end{theo}

In particular, this theorem implies that the $C^{1,1}$ solution
from \cite{chen991} is smooth on some open subset of $\Sigma\times
M$ which is also dense in $\p \Sigma\times M.\;$ We expect this
density property holds in interior as well. Indeed, it is true for
generic boundary values, moreover, we can estimate the size of the
singular set $\cS_\phi= \Sigma\times M\backslash \cR_\phi$.

\begin{defi}
\label{def:almostsmooth} We say that a solution $\phi$ of
(\ref{eq:hcma0}) is almost smooth if
\begin{enumerate} \item it is partially smooth,
\item The distribution ${\cal D}_\phi$ extends to a continuous
distribution in a saturated set $\tilde \cV\subset \Sigma \times
M$, such that the complement $\tilde \cS_\phi$ of $\tilde \cV$ is
locally extendable\footnote{A closed subset $S\subset \Sigma
\times M $ of measure $0$ is {\it locally extendable} if for any
continuous function in $\Sigma\times M$ which is $C^{1,1}$ on
$\Sigma\times M \setminus S$ can be extended to a $C^{1,1}$
function on $\Sigma\times M$. Notice that any set of codimension 2
or higher is automatically locally extendable.} and $\phi$ is
$C^1$ continuous on $\tilde \cV$. The set $\tilde \cS_\phi$ is
referred as the singular set of $\phi$. \item The leaf vector
field $X$ is uniformly bounded in $\cD_\phi.\;$
\end{enumerate}
\end{defi}

A smooth solution is certainly an almost smooth solution of
(\ref{eq:hcma0}). For a sequence of almost smooth solutions whose
boundary values converge in certain smooth topology, then it
converges to a partially smooth solution in $C^{1,\beta}$-topology
($0< \beta < 1$).

\begin{theo}
\label{th:almostsmooth} Suppose that $\Sigma$ is a unit disk. For
any $C^{k,\alpha}$ map $\phi_0: \p \Sigma \rightarrow \cH$ ($k \ge
2$, $0< \alpha < 1$) and for any $\epsilon
> 0$, there exists a $\phi_\epsilon: \p\Sigma\rightarrow \cH$
in the $\epsilon$-neighborhood of $\phi_0$ in
$C^{k,\alpha}(\Sigma\times M)$-norm, such that (\ref{eq:hcma0})
has an almost smooth solution with boundary value $\phi_\epsilon$.
\end{theo}

This partial regularity is sharp since we do have examples (at
least implicitly) where the solution for (\ref{eq:hcma0})
corresponds to a singular metric (cf. \cite{Dona01}).

Theorem \ref{th:partiallysmooth} will be derived from Theorem
\ref{th:almostsmooth} by using estimates developed in late
sections.\footnote{It is believed that for any smooth boundary map
$\phi_0:\p \Sigma \mapsto \cH$, the corresponding $C^{1,1}$
solution is almost smooth. It is also interesting to estimate
precise size of $\cS_\phi$.} The importance of Theorem
\ref{th:almostsmooth} lies in the following theorem, which implies
that an almost smooth solution is as good as a smooth solution
when one concerns the K-energy.
\begin{theo}
\label{th:hessianofK-energy} Suppose that $\phi$ is a partially
smooth solution to (\ref{eq:hcma0}). For every point $z \in
\Sigma, $ let $\bE(z)$ be the K-energy (or modified K energy)
evaluated on $\phi(z,\cdot)\in \overline{\cH}$. Then $\bE$ is a
bounded subharmonic function on $\Sigma$ in the sense of
distribution, moreover, we have the following
\[\int_{\cR_\phi} | {\bD}
{{\p \phi}\over {\p z}}|_{{\omega_{\phi(z,\cdot)}}}^2\,
\frac{\sqrt{-1}}{2} dz\wedge d\bar z \wedge
{\omega_{\phi(z,\cdot)}}^n\, d\,z d\, \b z \leq \displaystyle
\int_{\p \Sigma}\;{{\p \bE}\over{\p {\bf n}}} \big |_{\p \Sigma}
ds,
\]
where $ds$ is the length element of $\p \Sigma$ and for any smooth
function $\theta$, $\bD \theta$ denotes the (2,0)-part of
$\theta$'s Hessian with respect to the metric
$\omega_{\phi(z,\cdot)}$. The equality holds if $\phi$ is almost
smooth.
\end{theo}
Theorem \ref{intro:th:minimum} follows from this theorem and the
observation on subharmonic functions over an infinite strip at the
beginning of this subsection.

Another corollary of this theorem is
\begin{prop}
\label{prop:geodesics} If there are two constant scalar curvature
metrics (resp. two extremal K\"ahler metrics), there exists a path
in $\overline{\cH}$ of $C^{1,1}$-functions $\phi_t$ ($0\le t\le
1$) which connects those two metrics, such that the K-energy
(resp. modified K energy) achieves its minimum at every $\phi_t$
along the path.
\end{prop}

It was conjectured by the first named author that any $C^{1,1}$
function which attains the absolute minimum of the K-energy (or
the modified K-energy) is actually smooth (cf. \cite{chen991},
\cite{chen993}). In this paper, we will confirm this conjecture in
the case that the function arises from Proposition
\ref{prop:geodesics}. This is sufficient for completing the proof
of Theorem \ref{th:uniqueness}, that is, the uniqueness theorem of
extremal K\"ahler metrics in any given K\"ahler class. The
conjecture will be proved by first establishing a partial
$C^1$-regularity for the vertical volume form of any $C^{1,1}$
K-energy minimizers.

The major part of this paper will be devoted to proving these main
technical results as we described above. We believe that our
techniques developed in this paper can be applied to studying
regularity of a much wider class of degenerate complex Monge-Ampere
equations.

\subsection{Ideas for Proof of Theorem \ref{th:almostsmooth}}

It has been known for a long time that solutions of the
homogeneous complex Monge-Ampere equation are closely related to
foliations by holomorphic curves (cf. \cite{Lempt83},
\cite{Semmes92}, \cite{Dona01}). In \cite{Semmes92}, S. Semmes
formulated the Dirichlet problem for (\ref{eq:hcma0}) in terms of
a foliation by holomorphic curves with boundary in a totally real
submanifold of the complex cotangent bundle of the underlying
manifold.

In \cite{Semmes92} (also see \cite{Dona01}), Semmes constructed a
complex manifold $\cW_{[\omega]}$ with a holomorphic (2,0)-form $
\Theta)$ to each K\"ahler class $[\omega].\; $ There is  a natural
projection $\pi: \cW_{[\omega]}\mapsto M$. He observed that for
any $\phi\in \cH$, we can associate a complex submanifold
$\Lambda_\phi$ in $\cW_{[\omega]}$ such that
\begin{equation}
\label{eq:sympl-lang} \Theta|_{\Lambda_\phi} = -\sqrt{-1}
\omega_\phi,
\end{equation}
that is, ${\rm Re}(\Theta)|_{\Lambda_\phi}=0$ and $-{\rm
Im}(\Theta)|_{\Lambda_\phi} = \omega_\phi > 0$. Locally,
$\Lambda_\phi$ is simply the graph of $\p(\rho + \phi)$ where
$\omega = \sqrt{-1}\p\b{\p} \rho$. This means that
${\Lambda_\phi}$ is an exact Lagrangian symplectic submanifold of
$\cW_{[\omega]}$ with respect to $\Theta$. Conversely, given an
exact Lagrangian symplectic submanifold $\Lambda$ of
$\cW_{[\omega]}$, one can construct a smooth function $\phi$ such
that $\Lambda=\Lambda_\phi$. Hence, K\"ahler metrics in the
K\"ahler class $[\omega]$ are in one-to-one correspondence with
exact Lagrangian symplectic submanifolds of $\cW_{[\omega]}$.

Given $\phi_0: \partial \Sigma \mapsto \cH$, define
\begin{equation}
 \b \Lambda_{\phi_0} = \{ (\tau, v) \in
\partial \Sigma\times \cW_{[\omega]}~|~v\in
\Lambda_{\phi_0(\tau)}~\}. \label{intro:boundarydefi}
\end{equation}
One can show that $\b \Lambda_{\phi_0}$ is a totally real
sub-manifold in $\Sigma\times\cW_{[\omega]}$. So it makes sense to
study the {\it moduli} space of holomorphic disks with boundary in
$\b \Lambda_{\phi_0}$. Its significance is clear from the
following result from \cite{Semmes92} (also see \cite{Dona01}).
\begin{prop}
\label{prop:semmes92} Assume that $\Sigma$ is simply connected.
For any boundary map $\phi_0: \p \Sigma \rightarrow \cal H,\;$
there is a solution $\phi \in C^{\infty}(\Sigma, {\cal H}) $ of
(\ref{eq:hcma0}) with boundary value $\phi_0$ if and only if there
is a smooth family of holomorphic maps $h_x:\Sigma\mapsto
\cW_{[\omega]}$ parametrized by $x\in M$ satisfying: (1)
$\pi(h_x(z_0))=x$, where $z_0$ is a given point in
$\Sigma\backslash \partial \Sigma$; (2) $h_x(\tau)\in
\Lambda_{\phi_0(\tau)}$ for each $\tau \in \partial \Sigma$ and
$x\in M$; (3) For each $z\in \Sigma$, the map
$\gamma_z(x)=\pi(h_x(z))$ is a diffeomorphism of $M$.
\end{prop}

In \cite{Dona01}, S. Donaldson used this fact to study
deformations of smooth solutions for (\ref{eq:hcma0}) as the
boundary value varies. This inspired us to study foliations by
holomorphic disks in order to have a partial regularity theory for
(\ref{eq:hcma0}). Theorem \ref{th:almostsmooth} will be proved by
establishing existence of foliations by holomorphic disks with
relatively mild singularities. More precisely, we will show that
for a generic boundary value, there is an open set in the moduli
space of holomorphic disks which generates a foliation on
$\Sigma\times M\backslash S$ for a closed subset $S$ of
codimension at least one.

Now let us fix a generic boundary value $\phi_0$ and study the
corresponding {\it moduli} $\cM_{\phi_0}$ of holomorphic disks.
First it follows from the Index theorem that the expected
dimension of this moduli is $2n$. Recall that a holomorphic disk
$u$ is regular if the linearized $\overline\partial$-operator
$\overline\partial _u$ has vanishing cokernel. The moduli space is
smooth near a regular holomorphic disk. Following \cite{Dona01},
we call $u$ super-regular if there is a basis $s_1,\cdots ,s_{2n}$
of the kernel of $\overline\partial _u$ such that
$d\pi(s_1)(x),\cdots, d\pi(s_{2n})(x)$ span $T_{u(x)}M$ for every
$x\in \Sigma$, where $\pi: \Sigma\times \cW_{[\omega]} \mapsto
\Sigma\times M$ is the natural projection. We call $u$ almost
super-regular if $d\pi(s_1)(x),\cdots, d\pi(s_{2n})(x)$ span
$T_{u(x)}M$ for every $x\in \Sigma\backslash
\partial \Sigma$. Clearly, the set of super-regular disks is open.

One of our crucial observations is that Semmes' arguments can be
made local along super-regular holomorphic disks. First let us
introduce the notion of a nearly smooth foliation.

\begin{defi}
\label{def:nearlyfoliation} A nearly smooth foliation
$\cF_{\phi_0}$ associated to a boundary value $\phi_0$ is given by
an open subset $\cU_{\phi_0}\subset\cM_{\phi_0}$ of super-regular
disks whose images in $\Sigma\times M$ give rise to a foliation on
an open-dense set $\cV_{\phi_0}$ of $\Sigma\times M$ such that
\begin{enumerate}
\item this foliation can be extended to be a continuous foliation
by holomorphic disks in an open set
$\tilde\cV_{\phi_0}\subset(\Sigma\backslash\p\Sigma)\times M $
such that it admits a continuous lifting in $\Sigma\times \cW_M$;
\item the complement of $\tilde \cV_{\phi_0}$ in $\Sigma\times M$
is locally extendable. \item The leaf vector field (cf. definition
below) induced by the foliation in $\cV_{\phi_0}$ is uniformly
bounded.
\end{enumerate}
\end{defi}

\begin{defi} \label{def:leafvector1}  For each $(z,x) \in \cV_{\phi_0}$, the leaves of $\cV_{\phi_0}$ give rise to a vector
field of the form
\[
{\p \over {\p z}} + X,\qquad {\rm where} \;\; X \in
T_{z,x}^{1,0}M.
\]
This $X$ is called the leaf vector field in $\cV_{\phi_0}.\;$
\end{defi}

Proposition \ref{prop:semmes92} has the following generalization.

\begin{theo}
\label{th:almostsmooth=nearlyfoliation} Almost smooth solutions of
(\ref{eq:hcma0}) are in one-to-one correspondence with nearly
smooth foliations. Moreover, if $\phi_0$ is generic, the
corresponding almost smooth solution $\phi$ has additional
properties: $\omega_\phi$ is a smooth $(1,1)$ form in
$\Sigma\times M \setminus \tilde \cS_\phi$ and the singular set
$\tilde \cS_\phi$ has codimension at least 2 in each slice
$\{z\}\times M, \;\forall \; z \in \Sigma^0$.\footnote{The
corresponding nearly smooth foliations have additional properties
and will be called almost super-regular foliations (cf. Section
3.3).}
\end{theo}

Thus, in order to prove Theorem \ref{th:almostsmooth}, we only
need to show the following

\begin{theo}
\label{th:nearlyfoliation} For a generic boundary value $\phi_0$,
there is a nearly smooth foliation associated to $\phi_0$
generated by an open set $\cU_{\phi_0}$ of the corresponding
moduli space $\cM_{\phi_0}$. Moreover, the set of holomorphic
disks which are neither super-regular nor almost super-regular has
codimension at least two in the closure of $\cU_{\phi_0}$ in
$\cM_{\phi_0}$.
\end{theo}

The idea for proving Theorem \ref{th:nearlyfoliation} is outlined
as follows. Let $\phi_0$ be a generic boundary value such that
$\cM_{\phi_0}$ is smooth. This follows from a result of Oh on
transversality. By the same transversality argument, one can show
that there is a generic path $\phi_t$ ($0\le t\le 1$) such that
$\phi_1=0$ and the total moduli $\tilde \cM = \bigcup_{t\in [0,1]}
\cM_{\phi_t}$ is smooth. Moreover, we may assume that
$\cM_{\phi_t}$ are smooth for all $t$ except finitely many
$t_1,\cdots ,t_N$ where the {\it moduli} space may have isolated
singularities. It follows from Semmes and Donaldson's
work---Proposition \ref{prop:semmes92} that 
 $\cM_1$ has at least one component which
gives a foliation for $\Sigma\times M$. We want to show that this
component will deform to a component of $\cM_{\phi_0}$ which
generates a nearly smooth foliation. We will use the continuity
method. Assume that $\phi$ is the unique $C^{1,1}$-solution of
(\ref{eq:hcma0}) with boundary value $\phi_t$ for some $t\in
[0,1]$. Let $f$ be any holomorphic disk in the component of
$\cM_{\phi_t}$ which generates the corresponding foliation.

Using the $C^{1,1}$ bound on $\phi$, one can have a uniform
area\footnote{We actually calculate area of the image of disks in
$\Sigma\times M.\;$} bound on holomorphic disks in $\cM_{\phi_t}$.
It follows from an extension of Gromov's compactness theorem that
any sequence of such holomorphic disks has a subsequence which
converges to a holomorphic disk together with possibly finitely
many bubbles. These bubbles which occur in the interior are
holomorphic spheres, while bubbles in the boundary might be
holomorphic spheres or disks. We will show that no bubbles can
actually occur.  According to E. Calabi and X. Chen
\cite{chen992}, this infinite dimensional space $\cH$ is
non-positively curved in the sense of Alexanderov. Heuristically
speaking, we can exploit this curvature condition to rule out the
existence of interior bubbles. One can also rule out boundary
bubbles  by using the non-positivity and totally real property of
the boundary condition. Since there are no bubbles, the Fredholm
index of holomorphic disks is invariant under the limiting
process. This is an important fact needed in our doing deformation
theory.

In order to get a nearly smooth foliation, we need to prove that
the {\it moduli} space has an open set of super-regular
holomorphic disks for each $t$. First we observe that the set of
super- regular disks is open. Moreover, using the transversality
arguments, one can show that for a generic path ${\phi}_t$, the
closure of all super-regular disks in each $\cM_{{\phi}_t}$ is
either empty or forms an irreducible component. This implies the
openness. It remains to prove that each {\it moduli} has at least
one super-regular disk. It is done by using capacity estimates and
curvature estimate along super-regular holomorphic disks (cf.
Sections 4 and 5 for details).

\subsection{Organization}
In Section 2, we establish the correspondence between homogeneous
complex Monge-Ampere equations and foliations by holomorphic
curves. The goal is to prove Theorem
\ref{th:almostsmooth=nearlyfoliation}. The proof is based on a
local version of Semmes' construction. Semmes's construction is
global in nature and was re-diskovered in Donaldson's work
\cite{Dona01}. In Section 3, we show necessary transversality
results. In particular, we show that the set of boundary values
such that corresponding {\it moduli} space $\cM$ induces an almost
super regular foliation is generically open. In Section 4, we
study the defomation of holomorphic disks arising from a smooth
solution to a homogenous complex Monge-Ampere equation. This is a
local theory which is used in Section 2,3 and late sections as
well. In Section 5, we prove the set of boundary values such that
corresponding {\it moduli} space $\cM$ induces an almost super
regular foliation is closed. This will be done by proving a volume
ratio estimate via a capacity argument.  In Section 5, we prove
that the K energy function is subharmonic when restricted to a
disk family of almost smooth solutions, which in turns implies
that the K energy function is always bounded from below. In
section 7, we derive a partial $C^1$-regularity for the vertical
volume form of any $C^{1,1}$ K-energy minimizer.  We need to
introduce a notion of weak K\"ahler Ricci flow to derive this {\it
a prioir} estimate. In Section 8, we prove the uniqueness result
for extremal K\"ahler metrics.  In the Appendix, we discuss the
Loop
 groups of $GL(n,C)$ and their relation to holomorphic discs.  This is
essentially written by Professor E. Lupercio.

\indent {\bf Acknowledgement}  The first author is grateful to S.
K. Donaldson for insightful discussions and constant
encouragements during this project and over years, in particular,
the first author is deeply indebted to him during the period of
time when both of us were visiting Stanford University in 1998.

\section{Foliations and the Homogenous comlex Monge Ampere Equation}

In this section, we diskuss the correspondence between homogeneous
complex Monge-Ampere equations and foliations by holomorphic
disks. We will prove Theorem
\ref{th:almostsmooth=nearlyfoliation}.

\subsection{Semmes' construction}

In \cite{Semmes92}, Semmes associated a complex manifold
$\cW_{[\omega]}$ to each K\"ahler class $[\omega]$: Let $\{U_i, i
\in \cI\}$ be a covering of $M$ such that $\omega|_{U_i} =
\sqrt{-1}\partial\overline{\partial} \rho_i,\;$ where $\cI$ is an
index set. For any $x=y \in U_i\cap U_j (i, j\in \cI),$ we
identify $(x,v_i)\in T^*U_i$ with $(y, v_j)\in T^*U_j$ if $v_i=v_j
+
\partial (\rho_i-\rho_j).\;$ Then $\cW_{[\omega]}$ consists of all
these equivalence classes of $[x,v_i]$. There is a natural map
$\pi: \cW_{[\omega]}\mapsto T^*M$, assigning $(x,v_i)\in T^*U_i$
to $(x, v_i - \partial \rho_i)$. Then the complex structure on
$T^*M$ pulls back to a complex structure on $\cW_{[\omega]}$,
moreover, there is also a canonical holomorphic 2-form $\Theta$ on
$\cW_{[\omega]}$, in terms of canonical local coordinates
$z_\alpha, \xi_\alpha$ ($\alpha = 1,\cdots, n$) of $T^*U_i$,
$$\Theta = d z_\alpha \wedge d\xi_\alpha .$$
Now for any $\phi\in \cH$, we can associate a complex submanifold
$\Lambda_\phi$ in $\cW_{[\omega]}$: For any open subset $U$ on
which $\omega$ can be written as $\sqrt{-1}\partial\overline
\partial \rho$, we define $\Lambda_\phi|_U\subset \cW_{[\omega]}$ to be the graph of
$\partial (\rho + \phi)$ in $T^*U$. Clearly, this $\Lambda_\phi$
is independent of the choice of $U$. A straightforward computation
shows
\begin{equation}
\label{eq:sympl-lang2} \Theta|_{\Lambda_\phi} = -\sqrt{-1}
\omega_\phi,
\end{equation}
that is, ${\rm Re}(\Theta)|_{\Lambda_\phi}=0$ and $-{\rm
Im}(\Theta)|_{\Lambda_\phi} = \omega_\phi > 0$. This means that
${\Lambda_\phi}$ is an exact Lagrangian symplectic submanifold of
$\cW_{[\omega]}$ with respect to $\Theta$. Conversely, given an
exact Lagrangian symplectic submanifold $\Lambda$ of
$\cW_{[\omega]}$, we have a smooth function $\phi$ such that
$\Lambda=\Lambda_\phi$. Hence, K\"ahler metrics with K\"ahler
class $[\omega]$ are in one-to-one correspondence with exact
Lagrangian symplectic submanifolds in $\cW_{[\omega]}$.

This was briefly diskussed in our introduction. We refer the
readers to both \cite{Semmes92} and \cite{Dona01} for more
details. For the readers' convenience, let us   briefly explain
the proof of Proposition \ref{prop:semmes92}. Let $\phi$ be a
solution of (\ref{eq:hcma0}) on $\Sigma\times M$ such that $\phi
(z,\cdot)\in \cH$ for any $z\in \Sigma$. Recall that there is an
induced distribution $\cD_\phi\subset T(\Sigma\times M)$ by
\begin{equation}
\label{eq:distribution} \cD_\phi|_p = \{v\in T_p(\Sigma\times
M)~|~i_v(\pi_2^*\omega +\sqrt{-1}\partial\overline{\partial}
\phi)=0~\},~~~p\in \Sigma\times M.
\end{equation}
It is a holomorphic integrable distribution. If $\Sigma$ is
simply-connected and $\phi(z,\cdot)\in \cH$ for each $z\in
\Sigma$, then the leaf of $\cD_\phi$ containing $(z_0,x)$ is the
graph of a holomorphic map $f_x:\Sigma\mapsto M$ with
$f_x(z_0)=x$. If we write $f_x(z)=\sigma_z(x)$ we get a family of
diffeomorphisms $\sigma_z$ of $M$ with $\sigma_{z_0}={\rm Id}_M$.
Now for any fixed $z$ we have a K\"ahler form $\omega +
\sqrt{-1}\partial\overline\partial \phi(z,\cdot)$ on $M$ and hence
a section $s_z:M\mapsto \cW_{[\omega]}$ whose image is the exact
Lagrangian symplectic graph $\Lambda_{\phi(z,\cdot)}$. Then
$h_x(z)=\gamma_z(x)=s_z(f_x(z))$ as required. This process can be
reversed. Since we have to carry out this reversed process in the
proof of Theorem \ref{th:almostsmooth=nearlyfoliation}, we omit
details here and refer the readers to the next subsection.

\subsection{Local Uniqueness for Equation (\ref{eq:hcma0})}
We will prove uniqueness of compatible solutions for
(\ref{eq:hcma0}) near a super-regular disk.  We will define a
compatible solution as follows.

Given a boundary value $\phi_0$ on $\p\Sigma\times M$. Suppose
that $\cF_{\phi_0}$ is a nearly smooth foliation (cf. Definition
\ref{def:nearlyfoliation}). A solution $\phi$ of (\ref{eq:hcma0})
with boundary value $\phi_0$ in an open subset $\cV\subset\Sigma
\times M$ is called {\it compatible} with this foliation
$\cF_{\phi_0}$ if
\begin{enumerate}
\item  $\cV$ is saturated with respect to $\cF_{\phi_0}$; \item
$\omega + \sqrt{-1}\partial \bar \partial \phi(z,\cdot)$ is a
family of K\"ahler metrics on $M$; \item The kernel of $\pi_2^*
\omega + \sqrt{-1}\partial \bar
\partial \phi$ lies in the distribution induced by
$\cF_{\phi_0}$.
\end{enumerate}
\begin{theo}
\label{th:localuniqueness} Two compatible solutions of
(\ref{eq:hcma0}) with respect to the nearly smooth foliation
$\cF_{\phi_0}$ coincides along the intersection of their domains.
\end{theo}

We will adopt the notations from previous sections. First we
recall the integrable distribution
\begin{eqnarray}
\label{eq:cDovercR2} \cD_\phi|_{(z,x)} = \{ v \in T_z\Sigma\times
T_xM~|~i_v \left (\pi^*_2\omega + \sqrt{-1}\p\overline{\p}
\phi\right ) =0\}, ~~~(z,x)\in \cR_\phi.
\end{eqnarray}
Here $i_v$ denotes the interior product. Since $\cR_\phi$ is
saturated, every maximal integral submanifold of $\cD_\phi$ in
$\cR_\phi$ is a disk and closed in $\Sigma\times M$.

\begin{lem} \label{defo:localpotentialisharmonic}
For any $f \in  \cU_{\phi_0}, $  suppose that  ${\cal O}_f$ is a
saturated open neighborhood of the image of $f$ in $\Sigma\times
M$.  Suppose that $\phi_f$ is a solution of (\ref{eq:hcma0}) on
${\cal O}_f$ compatible with $\cF_{\phi_0}$. Then, for any $\tilde
f:\Sigma \mapsto \cV_{\phi_0}$ such that the image of $\pi\cdot
\tilde f$ lies completely in ${\cal O}_f$, we have
\begin{equation}
{\p^2 \over {\p z \; \p \b z}} \left( \phi _f (\pi\circ \tilde f
(z)) \right)
  =  - |\p (\pi\circ \tilde f )|_\omega^2( \pi\circ \tilde f (z)),\qquad \forall\;z \in \Sigma, \label{eq:folid1}
\end{equation}
and
\begin{equation} \phi _f(\pi\circ \tilde f (z)) = \phi_0 (\pi\circ \tilde f (z)),
\qquad \forall\;\; z \in \p \Sigma. \label{eq:folid2}
\end{equation}
\end{lem}
This lemma implies that $\phi_f$ is uniquely determined by only
the geometric conditions along the image of each leaf. Theorem
\ref{th:localuniqueness}
follows from this lemma.\\

Using this theorem, we can patch together $\phi_f$, initially
defined on small open, saturated sets around a super regular leaf
, to obtain a smooth potential function $\phi$ in $\cU_{\phi_0}$.
This function $\phi$ solves (\ref{eq:hcma0}) and satisfies
(\ref{eq:folid1}) - (\ref{eq:folid2}).

\subsection{Almost Smooth Solutions
$\Leftrightarrow$ Nearly Smooth Foliations}

In this subsection, we establish the equivalence between almost
smooth solutions of (\ref{eq:hcma0}) and nearly smooth foliations.
This generalizes Semmes' construction. We will adopt notations
from previous subsections.

\begin{prop}
\label{prop:almostsmooth>nearlyfoliation} An almost smooth
solution to the equation (\ref{eq:hcma0}) with boundary value
$\phi_0$ induces a nearly smooth foliation associated to $\phi_0$.
\end{prop}
\begin{proof} Let $\phi$ be an almost smooth solution with boundary value $\phi_0:\p\Sigma
\mapsto \cH$. For every point $(z,x) \in {\cR}_\phi$, there is a
unique holomorphic map $f\in \cM_{\phi_0}$ whose corresponding map
$\pi\circ f: \Sigma \mapsto \Sigma\times M$ passes through
$(z,x)$. The property that ${\cR}_\phi$ is saturated implies that
$\pi\circ f$ is a holomorphic disk and extends to the boundary of
$\Sigma\times M$. According to Donaldson\cite{Dona01}, such a
holomorphic disk is super-regular. All these super-regular disks
from $\cR_{\phi_0}$ give rise to this open set $\cU_{\phi_0}
\subset \cM_{\phi_0}.\;$ The other two conditions of a nearly
smooth foliation can be verified in a straightforward fashion. In
other words, an almost smooth solution indeed induces a nearly
smooth foliation $\cF_{\phi_0}$.
\end{proof}

Theorem \ref{th:almostsmooth=nearlyfoliation} follows from the
above proposition and the following.

 \begin{theo}
 If $\cF_{\phi_0}$ is a nearly smooth foliation (cf. Definition
 \ref{def:nearlyfoliation}) associated to
 a boundary value $\phi_0:\p \Sigma \rightarrow \cH$, then
 there is an almost smooth solution $\phi$ to (\ref{eq:hcma0}) with boundary value $\phi_0$.
 \end{theo}

The rest of this subsection is devoted to proving this theorem.
Let $\cU_{\phi_0}$ and $\cV_{\phi_0}$ be the open subsets of
$\cM_{\phi_0}$ and $\Sigma\times M$ as given by $\cF_{\phi_0}$. By
definition, the induced foliation in $\cV_{\phi_0}$ can be
extended to be  a continuous foliation by holomorphic disks in an
open and dense subset $\tilde \cV_{\phi_0}$ such that it admits a
continuous lifting to $\Sigma\times \cW_M.\;$ Moreover,
$\Sigma\times M \setminus \tilde \cV_{\phi_0}$ is {\it locally
extendable}.

\begin{prop}\label{prop:omegaform} There is a smooth family of non-degenerate,
closed (1,1) forms $\tilde \omega(z,\cdot)$ defined on
$\cV_{\phi_0}\cap\{z\}\times M$ and a closed (1,1) form $\Omega$
in $\cV_{\phi_0}$ such that
\begin{enumerate} \item $\tilde \omega =\omega_{\phi_0}$ in $\p \Sigma\times M, $
wherever $\tilde\omega$ is defined; \item  The restriction of
$\;\tilde \omega$ to each leaf is a constant form; \item $\Omega$
is defined by the following conditions:
\[ \Omega\mid_{\{z\}\times M} = \tilde{\omega}, \qquad
{\rm and}\qquad i_{{\p \over {\p z}} + X} \Omega = 0,\] where $X$
is the leaf vector field in $\cV_{\phi_0}$ induced by the {\it
nearly smooth foliation} (cf. Defi. \ref{def:leafvector1}).
\end{enumerate}
\end{prop}
\begin{proof} This is a local theorem. The proof can be found in
\cite{Semmes92} \cite{Dona01}.\end{proof}

Next we want to show that $\Omega=\pi_2^*\omega_0 + i \dbd \phi$
for a {\it compatible} function $\phi$ in $\Sigma \times M.\;$ We
want to define this potential function in a small open and
saturated neighborhood of any super regular leaf. Theorem
\ref{th:localuniqueness} implies that two different locally
defined {\it compatible} solutions of (\ref{eq:hcma0}) must agree
with each other on the overlap of their domains of definition.
Since $\cV_{\phi_0}$ is dense in $\Sigma\times M$, this defines
$\phi$ in $\Sigma \times M.\;$ The final step is to show that
$\phi$ is uniformly $C^{1,1}$ and solve (\ref{eq:hcma0}). Note
that this approach is different from \cite{Dona01} since we are
not dealing with  a super regular foliation.


\begin{prop}  \label{defo:localpotentialforsuperregular}For any super regular leaf $f, $  there exists a smooth function
$\phi$ defined in a small tubular neighborhood ${\cal O}_f$ (which
is saturated with respect to $\cF_{\phi_0}$) of $\pi\circ ev
(f)\subset \Sigma \times M$ such that
\begin{eqnarray}
\Omega  &=  \pi_2^*\omega_0 + i \dbd \phi,\qquad &{\rm on}
\;\; {\cal O}_f\subset \Sigma \times M \label{eq:foli3}\\
\;\;\phi  &=   \phi_0, \qquad  \qquad ~~~&{\rm on}\;\; {\cal O}_f
\bigcap \left(\p \Sigma \times M\right). \label{eq:foli4}
\end{eqnarray}
\end{prop}

\begin{rem} \label{defo:localpotentialfor  super regular1}
The potential function $\phi$ can be also defined in a tubular
neighborhood $\cO_f$ of an almost super regular leaf $f$. The
function is smooth in ${\cal O}_f \bigcap (\Sigma^0\times M)$.
\end{rem}
\begin{proof}
Recall that $M = \displaystyle\; \bigcup_{i\in {\cal I}}\; U_i$ and
 $\rho_i (i \in \cI) $ is the local defining potential function
 for $\omega_0.\;$ For any point $(z,x) \in \cV_{\phi_0}, $ suppose that $x\in U_i$ for some $i \in \cal I.\;$
  In local coordinates, write
\[
\omega_0 = \displaystyle \sum_{\alpha, \beta=1}^n\; {{\p^2 \rho_i}
\over {\p  w_\alpha  \p w_{\b \beta}}} \; d\, w^\alpha \wedge
\;d\, w^{\b \beta}, \qquad
  \tilde \omega = \displaystyle \sum_{\alpha, \beta=1}^n\; {{\p^2 (\rho_i+\phi)}\over {\p w_\alpha \p \b w_\beta}}
  \; d\, w^\alpha \wedge \;d\, w^{\b \beta}.
\]

We can express the image of $\Sigma$ as\footnote{Here
$x=\pi_2\circ \pi \circ f(z)$  in the formula. However, for
notation simpliticity, we simplify it as $x=f(z).\;$ This
convention will be used later.  This should cause no confusion.}
\[
 \{\left(z, x = f(z), \xi(x)\right), \;\;\forall z\in \Sigma\},
\]
where
\[
  \xi(x) = \p (\phi +\rho_i).
\]
Since $\xi(x)$ is uniquely detrmined by image of an open set of
super regular disks in $\Sigma\times \cW_M, $ then $\phi$ is
uniquely determined by $\xi$, or by the structure of $\cW_M,\;$ up
to a constant in $U_i\subset M.\;$ In particular, in $\left(\p
\Sigma \times M\right) \bigcap
\cV_{\phi_0},\;$ we have $\phi = \phi_0\;$ modular some function in $z$ locally.  \\

By definition,   we write
 \[
 \begin{array}{l} \Omega =
   \omega_0 + \displaystyle
\sum_{\alpha,\beta = 1}^n\; {{\p^2 \phi}\over {\p w_\alpha \p \b
w_\beta}} d\,w_\alpha \;d\,\b w_\beta +\displaystyle
\sum_{\alpha=1}^n\; \zeta^\alpha\;  d\,w_\alpha \;d\,\b z \\
\qquad \qquad + \sum_{\beta=1}^n\; \zeta^{\b \beta}\; d\,z
\;d\,w_{\b \beta} + h_f d\,z \;d\,\b z.
\end{array}
 \]
Since
\[i_{{\p \over {\p z}} + X} \Omega =0~~{\rm and}~~X= \sum_1^n\; \eta^\alpha \;
{\p \over {\p w^\alpha}},\] we have
\[
\zeta^\alpha + \left(g_{0,\alpha \b \beta} + {{\p^2 \phi}\over {\p
w^\alpha \p w^{\b \beta}}}\right)  \eta^{\b \beta} = 0, \qquad
\forall\; \alpha =1, 2, \cdots n.
\]
On the other hand,  since $\xi(f(z))$ is a holomorphic function of
$z$, we have
\[
  \begin{array}{lcl} 0 & = & {{\p \xi^\alpha}\over {\p \b z}} \\
  & = &  {{\p^2 \phi}\over {\p w^\alpha \p \b z}}  + {{\p^2 (\phi +\rho_i)}\over {\p w^\alpha \p w^{\b \beta}}}  {{\p f^{\b \beta}}\over {\p \b z}} \\
  & = &  {{\p^2 \phi}\over {\p w^\alpha \p \b z}}  + \left(g_{0,\alpha \b \beta} + {{\p^2 \phi}\over {\p w^\alpha \p w^{\b \beta}}}\right)  \eta^{\b \beta}.  \end{array}
\]
Then
\[ \zeta^\alpha = {{\p^2 (\phi + \rho_i)}\over {\p w^\alpha \p \b z}} , \qquad \forall \alpha=1,2,\cdots n.\]
Consequently, $\Omega$ takes the form:
 \[
 \begin{array}{l}
   \omega_0 + \displaystyle
\sum_{\alpha,\beta = 1}^n\; {{\p^2 \phi}\over {\p w_\alpha \p \b
w_\beta}} d\,w_\alpha \;d\,\b w_\beta +\displaystyle
\sum_{\alpha=1}^n\; {{\p^2 \phi}\over {\p \b z \p
   w_\alpha}} d\,w_\alpha \;d\,\b z \\ \qquad \qquad + \sum_{\beta=1}^n\;
{{\p^2 \phi}\over {\p z \p
   w_{\b \beta}}} d\,z \;d\,w_{\b \beta} + h_f d\,z \;d\,\b z.
\end{array}
 \]
Moreover, the following hold about the leaf vector $X:$
\[
{{\p^2 \phi }\over {\p w^\alpha \p \b z}} = - \left(g_{0,\alpha \b
\beta} + {{\p^2 \phi}\over {\p w^\alpha \p w^{\b \beta}}}\right)
\eta^{\b \beta}, \qquad \forall \alpha = 1,2, \cdots n.
\]
Using this equation and the fact that $\;\Omega^{n+1} = 0, $ we
have
\[\begin{array}{lcl} 0 & = & \left(h_f - {g_{\varphi}}^{,\alpha \b \beta} \cdot {{\p^2 \phi }\over {\p w^\alpha \p \b
z}} \cdot {{\p^2 \phi}\over {\p w^{\b \beta} \p z}} \right)\cdot \tilde{\omega}^n\\
&= & \left(h_f - \left(g_{0,\alpha \b \beta} + {{\p^2 \phi}\over
{\p w^\alpha \p w^{\b \beta}}}\right)  \eta^{\b \beta} \eta^\alpha
\right) \cdot \tilde{\omega}^n. \end{array}
\]
Since $\tilde{\omega} $ is non-degenerate in $\cO_f$, then
\[h_f =  \left(g_{0,\alpha \b \beta} + {{\p^2 \phi}\over
{\p w^\alpha \p w^{\b \beta}}}\right)  \eta^{\b \beta} \eta^\alpha
\]
is uniquely determined as well.\\

 However (note $z= w_{n+1}$),
\[
  \Omega - \displaystyle \sum_{\alpha, \beta=1}^{n+1} \; \sqrt{-1} \;{{\p^2 (\rho_i + \phi)}\over {\p w_\alpha \p
   w_{\b \beta}}} d\,w_\alpha \;d\,w_{\b \beta}
\]
is a closed form. Consequently,
\[
  l_f d\,z \;d\,\b z = \left(h_f - {{\p^2 \phi}\over {\p z \p \b z}}
  \right) d\,z \;d\,\b z
\]
is a closed form on $\cO_f.\;$ Therefore, $l_f$ is a function of
$z$ only.  Locally, we can replace $\phi$  by $\phi+ K_f(z)$ where
\[
{{\p^2 K_f}\over {\p \;z\; \p \b z\;}} = - l_f.
\]

After such a replacement, in each local coordinate chart in
$\cO_f$, we can choose the potential function $\phi$ uniquely, up
to a harmonic function on $z$ only. This follows from the fact
that
\[
 {{\p^2 \phi}\over {\p z \p \b z}} = h_f
\]
 where $h_f$ is uniquely determined by geometric data of super regular disks in $\cO_f$.\\

Because of the unique extension property of harmonic functions,
$\phi$ is uniquely determined in $\cO_f$ by a global harmonic
function in the $z$ direction. Choose such a potential function in
$\cO_f $ now.  Notice that in $\left(\p \Sigma\times M\right)
\bigcap \cO_f$, we have $\omega_z\mid_{\p \Sigma} =
\omega_{\phi_0}.\;$ Then, we can set
\[
  \phi (z,\cdot) = {\phi_0} (z,\cdot) + m_f(z), \qquad \forall\; z\in \p \Sigma
\]
where $m_f(z)$ is a function of $z$ in $\Sigma.\;$ Note that $m_f$
may not be harmonic function.

 Choose a function $K_{ff} $
as a function of $\;z\;$ only such that \[ {{\p^2 K_{ff}}\over {\p
\;z\; \p \b z\;}} = 0.
 \]
 and
 \[
 K_{ff} \mid_{\p \Sigma} = m_f.
 \]
 This Drichelet problem has a unique solution. Now
 replace $\phi$ by $\phi - K_f.$ Then
 $\Omega $ can be re-written as $\pi_2^*\omega_0 + \sqrt{-1}
\dbd \phi$ in a tubular neighborhood of $\pi\circ f(\Sigma) $ in
$\Sigma \times M$ such that $\phi =\phi_0$ in
$\left(\p\Sigma\times M\right) \bigcap \cO_f.\;$
\end{proof}

Now $\phi$ satisfies equation (\ref{eq:hcma0}) on
$\cV_{\phi_0}.\;$ Now we wish to extend it to solve the same
equation in $\Sigma\times M.\;$ The key step  is to prove
a global $C^{1,1}$ bound for $\phi$ in $\cV_{\phi_0}.\;$\\

\begin{prop} \label{defo:omegaispositive} We have $\omega_\phi > 0$ when restricted to $M$  in $\cV_{\phi_0}$.
 \end{prop}
 \begin{proof} For any point $(z,x') \in \cV_{\phi_0},$  there exists
 a holomorphic leaf $f \in \cU_{\phi_0}$  such that $\pi\circ ev(z,
f) = (z, x').\; $  For any $z \in
\partial \Sigma,\;$
\[
\omega_\phi =  \omega_0 + \sqrt{-1} \p \b \p \phi_0(z, \cdot)
> 0.
\]
In other words, $\omega_\phi\mid_{\pi\circ ev(\p \Sigma, f)} >
0.\; $ However, $\omega_\phi\mid_{\pi\circ ev( \Sigma, f)}$ is a
constant form. Thus, $\omega_\phi$ is strictly positive for any
$(z,x')\in \cV_{\phi_0}.\;$
 \end{proof}
It follows from this proposition that $\omega_\phi$ defines a
smooth K\"ahler metric in $\cV_{\phi_0}$. We want to show that
this metric has a uniform $L^\infty$ bound in $\cV_{\phi_0}.\;$ In
a local coordinate chart, write
\[\left\{\begin{array}{lcl}
\omega_0 &=  & \displaystyle\sum_{\alpha,\beta=1}^n\;
g_{0,\alpha\b\beta}\; d\,w^\alpha\;d\,w^{\b \beta},\qquad
\omega_\phi = \displaystyle\sum_{\alpha,\beta=1}^n\;
g_{\phi,\alpha\b\beta}\; d\,w^\alpha\;d\,w^{\b \beta},\\
g_{\phi,\alpha \b \beta} & = & g_{0,\alpha\b\beta} + {{\p^2
\phi}\over {\p w^\alpha \p w^{\b \beta}}},\;\forall\;\alpha,
\beta=1,2 \cdots n.
\end{array} \right.
\]
For any super regular leaf $f$ and for any $z\in \Sigma,\;$ the
restricted bundle $T_{f(z)}M$ at $f(z)$ is a trivial holomorphic
bundle over $\Sigma$ with complex rank $n$. Restriction of
$g_\phi$ to this $TM$ bundle induces an hermitian metric on this
bundle. Denote by ${F^{\alpha}}_{ \beta} (1\leq \alpha,\beta\leq
n)$ the curvature of this hermitian metric. We have the following
formula (cf. Section 4),
\begin{eqnarray}
{F^{r }}_{ \alpha} & = & -  \p_{\b z} \left(g_{\phi}^{ r \overline{\delta}} \;\partial_{z} g_{\phi,\alpha
\overline{\delta}}\right) \\
& = & - {{\partial \eta^{r}}\over {\partial w_{\b i}}} \cdot
{{\partial \eta^{\b i}}\over {\partial w_\alpha}} \leq 0.
\label{eq:curv} \end{eqnarray} In particular, the curvature is
always semi-negative.
Now we are ready to state the maximum principle.

\begin{prop} (Maximum Principle along leaves) The K\"ahler metric $g_\phi$ is uniformly bounded from
above in each leaf in $\cU_{\phi_0}.\;$
\end{prop}
\begin{proof} First give a finite covering of $M = \bigcup_{i\in \cI} U_i..\;$ There exists a uniform constant $C(g_0)$ which depends
only on this covering such that in each of these coordinate chart,
we have
\[
\|{{\p g_{0,i\b j}}\over {\p w_k}}\|_{g_0} \leq C(g_0), \qquad
\forall i,j,k = 1,2,\cdots n.
\]
For any $f\in \cU_{\phi_0},$ the restricted $TM$ bundle on $\pi\circ
f(\Sigma)$ is a trivial holomorphic vector bundle over $\Sigma.\;$
For any $(z,p)$ in this leaf, pick any $n-$frame $s_1,s_2,\cdots,
s_n \in T_p^{1,0} M$ and extend these vectors over the disk as a
frame of $n$ constant holomorphic sections in this $T^{1,0}M$
bundle.  We still denote as $\{s_1,s_2,\cdots, s_n\}.\;$ Note that
in any coordinate charte $U_i$ chosen above, we have
\[{\p \over {\p z}} g_0(s,s) =
{\p \over {\p z}} g_{0,i\b j} s^i s^{\b j} = {{\p g_{0,i\b j}}\over
{\p w_k}} X^k s^i s^{\b j} \leq C(g_0) \displaystyle \max |X| \cdot
g_0(s,s).
\]
In other words, for these constant holomorphic section of the $TM$
bundle over $f(\Sigma)$, it is easy to see that
\begin{equation}
 \displaystyle \inf_{z\in \Sigma} \; g_0(s, s) \geq c_0 (f)  \displaystyle \max_{z\in
   \Sigma} \; g_0(s, s)\label{def:g0control}
   \end{equation}
where $c_0(f)$ depends on $g_0$ and the embedding of $f$ only. \\

To prove the maximum principle for metric $g_\phi$ along the leaf,
we just need to show that $g_\phi(s, s) (1\leq \alpha \leq n) $ has
a uniform upper bound for any constant holomorphic leaf the tangent
bundle. If the upper bound is achieved on the boundary, then the
claim is proved since $g_\phi= g_{\phi_0}$ in $\p\Sigma\times M.\;$
If the maximum is attained at some interior point $(z,p)\in \Sigma^0
\times M,\;$ choose an appropriate coordinate in $T_pM.\;$ We may
assume
\[
 g_{\phi,\alpha \b \beta} (z,p) =
 \delta_{\alpha \beta} (\forall \;\alpha,\beta=1,2,\cdots n),
 \qquad \p_z g_\phi \mid_{(z,p)} = \b \p_z g_\phi \mid_{(z,p)} =
  0.
\]

At this point $(z,p),$ we have
\[\begin{array}{lcl}
  \b \p_z \p_z g_\phi(s(c), s(c))& = & g_\phi(\p_z
  s(c), \p_z s(c)) + (\b \p_z \p_z g_\phi)(s(c),
  s(c))\\
  & = & g_\phi(\p_z
  s(c), \p_z s(c)) -F (s(c),
  s(c))\geq 0.\end{array}
\]
This shows that the maximum must achieve in the boundary.   This
along with inequality \ref{def:g0control} gives desired upper bound
on $g_\varphi.\;$
\end{proof}

Using this proposition, we can prove
\begin{theo} There exists a uniform upper bound for the K\"ahler
metric $g_\phi$ in $\cV_{\phi_0}.$
\end{theo}
\begin{proof} This follows from the proof of the previous theorem as
well as the fact that any sequence of holomorphic discs in
$\cV_{\phi_0}$ always have a subsequence which converges to an
embedded holomorphic discs in $\Sigma\times M.\;$
\end{proof}
Finally, we have
\begin{theo}  
$\phi$ is uniformly $C^{1,1}$ in $\Sigma\times M$, smooth in
$\cV_{\phi_0}\; $ such that it solves

\begin{equation}
\left(\pi_2^*\omega_0 + \sqrt{-1}\p \b \p \phi\right)^{n+1} = 0,
\qquad {\rm in}\qquad \Sigma\times M. \label{defoeq:folid4}
\end{equation}
Moreover, this solution coincides with the solution established in
\cite{chen991}.
\end{theo}
\begin{proof} By Definition
\ref{def:nearlyfoliation}, $\phi$ can be extended to be a $C^1$
continuous function in $\Sigma\times M \setminus \tilde
\cS_{\phi_0}$. Since
$\tilde \cS_\phi$ is locally extendable and $\phi$ is uniformly
$C^{1,1}$ bounded in $\cV_{\phi_0}$, we can show that $\phi$ can be
extended as a global $C^{1,1}$ function in $\Sigma\times M$.    It
follows that there is a sequence of K\"ahler potential $\{\phi_m \in
\cH, m \in \NN\}$ such that $\phi_m \rightarrow \phi$ in weak
$C^{1,1} (\Sigma\times M)$ topology.  Moreover, the convergence is
smooth in any compact subset of $ \cV_{\phi_0}.\;$  It follows that
for any smooth test function $\psi,$ we have
\[
\displaystyle \lim_{m\rightarrow \infty} \int_{\Sigma\times M} \; \psi \omega_{\phi_m}^{n+1} = 0.
\]
Suppose $\varphi$ is the $C^{1,1}$ solution given by the first author in \cite{chen991}, we have

\[
\begin{array}{lcl}  0 & = & \displaystyle \lim_{m\rightarrow \infty} \int_{\Sigma\times M} \; (\varphi-\phi) \left(\omega_{\varphi}^{n+1}  -   \omega_{\phi_m}^{n+1}\right) \\
& =  & \displaystyle \lim_{m\rightarrow \infty} \int_{\Sigma\times M} (\varphi- \phi) (\omega_\varphi-  \omega_{\phi_m})\left(\displaystyle \sum_{i=0}^n \omega_\varphi^i  \wedge \omega_{\phi_m}^{n-i} \right)
\\ & = & - \displaystyle \lim_{m\rightarrow \infty} \int_{\Sigma\times M} \;\sqrt{-1} \p (\varphi- \phi) \wedge \b \p (\varphi- \phi_m)\left(\displaystyle \sum_{i=0}^n \omega_\varphi^i  \wedge \omega_{\phi_m}^{n-i} \right)\\
& = & - \int_{\cV_{\phi_0}} \;\sqrt{-1} \p (\varphi- \phi) \wedge \b \p (\varphi- \phi)\left(\displaystyle \sum_{i=0}^n \omega_\varphi^i  \wedge \omega_{\phi}^{n-i} \right).
\end{array}
\]
It follows that $\phi = \varphi$ in $\cV_{\phi_0}$ since $\omega_\phi$ is smooth in $
\cV_{\phi_0}.\;$  Since $\cV_{\phi_0}$ is dense in $\Sigma\times M$, it follows that the solution
we constructed conincide with solution established in the first author's paper \cite{chen991}.

\end{proof}

\section{Deformation of holomorphic disks with totally real boundary}

\subsection{Local analysis of holomorphic disks}

For any boundary map ${\phi_0}:\partial \Sigma \rightarrow \cal
H$, there is a $2n+1$-dimensional totally real submanifold $\b
\Lambda_{\phi_0} = \displaystyle \bigcup_{z \in \p \Sigma}\;
\left(\{z\}\times \Lambda_{{\phi_0}(z,\cdot)}\right)$ in $\Sigma
\times \cW_M$. Consider the {\it moduli} space ${\cal M}_{\phi_0}$
of all of holomorphic disks $ \rho: \;(\Sigma, \p \Sigma)
\rightarrow \left(\Sigma \times \cW_M, \b\Lambda_{\phi_0}\right)$
with vanishing normal {\it Maslov} index. Note that the normal
bundle over $\rho(\Sigma)$ in $\Sigma \times \cW_M\;$ is always
trivial and we will denote it by
$$ \left(\begin{array}{c} \CC^{2n}\\
\downarrow \pi \\  \rho(D) \end{array}\right).\;$$ For any $z =
e^{i \theta}\; (0\leq \theta \leq 2\pi)$, let
$\RR^{2n}(e^{i\theta})$ be the totally real subspace  $\rho^* \;
T_{\rho(e^{i\theta})}\; \left(\b \Lambda_{{\phi_0}}\right)$ of $
T_{\rho(e^{i\theta})} \CC^{2n} = \CC^{2n}$. Consider all
$H^{1,2}$-sections $s: \Sigma\rightarrow \CC^{2n}$ such that
$s(e^{i\theta}) \in \RR^{2n}(e^{i \theta}).\;$  The linearized
operator of $\rho$ is given by
 \[
 \bar \p_z: H^{1,2}(\Sigma, \CC^{2n})\rightarrow L^{2}(\Sigma, \CC^{2n}).
 \]
This is a Fredholm operator, so we can compute its index
\[
{\rm index}(\bar \p_z) = {\rm dim\; Ker}(\bar \p_z) - {\rm dim\;
Coker}(\bar \p_z).
\]
This indice is invariant under deformation of holomorphic disks.
If the normal {\it Maslov} indice of $\rho$ is denoted by
$\mu(\b\p_z)$, then we have the following (cf. \cite{ve67}):
\[
  {\rm indice}(\bar \p_z) = \mu(\b\p_z) + 2n = 2n.
\]
Thus the kernel of $\bar \p_z$ is of dimension at least $2n$.
Recall that a holomorphic disk $\rho$ is regular in the sense of
the Fredholm theory if the cokernel of $\bar\p_z$ vanishes,
equivalently, the kernel has dimension exactly equal to $2n$.

For every disk $\rho: (\Sigma,\p \Sigma)\rightarrow (\Sigma\times
\cW_M, \b \Lambda_{\phi_0})$, we have a loop of $2n$-dimensional
real subspaces $\{\RR^{2n}(e^{i\theta})\,|\,0\leq \theta\leq
2\pi\}$ in $\CC^{2n}$. By fixing a real $\RR^{2n}$ subspace in
$T_{\rho(1)} \; \cW_M$, this induces a map from $\p \Sigma$ to
$GL(2n, \CC)/GL(2n, \RR)$. In general, this map may not lift to a
map from $\p \Sigma$ to $Gl(2n, {\CC})$. However, this property of
being able to be lifted to $C^\infty(\p \Sigma, GL(2n, {\CC}))$ is
invariant under continuous deformation of the disk or the boundary
conditions. A disk is called {\it trivial} if all real subspaces
$\RR^{2n}(e^{i \theta})$ are equal to a constant $2n$-dimensional
real subspace (independent of $e^{i\theta}$). For a trivial disk,
its induced loop always admits a lifting to $GL(2n, {\CC}).\;$
Therefore, if a disk is path-connected to a trivial disk, then its
induced loop must admit a lifting to an map $C^\infty(\p
\Sigma,GL(2n, {\CC})).\;$ We call this an {\it associated loop} of
the disk $\rho.\;$  It is clear that {\it associated loop} is
defined up to multiplication by ${\cal L}^+ GL(2n,\CC)$ on the
left. Here ${\cal L} GL(2n,\CC) = C^\infty(\p \Sigma, GL(2n,
{\CC})),$
 and ${\cal L}^+ GL(2n,\CC) \subset C^\infty(\p \Sigma, GL(2n, {\CC}))$  is the set of loops which can
 be extended to a holomorphic map $ C^\infty(\Sigma, GL(2n, {\CC})).\;$  In this paper, we only consider discs which are path
connected to a trivial disc. For these holomorphic discs, it is
natural to consider the partial indices which are independent of
the lifting. According to \cite{ve67}, \cite{Gl94} and
\cite{Oh952}, using a special form of Birkhoff factorization, we
have
\\
\noindent {\bf Theorem A} \footnote{For the generic {\it Maslov }
indice, this theorem was first obtained by \cite{Franc87} in
complex surface. It was generalized to all dimensions in
\cite{Gl94} with the assumption that all of the partial indices
are non-negative. This last restriction was removed in
\cite{Oh952}. The present statement follows the format in
\cite{Oh952}. }: {\it Let $\tilde{\rho}: S^1 = \p \Sigma
\rightarrow \RR^{2n}(\theta)$ be a loop of totally real $2n$
dimensional sub-spaces in $\CC^{2n}, \;$ which is induced by some
holomorphic disk $\rho: \;(\Sigma, \p \Sigma) \rightarrow
\left(\Sigma \times \cW_M, \b \Lambda_{{\phi_0}}\right).\;$ Then
one can represent this loop map as
\[ \tilde{\rho}(z) = \Theta(z) \Lambda(z)^{1\over 2} \cdot
\RR^{2n},\qquad z \in \p \Sigma = S^1,
 \]
where $\Theta(z) \in {\cal L} GL(2n,\CC)$ and $\Lambda(z)$ is a
diagonal matrix:
\[
  \Lambda(z) = [ z^{k_1}, z^{k_2},\cdots, z^{k_{2n}}], \forall \;z \in \p \Sigma.
\]
Here $(k_1,k_2, \cdots ,k_{2n})$ is called partial indices of the
loop $\rho.\;$ Moreover, these partial indices have the following
properties:
 \begin{enumerate}
 \item Each individual partial indice is not invariant under
 continuous deformation. However, the total sum of all partial indices is precisely
 the normal Maslov invariant, so it is invariant under any continuous defomation. Thus,
 \[
 \displaystyle \; \sum_{i=1}^{2n}\; k_i = \mu = 0.
 \]
 \item A disk is Fredholm regular if and only if all of its partial
 indices $\geq -1.\;$
 \end{enumerate}}

 Using this theorem, Oh was able to reduce the equation for kernel
 vectors into a scalar equation:
 \[ u = \left\{
\begin{array}{lcll} {{\p \xi}\over {\p \bar z}} &  =  & 0, & \qquad
\forall z\in \Sigma,
\\ \xi(z)  & =  & z^{k_i\over 2} \cdot {\RR}, & \qquad\forall z
\in S^1.
\end{array}\right.
\]

 This equation has no solution when $k_i \leq -1.\;$ For $k_i\geq
 0,\;$ this equation has exactly $k_i + 1$ linearly independent
 solution while each solution is a polynomial in $z$ with degree
 $k_i.\;$\\

\noindent {\bf Theorem B} {\it Suppose $f$ is a regular disk whose
partial indice decomposition $(k_1, k_2, \cdots ,k_{2n})$ contains
exactly $l ( \in [0,n])$ number of partial indices which equals $-1.\;$
Then the kernel matrix of this disk has co-rank at least
$l$ everywhere in the interior of this disk.}\\

This can be easily derived from \cite{Gl94} and \cite{Oh952}.

\subsection{The universal moduli space is regular}

Define
\[
{\cal G} = \displaystyle \bigcup_{{\phi_0} \in C^\infty(\p \Sigma,
{\cal H})}\,\b \Lambda_{\phi_0},
\]
and
\[ \Upsilon =\displaystyle \bigcup_{{\phi_0} \in C^\infty(\p
\Sigma, {\cal H})}\,{\cal M}_{\phi_0},
\]
where $ {\cal M}_{\phi_0} $ is the moduli space of all holomorphic
disks with vanishing normal {\it Maslov} indice:
\[
\rho: \left(\Sigma, \p \Sigma \right) \rightarrow \left(\Sigma
\times \cW_M, \b \Lambda_{{\phi_0}}\right).
\]
Clearly, ${\cal G}$ is an infinite dimensional manifold and there
is a natural projection $p: \Upsilon \rightarrow {\cal G}$ such
that for any ${\phi_0} \in C^\infty(\p \Sigma, {\cal H}),$ the
moduli space $\cM_{\phi_0}$ is mapped to $\b \Lambda_{\phi_0}$.

Recall that the {\it Moduli} space is {\bf smooth} if every
holomorphic disk in this {\it Moduli} space is regular. It follows
from the following lemma and the Sard-Smale transversality theorem
that $\cM_{\phi_0}$ is smooth for a generic boundary value
$\phi_0$.

\begin{lem} \footnote{This was first
carried out in \cite{Oh96} in the context of Lagrange/totally real
submanifold. For convenience of readers, we include a proof of
this transversality below.} The universal {\it moduli} space  $p:
\Upsilon \rightarrow {\cal G}$  is smooth.
\end{lem}
\begin{proof} The tangent space of ${\cal G}$ at $\phi_0$ can be considered as
\[
  T_{\phi_0}  {\cal G } = C^\infty( \p \Sigma \times M), \qquad
  \forall \; {\phi_0} \in C^\infty(\p \Sigma, {\cal H}).
\]
Let $\epsilon_k\rightarrow 0$ be a sequence of positive numbers
which converges to zero.  Denote $\b \epsilon =
(\epsilon_1,\epsilon_2,\cdots ).\;$ Set
\[
   \|f\|_{\bar \epsilon} = \displaystyle\; \sum_{k=0}^\infty \epsilon_k\; \displaystyle
   \max_{x \in \p\Sigma\times M}\;
   |D^k\; f(x)|.
\]
This defines an $\b \epsilon$-norm on \[ C^{\b
\epsilon}(\Lambda_{\phi_0}) = \{ f \in C^{\infty}(\p \Sigma \times
M) \mid \|f\|_{\b \epsilon} < \infty\}.
\]
This norm has been introduced by Floer in a different context and
this is a Banach space. We can choose $\b \epsilon$ so that $C^{\b
\epsilon}(\Lambda_{\phi_0})$ is dense in $C^\infty(\p \Sigma
\times M)$ with respect to the $L^2 $ norm.

Now fix $s> 1$ and define \[ {\cal F} = {\cal F}^s =
H^{s+1}(\Sigma, \Sigma \times \cW_M)\] which is a Sobolev space of
all maps $\omega: \Sigma \rightarrow \Sigma \times \cW_M$ whose
$(s+1)^{th}$ derivatives are in $L^2.\;$ For any boundary map
${\phi_0}: \p \Sigma \rightarrow \cal H$, the totally real
submanifold $\Lambda_{\phi_0} $ of $\Sigma \times \cW_M $ is a
point in $\cal G.\;$ For any small $r$ positive, we define an $r$
neighborhood of this point in $\cal G$ as:
$${\cal N}(\Lambda_{\phi_0}) =  \{ \Lambda_{{\phi_0} + f}\mid
\|f\|_{\bar \epsilon} < r\;{\rm and} \; f\in C^\infty(\p \Sigma
\times M)\}.\;$$
 The corresponding neighborhood of holomorphic disks is
\[\begin{array}{lcl}
 \bar {\cal M} & = & {\cal M}({\cal N}(\Lambda_{\phi_0}))\\
  & =&\{ (\rho, \Lambda_{{\phi_0} + f} ) \,\mid \,\bar \p \rho = 0,\; \rho\mid_{\p \Sigma} \subset \Lambda_{{\phi_0} + f}, \;
  \|f\|_{\bar \epsilon} < r\}.
\end{array}\]

For each $\rho \in {\cal F}$, define the pulled back bundle as
\[ {\cal B}_\rho = H^s(\rho^* T \cW_M)
\]
consisting of all $H^s$-sections of $\rho^* T\cW_M$ on $\Sigma$.
Set
\[ {\cal B} = \displaystyle \bigcup_{\rho \in {\cal F}} \; {\cal
B}_\rho =  \displaystyle \bigcup_{\rho \in {\cal F}} \; H^s(\rho^*
T \cW_M).
\]
This is a smooth bundle over ${\cal F}$. We further set
\[ \Omega(\Lambda_{\phi_0})
:= H^{s+{1\over 2}} (\p \Sigma, \Sigma \times \cW_M) \cap C^0(\p
\Sigma, \Lambda_{\phi_0}).
\]
This is simply the space of $H^{s+{1\over 2}}$ maps from $\p
\Sigma$ to $ \Lambda_{\phi_0}.\;$ By the trace theorem, for each
map $\rho \in H^{s+1} (\Sigma, \Sigma \times \cW_M), $ its
boundary map $\rho\mid_{\p \Sigma}$ lies in $\Omega(\Sigma\times
\cW_M) = H^{s+{1\over 2}}(\p \Sigma, \Sigma\times \cW_M)$. Now we
define a map
\[ \triangle: {\cal F}\times {\cal N}({\Lambda_{\phi_0}}) \rightarrow
{\cal B}  \times \Omega(\Lambda_{\phi_0})
\]
by
\[
  \triangle (\rho, \Lambda_{{\phi_0}+f}) = \left(\bar \p \rho, \phi_{{\phi_0}+f}^{-1} \left(\rho\mid_{\p
  \Sigma}\right)\right)
\]
where $\phi_{{\phi_0}+f}: \Lambda_{{\phi_0}}\rightarrow
\Lambda_{{\phi_0}+f}$ identifies the small perturbation
$\Lambda_{{\phi_0}+f}$ with $\Lambda_{\phi_0}$.  Denote by
\[ X_f
= {d\over {d\,t}}\left (\phi^{-1}_{{\phi_0}+ t
f}\right)\mid_{t=0}\, \in\, T_{\Lambda_{\phi_0}} \;{\cal
N}({\Lambda_{\phi_0}}).
\]
Consequently, $T_{\Lambda_{\phi_0}} \;{\cal
N}({\Lambda_{\phi_0}})$ consists of all such fields $X_f$ for
$f\in C^\infty(\p\Sigma\times M)$. Note that
\[
\bar {\cal M} =  \triangle^{-1} (\{0\}\times
\Omega(\Lambda_{\phi_0})).
\]

The goal here is to show that the map $\triangle$ is transverse to
the submanifold at $(\{0\}\times \Omega(\Lambda_{\phi_0})) \subset
\;{\cal B}  \times \Omega(\Sigma \times \cW_M)$. Then it follows
that $\bar {\cal M}$ is a smooth Banach submanifold of ${\cal F}^s
\times {\cal N}(\Lambda_{{\phi_0}})$. Moreover, by the elliptic
regularity theory, $\bar {\cal M}$ is a smooth Banach submanifold
of ${\cal F}^s \times {\cal N}(\Lambda_{{\phi_0}})$ for all $s >
1$.

For any small $f\in C^\infty(\p \Sigma \times M)$, we set $\phi =
  \phi_{{\phi_0} +f}$ for simplicity. To show the transversality, we need to show
\begin{equation}
 Im_\triangle\;T_{\rho,\Lambda_{{\phi_0}+f}} \left({\cal F}\times {\cal N}({\Lambda_{\phi_0}})\right)+
\{0\} \oplus  T_{\phi^{-1}(\rho\mid_{\p\Sigma})}
\Omega(\Lambda_{\phi_0}) = T_{0,\phi^{-1}(\rho\mid_{\p \Sigma})}
\left({\cal B} \times \Omega(\p \Sigma \times
\cW_M)\right),\label{eq:transversality}
\end{equation}
where $(\rho, \Lambda_{{\phi_0}+f}) \in {\cal F}\times {\cal
N}({\Lambda_{\phi_0}}).\;$ If $(\xi, X_f) \in
T_{\rho,\Lambda_{{\phi_0}+f}} \left({\cal F}\times {\cal
N}({\Lambda_{\phi_0}})\right),$ then a straightforward calculation
shows that
\[
Im_\triangle\;(\xi,X_f) = \left(\bar \partial \xi, X_f -
\xi\mid_{\p\Sigma}\right).
\]
Clearly, the LHS (left hand side) of (\ref{eq:transversality}) is
a subspace of the RHS (right hand side). We need to show that the
normal space to LHS in the RHS of (\ref{eq:transversality}) is
null. Suppose that $(r,\alpha)$ is in such an normal space, that
is, $(r,\alpha) \perp Im_\triangle\;(\xi,X_f) $ and $(r,\alpha)
\perp \left( \{0\} \bigoplus T_{\phi^{-1}(\rho\mid_{\p\Sigma})}
\Omega(\Lambda_{\phi_0})\right)$. The second condition implies
\[\alpha \in
\left(T_{\phi^{-1}(\rho\mid_{\p\Sigma})}
\Omega(\Lambda_{\phi_0})\right)^{\perp}.\] In other words,
$\alpha$ represents some variation normal to
$T_{\phi^{-1}_{\phi_0}(\rho\mid_{\p \Sigma})}
\Omega(\Lambda_{\phi_0})$. The first condition implies that
\[
\displaystyle \int_\Sigma\;(\bar \p\;\xi,r) + \displaystyle
\int_{\p \Sigma}\;(X_f-\xi\mid_{\p\Sigma},\alpha) = 0.
\]
Integrating by parts, we have
\[
\displaystyle \int_\Sigma\;(\xi,\nabla_J\;r) +\displaystyle
\int_{\p \Sigma}\;(\xi\mid_{\p\Sigma}, e^{-i\theta} r -\alpha )
d\,\theta+ \displaystyle \int_{\p
\Sigma}\;(X_f-\xi\mid_{\p\Sigma},\alpha) d\,\theta = 0.
\]
Thus
\begin{eqnarray}
\nabla_J r & = & 0,\label{modu:paralell}\\
-\alpha\mid_{\p \Sigma} + e^{-i\theta} r\mid_{\p\Sigma} & = & 0,\label{modu:norminboundary}\\
\alpha^{\perp} & = & 0.\label{modu:norminterior}
\end{eqnarray}

Equation (\ref{modu:norminterior}) shows that $\alpha$ must be
tangent to $T_{\phi^{-1}(\rho\mid_{\p\Sigma})}
\Omega(\Lambda_{\phi_0})$, but $\alpha$ must be also normal to
this space, so $\alpha = 0$. Consequently, $r\mid_{\p \Sigma} = 0$
by (\ref{modu:norminboundary}). This, together with
(\ref{modu:paralell}), implies that $r=0$ in $ \Sigma$. This
completes the proof of transversality.
\end{proof}

The same arguments also show that for a generic path
$\psi:[0,1]\times \p\Sigma \mapsto \cH$, the total moduli
$\bigcup_{t\in [0,1]} \cM_|{\psi(t,\cdot)}$ is smooth.

\subsection{Selection of a generic path}

Next we turn our attention to variations of an arbitrary disk $f$
in the universal {\it moduli} space of holomorphic disks.  As
before, for every disk, it induces a map from $\p\Sigma $ to the
space of totally real $2n$ plane in $\CC^{2n}.\;$  Since all disks
concerned are path connected to a trivial disk, this induced map
can be lifted to a map from the universal {\it moduli} space of
holomorphic disks to the loop space ${\cal L} GL(2n, \CC).\;$ It
is well defined up to some normalization of the induced normal
bundle of $\rho^*T_w \cW_M$ over $\Sigma.\;$ In other words, it is
a map from a holomorphic disk to ${\cal L} GL(2n, \CC)/ {\cal L}^+
GL(2n,\CC)$. Define a fiber bundle ${\cal C}$ over ${\cal F}$ such
that each fibre is isomorphic to ${\cal L} GL(2n, \CC)/ {\cal L}^+
GL(2n,\CC)$. This defines a natural map from the universal {moduli
space} $\bar {\cal M}$ to this fibre bundle
\[
\star:{\cal G} \rightarrow {\cal C}
\]
which simply maps each holomorphic disk to its associated loop in
${\cal L} Gl(2n, \CC)/ {\cal L}^+ Gl(2n,\CC).\;$ \\

It is well known that ${\cal L} GL(2n, \CC)/ {\cal L}^+
GL(2n,\CC)\;$ admits a smooth stratification of loops by its
partial indice $ k = (k_1,k_2,\cdots, k_{2n}).\;$
A somewhat lengthy calculation \footnote{Proof of this Lemma can be founded in Appendix.} shows\\

\begin{lem} \label{defo:codimension} {\it For the smooth stratification of ${\cal L}
GL(2n, \CC)/ {\cal L}^+ GL(2n,\CC)$ by its partial indices
$k=(k_1,k_2,\cdots k_{2n})$, the real codimension of each
component indexed by $k$ is
\[
 d = \displaystyle \sum_{i=1}^{2n}\;\sum_{j=i+1}^{2n}\; \left(k_i- k_j - \lambda_{ij}\right)
\]
where $k_1\geq k_2\geq\cdots \geq k_{2n}$ and
\[\lambda_{ij} =\{ 
\begin{array}{ll} 1 & {\rm if} \;k_i> k_j\; {\rm and}\; i <j,\\
 0 & {\rm Otherwise.} \end{array}
\]}
\end{lem}
 Let
$S_0, S_1, S_2, S$ be the set of loops whose partial indices
satisfy:
\begin{enumerate}
\item All partial indices in $S_0$ are equal to  $0;\;$ \item all
partial indices are of the form $(1,0,\cdots, 0,-1)$ in $S_1$;
\item at least two of the partial indices equal to $-1$ in $S_2,$
but no partial indice is $\leq -2;$ \item at least one partial
indice in $S$ is less or equal to $-2.\;$
\end{enumerate}

According to Lemma \ref{defo:codimension}, $S_0$ is in generic
position, while the real codimension for $S_1$ is $1$. For $S
\subset {\cal L} GL(2n, \CC)/ {\cal L}^+ GL(2n,\CC)$, suppose that
$k_j \geq  k_{2n}$ and $k_{2n} \leq -2, \forall i \in [1, 2n].\;$
Then the codimension is:
\[\begin{array}{lcl} d & = & \displaystyle \sum_{i=1}^{2n}\;\displaystyle \sum_{j=i+1}^{2n}\; \left(k_i- k_j -
\lambda_{ij}\right)\\ & \geq & \displaystyle \sum_{i=1}^{2n}\;
\left(k_i- k_{2n} - \lambda_{i(2n)}\right) \\
 & = & \displaystyle \sum_{i=1}^{2n}\; k_i  + \displaystyle
\sum_{i=1}^{2n}\;(-k_{2n}) -\displaystyle
\sum_{i=1}^{2n}\;\lambda_{i(2n)} \\ &\geq & 0 +
\displaystyle\sum_{i=1}^{2n}\; 2 - \displaystyle \sum_{i=1}^{2n-1}
1 = 2n+1.
\end{array}
 \]
For $S_2,$ we can assume $k_{2n-1}=k_{2n}=-1.\;$ Thus,
\[\begin{array}{lcl} d & = & \displaystyle \sum_{i=1}^{2n}\;\displaystyle \sum_{j=i+1}^{2n}\; \left(k_i- k_j -
\lambda_{ij}\right)\\ & \geq & \displaystyle \sum_{i=1}^{2n-1}\;
\left(k_i- k_{2n-1} - \lambda_{i(2n-1)}\right) + \displaystyle
\sum_{i=1}^{2n-2}\;
\left(k_i- k_{2n} - \lambda_{i(2n)}\right)\\
 & = & 2 \displaystyle \sum_{i=1}^{2n-2}\; k_i  + \displaystyle
\sum_{i=1}^{2n-2}\;(-k_{2n} - k_{2n-1}) -\displaystyle
\sum_{i=1}^{2n-2}\;(\lambda_{i(2n)} + \lambda_{i (2n-1)}) \\
& = & 2 (k_{2n-1} + k_{2n})  + \displaystyle
\sum_{i=1}^{2n-2}\;(-k_{2n} - k_{2n-1}) -\displaystyle
\sum_{i=1}^{2n-2}\;(\lambda_{i(2n)} + \lambda_{i (2n-1)}) \\
&\geq & 4 + 2(2n-2) - 2 (2n-2) =4.
\end{array}
 \]

According to Lemma \ref{defo:codimension}, we have
\[
 {\cal L} GL(2n, \CC)/ {\cal L}^+ GL(2n,\CC) = S_0\bigcup S_1\bigcup
 S_2 \bigcup S = S_0 \bigcup (S_1^{a.s.} \bigcup S_1^{n.a.s.})
 \bigcup S_2 \bigcup S.
\]
Here $ S_1^{a.s.}$ denotes all of the holomorphic disks in $S_1$
which are super regular at $z=0$ and
\[
  S_1 = S_1^{a.s.}\bigcap S_1^{n.a.s}.
\]
It is straightforward to check that $S_1^{n.a.s.} $ has at least
codimension $1$ in $S_1.\;$\\

\begin{prop} This map $\star$ is a submersion at any embedded disk of
$\bar {\cal M}.\;$ \end{prop}
\begin{proof} We need to show that $\star$ is a submersion at a
regular disk or at a non-regular but embedded disk. The first
assertion follows from the fact that any regular holomorphic disk
$f$ with boundary in a totally real submanifold is stable under a
small deformation of the boundary map. To be more explicit, we let
$f: (\Sigma,\p\Sigma) \rightarrow \left(\Sigma\times \cW_M,
\Lambda_{\phi_0}\right)$ be a holomorphic disk with vanishing
normal {\it Maslov} indice.
If $f$ is regular in the sense of Fredholm theory, then there is a
holomorphic disk $f+\delta f$ such that its boundary lies in some
$ \RR^{2n}(e^{i\theta}) +\delta P(e^{i\theta})$. For this family
of holomorphic disks, the associated loop is exactly $\ell
+\delta\ell$. Thus, ``$\star$'' is an submersion at the image
(under $\star$) of every regular holomorphic disk.

Now suppose that $\star(f)=\ell \in S \subset \cC.\;$ Since $\cC$
is a smooth infinite dimensional manifold which admits a smooth
stratification by partial indices. More specifically, the space of
loop matrices in $Gl(2n, \RR),$  it can be decomposed as the union
of $ S_0 \bigcup S_1\bigcup S_2\bigcup S.\;$ Here we are only
considering $S$ which lie in a connected component of $S_1\bigcup
S_1\bigcup S_2.\;$ Therefore, there always exists a path
$\delta\ell (t) $ such that $\delta\ell(0)= 0;$  and $\ell
+\delta\ell(t) \in {\cal C}\setminus S$ when $t \neq 0.\;$
Consider
\[
 f:(\Sigma, \p \Sigma) \rightarrow (\Sigma\times \cW, \displaystyle \bigcup_{\theta\in S^1}\; (\theta, L(\theta))).
\]
Here $L(\theta)$ is a totally real sub-manifold in $(\theta, \cW)$
for any $\theta \in S^1$. At the tangential level, $T_{f(\theta)}
\cW_M$ is a trivial $\CC^{2n}$ bundle over $\Sigma$. Using this
trivialization, we may assume
\[
   T_{f(\theta)} L(\theta) = \RR^{2n}(\theta) =  A(\theta) \cdot
   \RR^{2n}
\]
for some $\RR^{2n} $ fixed in $T_{f(1)} \cW_M.\;$  Here $A(\theta)
\in Gl(2n, \CC).\;$  Clearly,  $\ell$ can be lifted up to be a
loop in $Gl(2n,\CC)$ and
\[
   \ell(\theta) =  A(\theta), \forall \; \theta \in S^1.
\]
The tangent space of $LGL(2n, \CC)$ at $\ell$ can be represented
by a smooth 1-parameter family of loops of matrices ($-\epsilon
\leq t \leq \epsilon$):
  \[
  \ell(t,\theta) =  A(\theta) (I + t B(\theta)), \qquad \forall \,\theta \in S^1.
  \]
The surjectivity at $f$ is equivalent to there being a pre-image
of this path $\ell(t,\theta)$ for an arbitrary loop matrix $B$.
Near a small tubular neighborhood of $L(\theta)\in \cW_\theta$, we
define a product metric (so $L(\theta)$ becomes totally geodesic
in $\cW_M$). Call this metric $g_\theta$. Define
\[
L(t,\theta) =  exp_{f(\theta), g_\theta} ((A (\theta) (I +t
B(\theta)) )\cdot \RR^{2n}), \] where $(A(\theta)(I + t
B(\theta)))\cdot \RR^n$ represents the n-dimensional plane spanned
by it in $T_{f(\theta)} \cW_M.\;$  Clearly,  $L(0, \theta) =
L(\theta).\; $ Define $f(t)$ to be a family of disks in the total
{\it moduli} space
\[
f(t): (\Sigma, \p \Sigma) \rightarrow (\Sigma\times \cW,
\displaystyle \bigcup_{\theta\in S^1}\; (\theta, L(t, \theta)).
\]
such that the image of each $f(t)$ is identified with $f$, but
they represent a 1-parameter path of holomorphic disks in the
total {\it moduli} space.  Clearly,
\[
   \star(f(t)) = \ell(t).
\]
In other words,
 the map $\star$ is transversal to $S \subset
{\cal L} GL(2n, \CC)/ {\cal L}^+ GL(2n,\CC).\;$
\end{proof}

Next we want to use this submersion map $\star$  to calculate the
codimension of various components of the universal {\it moduli}
space.

 Note that
$\star^{-1} S_0, \star^{-1}S_1, \star^{-1} S_2$ and $\star^{-1}S$
are smooth manifold or submanifold in $\bar {\cal M},$ where
$\star^{-1} S_0$ are the set of all super regular disks which is
generic in $\bar {\cal M}$, $\star^{-1}S_1$, $\star^{-1}S_2$ are
submanifolds of regular holomorphic disks in $\bar {\cal M}$ with
real codimension at least $1$ and $4$. Finally, $\star^{-1}S$ is
the smooth submanifold of all irregular disks in $\bar {\cal M}$
with real codimension at least $2n+1$. This, together with the
remark at the end of last section, implies

\begin{theo}\label{defo:genericpath1}
For any path $\psi: [0,1] \rightarrow C^\infty(\p \Sigma, {\cal
H}) $ such that ${\cal M}_{\psi(0,\cdot)}$ contains a super
regular disk with vanishing normal {\it Maslov} invariant, there
exists a generic path (we still denote it by $\psi$), which is
arbitrarily close to the original path, such that the total moduli
$\bigcup_{0\leq s \leq 1}\, \{s\}\times {\cal M}_{\psi(s,\cdot)}$
is a smooth $2n+1$-dimensional manifold.
Moreover, there is a connected component $\cM_\psi^0$ of this
total moduli containing the super-regular disk in initial moduli
$\cM_{\psi(0,\cdot)}$, such that the followings hold:
\begin{enumerate} \item The set of disks with partial indices
$(0,0,\cdots, 0)$ in $\cM_\psi^0$ is open and dense in this
connected component; \item The set of disks with partial indices
$(1,0,0,\cdots, 0, -1)$ has codimension at least 1 in
$\cM_\psi^0$. The set of disks with partial indice $(1,0,\cdots,
0,-1)$ but not super-regular at $z=0$ has codimension 2 and
higher;\item The set of all other holomorphic disks has
codimension 2 and higher; \item There exist at most finitely many
non-regular disks in the total moduli.

\end{enumerate}
\end{theo}

\subsection{Almost super regular foliations} In this subsection, we
first introduce the notion of an {\it almost super regular
foliation}. This is a stronger notion compared to {\it the nearly
smooth foliation} introduced earlier. Recall the natural
projection $\pi:\cW_{[\omega]}\rightarrow M$. The tangent space of
$\cW_{[\omega]}$ naturally splits into $TM$ union with tangent
space to the fibre space of $ T^*M.\;$ A regular disk $f$ in
${\cal M}_{\phi_0}$ is called {\bf super regular} at $z \in
\Sigma\times M$ if the Jacobi map of $ \pi \circ ev$ is
non-singular at $z \in\Sigma$, where $ev: \Sigma\times
\cW_{[\omega]}\mapsto \cW_{[\omega]}$ maps $(z,f)$ to $f(z)$. It
is called {\bf super regular} if it is super regular at every
point of $\Sigma$. It is called {\bf almost super regular} if the
Jacobi map of $\pi\circ ev$ is nowhere vanishing in $\Sigma^0.\;$


\begin{defi}\label{defo:almostsuperregularfoliation} For any boundary map ${\phi_0} \in C^\infty(\p \Sigma,
 \cal H), $  an open and connected $2n$ dimensional subset $\cU_{\phi_0} \subset {\cal
M}_{\phi_0}$ is called an {\bf almost super regular foliation} if
\begin{enumerate}
\item it is a nearly smooth foliation (cf. Definition
\ref{def:nearlyfoliation});
 \item  Every disk in $\bar \cU_{\phi_0}$ is regular except
 perhaps at most a set of finitely many disks. Furthermore, the set
 of almost super regular disks has co-dimension at least 1, while
 the set of all disks of other types has at least co-dimension 2 or higher.
\end{enumerate}
\end{defi}

Clearly, for any boundary map ${\phi_0} \in C^\infty(\p \Sigma,
 \cal H)$,  an almost super regular foliation $\cF_{\phi_0}$ is super
 regular if $ \cU_{\phi_0} = \bar \cU_{\phi_0} = \cM_{\phi_0}$.

\begin{prop}
If $\cF_{\phi_0}$ is an almost super regular foliation, then
$\Sigma\times \bar \cU_{\phi_0}$ induces a foliation in $\Sigma^0
\times M$ via the composition map of $\pi$ and $ev$, except at
most a set of codimension 2 .
\end{prop}
\begin{proof} Self evident.
\end{proof}

\begin{cor} For an almost super regular foliation, two disks
intersect at most at subset of of $\Sigma \times M$ with
codimension 2 or higher.
\end{cor}

This in turns implies
\begin{cor} For any almost smooth solution $\phi$ of (\ref{eq:hcma0})
which corresponds to an almost super regular foliation,  the leaf
vector field $X$ which annulate the {\it Levi} form
$\pi_2^*\omega_0 + \sqrt{-1} \p \b \p \phi$ is smooth in $\tilde
\cV_{\phi_0}$ and is uniformly bounded in $\Sigma^0\times M.\;$
\end{cor}

\begin{prop} \label{defo:almostsuperregular3}For a generic boundary map $\phi_0:\p \Sigma\rightarrow
\cal H$ such that every embedded disk in $\cM_{\phi_0}$ is
regular, then a nearly smooth foliation is necessary almost super
regular. Moreover, the connected component $\bar \cU_{\phi_0}$ is
a smooth manifold without boundary.
\end{prop}
\begin{proof} For any sequence of holomorphic disks  $f_k\in
\cU_{\phi_0},$ the leaf vector field $X_k$ has uniform upper
bound.  It follows that there is a subsequence (which we still
denoted as $\{f_k, k\in \cI\}),\;$ such that  converges to a
limiting embedded disk $f_\infty \in \cM_{\phi_0}.\;$ By our
assumption, this limiting disk must be regular in the sense of
Fredholm theory. In particular, $f_\infty$ is an interior point of
$\cM_{\phi_0}.\;$Consequently, $\b \cU_{\phi_0}$ is compact
without boundary.
\end{proof}
Conversely,
\begin{prop} \label{defo:almostsuperregular4}For a boundary map $\phi_0:\p \Sigma\rightarrow
\cal H$ such that all but possibly finitely many disks in
$\cM_{\phi_0}$ are regular. Define $\tilde \cU_{\phi_0}$to be the
set of all super regular and all almost super regular disks.
Suppose
\begin{enumerate}
\item$\bar \cU_{\phi_0}\setminus \tilde \cU_{\phi_0}$ has
codimension 2 or higher; \item The evaluation map is continuous on
$\bar \cU_{\phi_0}$.
\end{enumerate}
Then $\cU_{\phi_0}$ defines an almost super regular foliation.  In
particular, the covering indice for evaluation map is 1.
\end{prop}
\begin{proof}  Consider the evaluation of $N$ in the central fibre
$\{0\}\times M.\;$ The evaluation map is locally covering map from
generic points in $N.\;$ Since $\p N$ is a set of isolated
singular disks and the evaluation map is continuous at this set,
then the image of $N$ must be $\{0\}\times M\;$ entirely. Since
$M$  is connected, the covering indice must be some positive
constant $k\geq 1$ for generic point.   In particular,
$\cU_{\phi_0}$ must be decomposed into $k$ disjoint connected
component. And each connected component gives rise to a nearly
smooth foliation.   However, a nearly smooth foliation is unique.
Thus, there is only one connected component and covering indice is
$1.\;$
\end{proof}

Now, we introduce the notion of {\it partially smooth foliation},
which arises from limits of almost super regular foliations under
convergence of boundary maps in suitable norms.

\begin{defi} \label{defo:inducedfoliation} For any boundary map ${\phi_0} \in C^\infty(\p \Sigma,
 \cal H), $  an open $2n$ dimensional subset $\cU_{\phi_0} \subset {\cal
M}_{\phi_0}$ and closed subset $\tilde \cU_{\phi_0} \subset {\cal
M}_{\phi_0} $ is called a {\bf partially smooth foliation} if the
following conditions are met:
\begin{enumerate}
\item $\cU_{\phi_0}\subset \overline{\cU}_{\phi_0}\subset \tilde
\cU_{\phi_0}.\;$ \item Every disk in $\cU_{\phi_0}$ is super
regular.
 \item 
The evaluation map $\pi\circ ev: \Sigma\times \bar \cU_{\phi_0}
\rightarrow \Sigma \times M$ is a continuous onto map into its
image where the image is dense in $\p\Sigma\times M\;$ Moreover,
the image of $\tilde \cU_{\phi_0}$ is $\Sigma\times M.\;$
 \item  Any disk in $\cU_{\phi_0}$
doesn't intersect with any other disk in $\tilde \cU_{\phi_0}$ in
$\Sigma^0 \times M.\;$

\end{enumerate}
\end{defi}
Recall that an almost smooth solution of the HCMA equation
(\ref{eq:hcma0}) corresponds to a nearly smooth foliation. One can
view a {\it partially smooth foliation} as a sequential limit of
{\it nearly smooth foliations}, while a partially smooth solution
can be viewed as a sequential limit of almost smooth solutions. In
this regard, a partially smooth solution  corresponds conceptually
to a partially smooth foliation.

\subsubsection{Open-denseness of almost super regular disks}
 In this subsection, we reformulate Theorem
\ref{defo:genericpath1} in  terms of a.s.r or s.r. holomorphic
disks. Let us first describe some properties of holomorphic disks
with either partial indices $(0,0,\cdots, 0)$ or $(1,0,\cdots
0,-1).\;$ 

\begin{theo}
\label{th:genericdisksaresuperregular} Given a connected component
of ${\cal M}_{{\phi_0}}$ which consists of holomorphic disks with
partial indices $(0,0,\cdots ,0)$. If there exists at least one
super regular disk in this connected component, then all disks in
this component are super regular.
\end{theo}
\begin{proof} The set of super regular disks is open
in the {\it moduli} space. Therefore, we just need to show that it
is close among disks with partial indices $(0,0,\cdots,0 )$. Let
$\{f_i: (\Sigma,\p \Sigma) \rightarrow (\Sigma\times \cW_M, \b
\Lambda_{\phi_0})\}$
be a sequence of super regular disks such that $f_i \rightarrow f$
smoothly (cf. in $C^{2,\alpha}(\Sigma, \Sigma\times \cW_M)$-
norm\footnote{This regularity assumption is not optimal.}) in
$\cM_{\phi_0}$. We want to prove that if $f$ has partial indices
$(0,0,\cdots,0)$, then $f$ is also super regular.

Since $f$ is regular, there is a small neighborhood $\cO_f \subset
 \cM_{\phi_0}$ of $f$ such that $ev: \Sigma\times \cO_f \rightarrow \Sigma\times
 \cW_M$ is smooth. Put $F = \pi\circ ev$. Let $t_1,t_2, \cdots, t_{2n}$ be local coordinates of
$\cO_f$, write
\[
s_k^{(i)}(z) = {{\p\;ev}\over {\p t_k}}\mid_{f_i(z)} \,\in\,
T_{f_i(z)}^{1,0} \cW_M, \qquad 1\leq k\leq 2n.
\]
Then, $\{s_k\}_{k=1}^{2n}$ is a basis of the kernel space of the $\b \partial $ operator.  Moreover, at each
image point $f(z)$,  the set of $2n$ vertical vectors $\{s_k\}_{k=1}^{2n}$  is also a  basis of the ``vertical" tangential subspace $T_{f(z)} \cM_{\phi_0}$.
Since $\cW_M$ is locally the same as $T^*M $, its tangent space
naturally splits into a $TM$ part and
the tangent space of the fibre direction. Denote by $\left(\begin{array}{l}  u \\
v\end{array}\right)$ the corresponding two components of any
kernel vector, where $u,v\in C^n$. Set the k-th kernel vector as
\[
s_k^{(i)}  =\left(\begin{array}{c} u_k^{(i)}\\v_k^{(i)}
\end{array}\right), \qquad 1\leq k \leq 2n.
\]
Note that $u_k^{(i)}$ is clearly a vector in $T^{1,0}M$, while $v_k^{(i)}$ may depends on
the choice of local K\"ahler potentials.
According to Proposition \ref{defo:localpotentialforsuperregular},
there exists a solution $\phi^{(i)}$ of equation (\ref{eq:hcma0})
in $\pi\circ ev (\Sigma\times \cO_f)$ with
${\phi^{(i)}}\mid_{\p\Sigma\times M} = {\phi_0}\mid_{\p
\Sigma\times M}$. By Proposition 2.3.8, there is a uniform $C$
such that
\[ \mid \p \b \p \phi^{i} \mid \leq C. \]

For any point $(z,x)$ in the image of $\pi\circ ev (z, f),$ let
$U$ be a small open set of $x$ in $M.\;$ Then
\[
  G^{(i)} = \left(g_{0,\alpha\b\beta} + {{\p^2 \phi^{(i)}}\over {\p w_\alpha\p
\b w_\beta}}\right)_{n\times n}>0, \qquad\; S^{(i)} = \left({{\p^2
(\rho+ \phi^{(i)})}\over {\p w_\alpha \p w_\beta}} \right),
\]
where $\omega_0 = g_{0,\alpha\b \beta}\;d w^{\alpha}\wedge d\,\b w^\beta
= \sqrt{-1} \p \b \p \rho $ in $U\;$.
By a straightforward calculation, we have that for any $(z,w) \in
f_i(\Sigma)$,
\begin{equation}
  \left(v_k^{(i)}\right)_{n\times 1} = G_{n\times n}^{(i)} \cdot
\left(\bar u_k^{(i)}\right)_{n\times 1} + S_{n\times n}^{(i)}
\cdot \left( u_k^{(i)}\right)_{n\times 1}. \label{eq:super0}
\end{equation}
>From this equation, it is clear that $v_k^{(i)}$ is not ``tensorial" in the usual sense since $S$ is not.
However,
\begin{equation}
\det \left(\begin{array}{llll} u_1^{(i)} & u_2^{(i)} & \cdots & u_{2n}^{(i)}\\
v_1^{(i)} & v_2^{(i)} &\cdots &  v_{2n}^{(i)}\end{array}\right) =
\det \left(\begin{array}{llll} u_1^{(i)} & u_2^{(i)} & \cdots &
u_{2n}^{(i)}\\\bar u_1^{(i)} & \b u_2^{(i)} &\cdots & \b
u_{2n}^{(i)}\end{array}\right) \cdot \det G^{(i)}
\label{eq:super2}
\end{equation}
is both real and holomorphic.  Note that the right hand side is a function independent of
the choice of local coordinate in $M.\;$  Thus, the left hand side is well defined in $\Sigma$
and it must be a global constant along the disc.   Suppose this constant
$c_i$ \footnote{When restricted to each $\{z_0\}\times M,$ our
normalization forces the first term on the right hand side to be
positive. } and
\begin{equation}
\det \left(\begin{array}{llll} u_1^{(i)} & u_2^{(i)} & \cdots &
u_{2n}^{(i)}\\\bar u_1^{(i)} & \b u_2^{(i)} &\cdots & \b
u_{2n}^{(i)}\end{array}\right) \cdot \det G^{(i)}=c_i.
\label{eq:super3}
\end{equation}
  A more global view of
(\ref{eq:super3}) is
\[
 { \left(\pi\circ ev\right)^* \omega_{\phi^{(i)}}^n}    = {\omega_0^n}\mid_{z=z_0\in
    \p\Sigma}.
\]

Since $f_i$ is a super regular disk, by definition, we have
\begin{equation}
\det \left(\begin{array}{llll} u_1^{(i)} & u_2^{(i)} & \cdots &
u_{2n}^{(i)}\\\bar u_1^{(i)} & \b u_2^{(i)} &\cdots & \b
u_{2n}^{(i)}\end{array}\right) \neq 0 \label{eq:super1}
\end{equation}
holds everywhere in $f_i(\Sigma)$. Note that the left side of
(\ref{eq:super1}) is exactly the Jacobian of $\pi\circ ev$. By our
assumptions, $f$ is a disk with partial indices $(0,0,\cdots, 0)$
, that is, that the kernel matrix is nowhere singular:
\[
\det \left(\begin{array}{llll} u_1^{(i)} & u_2^{(i)} & \cdots & u_{2n}^{(i)}\\
v_1^{(i)} & v_2^{(i)} &\cdots &  v_{2n}^{(i)}\end{array}\right)
\neq 0
\]
everywhere in $f(\Sigma).\;$  What we need to prove is that the
inequality (\ref{eq:super1}) hold everywhere in $f(\Sigma).\;$
 For $f_i(\Sigma),$ it is easy to see that
This sequence of constants $\{c_i, i \in \cN\}$ has both uniform
upper and lower bounds, provided that the limiting disk $f$ has
partial indices $(0,0,\cdots, 0)$. Since $\det G^{(i)} \leq C$, we
deduce
\[
 \det \left(\begin{array}{llll} u_1 & u_2 & \cdots &
u_{2n}\\\bar u_1 & \b u_2 &\cdots & \b u_{2n}\end{array}\right)
> 0
\]
along the limiting disk $f$. This completes the proof of this
theorem.
\end{proof}

However, we can squeeze a little more out from the arguments
above. Let $f$ be a holomorphic disk with partial indice
$(1,0,\cdots,0,-1)$ which is super regular at $z=0$. We claim that
$f$ is almost super regular.  In fact, the condition implies that
\[
\det \left(\begin{array}{llll} u_1^{(i)} & u_2^{(i)} & \cdots &
u_{2n}^{(i)}\\\bar u_1^{(i)} & \b u_2^{(i)} &\cdots & \b
u_{2n}^{(i)}\end{array}\right) \mid_{z=0} \neq 0.
\]
On the other hand, according to (\ref{eq:super2}), for each
$f_i$,we have
\[
 c_i =
\det \left(\begin{array}{llll} u_1^{(i)} & u_2^{(i)} & \cdots & u_{2n}^{(i)}\\
v_1^{(i)} & v_2^{(i)} &\cdots &  v_{2n}^{(i)}\end{array}\right) =
\det \left(\begin{array}{llll} u_1^{(i)} & u_2^{(i)} & \cdots &
u_{2n}^{(i)}\\\bar u_1^{(i)} & \b u_2^{(i)} &\cdots & \b
u_{2n}^{(i)}\end{array}\right) \cdot \det G^{(i)}
\]
Using local deformation theory in next subsection (Corollary
\ref{local:linecurvature}), $\log \det G^{(i)}$ is a subharmonic
function in $\Sigma$. Moreover, it is uniformly bounded from
above. Set
\[
h_i = \log \det \left(\begin{array}{llll} u_1^{(i)} & u_2^{(i)} &
\cdots & u_{2n}^{(i)}\\\bar u_1^{(i)} & \b u_2^{(i)} &\cdots & \b
u_{2n}^{(i)}\end{array}\right)
\]
along $f_i(\Sigma)$. Then, $\{h_i, i \in \NN\}$  is a seequence of subharmonic functions
 on $f_i(\Sigma)$ which is  uniformly
bounded at $z=0.\;$ Moreover, we have
(cf. Proposition \ref{local:deltazpotential})
\[
\mid \Delta_z \; h_i\mid =\mid - \Delta_z\;\log \det G^{(i)}\mid
\leq C,.
\]
holds in any compact subdomain of $\Sigma^0.\;$
Then Harnack inequality for harmonic function implies that either
$h_i$ approaches to $-\infty$ everywhere in any compact subset of
$\Sigma^0$ or stays uniformly bounded in any compact subset of
$\Sigma^0$.  Since $h_i(0)$ is uniformly bounded, we have
\[
\displaystyle \lim_{i\rightarrow \infty}\;h_i(z) = h(z) =\log \det
\left(\begin{array}{llll} u_1^{(i)} & u_2^{(i)} & \cdots &
u_{2n}^{(i)}\\\bar u_1^{(i)} & \b u_2^{(i)} &\cdots & \b
u_{2n}^{(i)}\end{array}\right)\mid_{f(\Sigma)}
> -\infty,\qquad \forall z \in \Sigma^0.
\]
Consequently, $f$ is almost super regular. Thus we have proved

\begin{theo} If a holomorphic disk with partial indices
$(1,0,\cdots, 0,-1)$ is super regular at one interior point and it
can be connected to disks of partial indices $(0,0,\cdots, 0)$,
then it is almost super regular.
\end{theo}

In view of these two theorems, we can reinterpret Theorem
\ref{defo:genericpath1} as
\begin{theo}\label{defo:genericpath2}
Given any path $\psi: [0,1] \rightarrow C^\infty(\p \Sigma, {\cal
H})$ such that ${\cal M}_{{\psi}(0,\cdot)}$ contains a super
regular disk with vanishing Maslov disk, there exists a generic
path (still denoted by $\psi$ for simplicity), which is
arbitrarily close to the given path, such that for this new path,
a connected component $\cM_\psi^0$ of the total moduli space
$\bigcup_{0\leq s \leq 1}\; \{s\}\times {\cal M}_{\psi(s,\cdot)}$,
which contains the initial super regular disk, is a smooth
$2n+1$-dimension manifold.
Moreover, the followings hold
\begin{enumerate} \item The set of super regular
disks is open and dense in this connected component; \item The set
of almost super regular disks has codimension at least 1 in this
component; \item The set of disks, which are neither super regular
nor almost super regular, has codimension at least 2; \item There
exist at most finitely many irregular disks in the total moduli
space.
\end{enumerate}
Moreover, there exist at most finitely many points $0 < \b t_1 <
\b t_2 < \cdots <\b t_l < 1$ such that for any $ t \neq \b t_i
(1\leq i\leq l)$, all disks in ${\cal M}_{{\psi}(t,\cdot)}$ are
regular and its subset of disks in $\cM_\psi^0$, which are neither
super regular nor almost super regular, has codimension at least
2. When $t=\b t_i$ for some $i$, ${\cal M}_{{\psi}(t,\cdot)}\cap
\cM_\psi^0$ may either contain some isolated irregular disks or a
subset of disks which are neither super regular nor almost super
regular which has exactly codimension $1.\;$
\end{theo}

\subsubsection{Almost super regular foliations along a generic path}
In this subsection, we prove openness and closedness of almost
super regular foliations along a generic path, which is assured by
Theorem \ref{defo:genericpath2}.

\begin{theo}
\label{defo:openclosenessforasf} Let $\psi:[0,1]\mapsto
C^\infty(\Sigma,\cH)$ be a generic path with properties specified
in Theorem \ref{defo:genericpath2}. Suppose that
$\cM_{\psi(0,\cdot)}\cap \cM_\psi^0$ is connected and gives rise
to an {\bf almost super regular foliation}, where $\cM_\psi^0$ is
the connected component defined in Theorem
\ref{defo:genericpath2}. Then for each $t$,
$\cM_{\psi(t,\cdot)}\cap \cM^0_\psi$ is connected and induces a
foliation in an open dense subset of $\Sigma\times M.$ Moreover,
this component gives rise to an almost super regular foliation
except at most a finite number of times.
\end{theo}
We firs prove
\begin{lem} \label{defo:closenessforasf}For a sequence of $\tau_i, i\in
\NN (\displaystyle \lim_{i\rightarrow \infty} \tau_i = \b t \in
(0,1])$ such that $\cF_{\psi(\tau_i,\cdot)}$ is a sequence of
almost smooth foliations. Suppose that $\phi_i$ is the
corresponding sequence of almost smooth solutions and
$\displaystyle \lim_{i\rightarrow\infty} \phi_i =\phi_\infty.\;$
Then $\phi_\infty$ is a partially smooth solution of
(\ref{eq:hcma0}) and $\cR_{\phi_\infty}$ is an open dense subset
of $\Sigma\times M.\;$ Moreover, there is a unique connected
component of $\cM_{\psi(\b t, \cdot)}$ which is the limit of
$\cF_{\psi(\tau_i,\cdot)}.\;$ Either this component is regular in
which case there is an almost super regular foliation
$\cF_{\psi(t,\cdot)}$ for $t>\b t;$ or this component is an almost
super regular foliation itself.

\end{lem}
\begin{proof} Following Theorem \ref{comp:foliationlimit}, after
passing to a subsequence, $\cF_{\psi(\tau_i,0)}$ converges to a
partially smooth foliation $\cF_{\psi(\tau_\infty,\cdot)}$, where
$\cU_{\psi(\tau_\infty,\cdot)}$ denotes the set of all its super
regular disks. Theorem \ref{comp:superregularlimit} implies that
there is at least one super regular disk which is the limit of a
sequence of super regular disks in $\cF_{\phi_0(\tau_i)}$ with
uniformly bounded capacity (cf. formula \ref{eq:capacity}).
Therefore, $\cU_{\psi(\tau_\infty,\cdot)}$ is non-empty.\\

For convenience, let $\cB$ be the union of all disks in ${\cal
M}_{\psi(\tau_\infty,\cdot)}$ which are sequential limits of disks
in $\cU_{\phi_0(\tau_i)}$. By definition,  any leaf in $\cB,$  is
the limit of some sequence of disks in $\cM_{\psi_{\tau_i}}.\;$
Following Theorem \ref{comp:nobubble},
 for any such sequence of disks, the corresponding sequence of leaf
vector field in $TM$ has a uniform upper bound on length.   In
particular, all leaves in $\cB$ have a uniform upper bound on the
length of their leaf vector fields.  Consequently, any sequence in
$\cB$ must have a convergent subsequence where the limit is an
embedded disk in $\cM_{\tau_\infty}.\;$ It follows that $\cB$ is a
closed, bounded set in the moduli space. On the other hand, by the
choice of our generic path in Theorem \ref{defo:genericpath2}, the
moduli space at $t=\tau_\infty$ admits at most finitely many
non-regular, embedded disks. Therefore, all disks in $\cB$, except
at most a finite number of disks,  are regular. Consequently, the
evaluation map is continuous everywhere in $\cB$ and
differentiable  except at most
a finite number of points (leafs). \\

Moreover, $\cU_{\psi(\tau_\infty,\cdot)}$ is an open dense, and
irreducible subset
of $\cB$. 
If $\tau_\infty \neq \b t_k, (1\leq k\leq l)$, then all disks in
$\cB$ are regular and the set of disks which are neither super
regular nor almost super regular has codimension at least $2$. In
this case, $\cB = 
\cF_{\psi(\tau_\infty,\cdot)}$ is an almost super regular
foliation. On the other hand, if $\tau_\infty = \b t_k$ for some
$k$, then either $\cB$ contains  a finite number of  singular
disks, or all disks in $\cB$ are regular where the codimension for
non-almost super regular or non-super regular disks may be  1. In
the case that $\cB$ contains a finite number of isolated disks,
the codimension of non super regular disks or non almost super
regular disks must have codimension 2 or higher. In this case,
$\cB$ defines an almost super regular foliation.  The last
remaining case is that $\cB$ is regular but the the set of
non-super regular or non-almost super regular disks may have
codimension 1. In this case,  we can perturb this component $\cB$
for $t-\b t>0$ small. Following Proposition
\ref{defo:almostsuperregular4}, the connected component after
perturbation
 defines an almost super regular foliation for $t>\b t.\;$\\

In all cases, it is easy to see that $\phi_\infty$ is smooth in an
open dense subset $\cR_{\phi_\infty}$ of $\Sigma\times M.\;$
Moreover, we can show that $\cB$ is unique since the corresponding
partially smooth solution in the limit  is unique.   
This in particular implies that the limit
$\cF_{\psi(\tau_\infty,\cdot)}$ is
independent of the time sequence $\tau_i\rightarrow \infty.\;$\\

\end{proof}

Now we return to prove our main theorem.

\begin{proof} To prove openness, we assume $ \cF_{\psi(\b t,\cdot)}$ is an almost super regular
foliation. Here we follow the notations in Theorem
\ref{defo:genericpath2}. Without loss of generality, we may assume
$\b t\leq \b t_1$. If $\b t < \b t_1$, then $\cM_{\psi(\b
t,\cdot)}$ is smooth. In particular, the connected component $\b
\cU_{\phi(\b t,\cdot)}$ is smooth without boundary. Following from
the standard deformation theory,   this component will deform to a
smooth component $\b \cU_{\psi(t,\cdot)}$ of $\cM_{\psi(t,\cdot)}$
for $t - \b t > 0$ small enough. By Theorem
\ref{th:genericdisksaresuperregular} and \ref{defo:genericpath2},
$\b \cU_{\psi(t,\cdot)}$ induces an almost super regular
foliation, so the openness follows in this case.

Now assume $\b t = \b t_1$. 
We want to show that for $\b t_1 = \b t < t < \b t_2 $, there is
an almost super regular foliation $\cF_{\psi(t,\cdot)}.\;$

By our choice of the generic path $\psi$, we may assume that there
are a finite number of embedded, non-regular disks in
$\b\cU_{\psi(\b t,\cdot)}\backslash \cU_{{\psi}(\b t,\cdot)}.\;$
Since we are interested in preserving this connected component
$\cU_{\phi(\b t,\cdot)}, $ we want to rule out the possibility of
either a ``merge in" or ``spin off" occuring.  In other words,
there might be another component of $\cM_{\psi(\b t,\cdot)}$
connecting with $\b \cU_{{\psi}(\b t,\cdot)}$ through these
isolated singular disks: Two components before $t = \b t$ may
merge locally into one smooth connected component after $t= \b t$.
The situation can also occur in the reverse order: an open set of
the moduli $\cM_{\psi(t,\cdot)}$ may pinch off of a ``neck
$S^{2n-1}$" at $t = \b t$ and go on to become two
separate components after $t = \b t,\;$ at least locally near this ``neck.." 
We call the first case ``merge-in" and the second case
``pinching-off''. If either one occurs, this ``good" component
will change after singular disks.   
The deformation of almost super regular foliations is impossible
if either of these phenomenon occur beyond the time when singular
disks appear. We will deal only with the ``merge-in'' case here,
since the other cases (like `the `pinching-off'' case) can be
handled in a similar fashion.

Note that the ``merge-in" of the {\it moduli} spaces occurs only
at singular disks. Since there is only finite number of singular
disks and ``merge in" only occur locally near singular disk, we
may assume without loss of generality, there is
only one non-regular disk $\b f$ in $\b \cU_{\psi(t)}.\;$\\

Without loss of generality, set $\b t = \b t_1$ and $\cM_t =
\cM_{\psi(t,0)}$ is an almost super regular foliation for any $ t
\in (0, \b t]$ . Suppose $\b f$ is the only isolated singular disk
at $t= \b t.\;$ Then the metric ball $B_r(\b f)$  in $\cM_{\b t}$
can be represented by a cone in $R^{2n+1}$:
\[
  \displaystyle \sum_{i=1}^k x_i^2 - \displaystyle \sum_{i=k+1}^{2n+1}
  x_i^2 = 0
\]
where
\[
\displaystyle \sum_{i=k+1}^{2n+1} x_i^2 < r^2.
\]
Here $(0,0,\cdots, 0)$ represent $\b f.\;$ The ``merge in" or
``spin off" case corresponds to $k=2n$ or $k = 1.\;$ We diskuss
the ``merge in" case here. For $ \b t - t
> 0 $ small, the corresponding metric ball in $\cM_t$ is
\[
  \displaystyle \sum_{i=1}^{2n} x_i^2 -
  x_{2n+1}^2 = t -\b t, \qquad x_{2n+1}^2 <
  r^2 + \b t- t,\;\; x_{2n+1} > 0.
\]
 For $  t -\b  t
> 0 $ small, the corresponding metric ball in $\cM_t$ is
\[
  \displaystyle \sum_{i=1}^{2n} x_i^2 -
  x_{2n+1}^2 = t -\b t, \qquad  \displaystyle \sum_{k=1}^{2n} x_k^2 \leq
  r^2 + t - \b t.
\]
Choose a continuous path of disks $f(t)\in \cM_{\psi(t,\cdot)}$
such that $f(t) = \b f$ and $f(t) (t \neq \b t)$ is either super
regular or almost super regular disk. For notational simplicity,
we denote $f(t)$ by $\b f.\;$ Note that for $r>0$ small enough,
the intersection $B_r(\b f)\cap \cM_t$ consists of two disjoint
disks for $ t\leq \b t, $ but is cylinder-like for $t
> \b t.\;$ We consider the intersection of this ball with the
central fibre $\{0\}\times M$. In this proof, we use $ev$ to
denote the map
$\pi\circ ev(0,\cdot).\;$  Set $ev(\b f) = p.\;$ \\

Note that for $t > \b t$, the boundary of $B_r(\b f)$ is made of
two components $ N_1 \approx N_2 \approx S^{2n-1}$ ($\approx$
means diffeomorphic to) which are homotopic to each other in $B_r
(\b f)$. These boundary spheres are perturbations of $\p B_r(\b
f)\cap \cM_{\psi(\b t, \cdot)}$. Pick up one of these spheres, say
$N_1$ for $t > \b t.\;$ Since each component of $ev(B_r(\b f)\cap
\cM_{\psi(\b t,\cdot)})$ bounds a deformation retractable domain
in the central fibre $\{0\}\times M$. By continuity, $ev(N_1)$
also bounds a domain $\Omega$ which is deformation contractible to
an interior point $q \in \Omega$ for $t-\b t$ sufficiently small.
Let us denote this contraction by $F: [0,1]\times ev(N_1)\mapsto
\Omega$ such that $F(0,p)=p$ and $F(1,p)=q$ for any $p\in
ev(N_1)$. Since the set of disks which are neither super regular
nor almost super regular has codimension 2 or higher, there is an
open subset $V\subset ev(N_1)$ such that $F([0,1)\times V)$ does
not intersect with the image of the set of disks which are neither
super regular nor almost super regular. Now we can lift this
$F([0,1)\times V)$ to $\cM_t$ since any point in the subset has
its pre-image covered by either an super regular disks or an
almost super regular disk. This implies that there is a subset
$N_3$ such that $ev(N_3)$ is a single point, where $N_3$ consists
of all limiting points of the lifting of $[0,1)\times V$. Clearly,
any disk in $N_3$ is neither super regular nor almost super
regular. Observing that $N_3$ has codimension one, we get a
contradiction to the fact that the set of all disks which are not
almost super regular has codimension at least two.

By similar arguments, we can prove that there is no
``pinching-off" at $t=\b t$.
\end{proof}


\section{\label{local:defomation}Basic curvature equations along leaves}
\subsection{Introduction}
In this section, we show some deformation results for solutions of the
homogenous complex Monge Ampere equations. In particular, we give
a basic curvature formula for the restriction of involved metric
to leaves. This formula plays a crucial role in deriving key {\it
a priori} estimate. Suppose that $\phi$ is a solution of
(\ref{eq:hcma0}). Suppose that the $\pi_2^*\omega + \sqrt{-1}
\partial\overline{\partial} \phi$, referred as the {\it Levi}
form, has constant co-rank $1$. This gives rise to a foliation of
the domain by holomorphic disks. We further assume that:\\

{\bf At each point of the domain, the
leaf vector field is always transversal to a $M$.}\\

Under this assumption, the {\it Levi} form restricts to a K\"ahler
metric in $M$ for each $z\in \Sigma$. In this way, a solution of
(\ref{eq:hcma0}) can be alternatively viewed as a disk family of
K\"ahler metrics satisfying certain geometric conditions.  To understand the geometry of this
family of K\"ahler metrics,  we will
study the restriction of the complex tangent bundle $TM$ over this
family of holomorphic disks. These bundles are equipped with
natural Hermitian metrics, so we obtain a family of Hermitian bundle
over disks. In this section, we will compute curvature of these
Hermitian bundles. The main results are
\begin{enumerate} \item The curvature of these Hermitian metrics is always non-positive
(Theorem \ref{local:curvaturepositive});
\item The foliation is holomorphic if and only if the ``trace of
the curvature" of these Hermitian metrics vanishes (Theorem
\ref{local:vanishing curvature}); \item The trace of the Hermitian
curvature is always super harmonic (Corollary
\ref{local:linecurvature}).
\end{enumerate}

The results in this section lay foundations for the global deformation
of almost super regular foliations in this paper.

\subsection{Hermitian Curvature formulas}
In local coordinate, we write
\[
\omega_0 = \sqrt{-1}\sum_{\alpha,\beta=1}^n g_{0,\alpha\b\beta}\;
d\,w^\alpha\;\wedge d\,w^{\b\beta}, \qquad \omega_\phi =
\sqrt{-1}\sum_{\alpha,\beta=1}^n g_{\phi,\alpha\b\beta}\;
d\,w^\alpha\wedge d\,w^{\b\beta}
\]
where
\[
g_{\alpha \overline{\beta}} = g_{0,\alpha \overline{\beta}} +
{{\partial^2 \phi}\over{\partial w_{\alpha}
\partial \bar w_{{\beta}}}} , \qquad \forall
\alpha,\beta=1,2 \cdots n.
\]
As before, $z$ denotes the coordinate variable of $\Sigma$. Then
(\ref{eq:hcma0}) can be re-written as

\begin{equation}
   {{\partial^2 \phi}\over{\partial z
{\partial \b z}}} - g_\phi^{\alpha \overline{\beta}} {{\partial^2
\phi}\over{\partial z \partial \b w_{\beta}}} {{\partial^2
\phi}\over{\partial \b z}
\partial w_{\alpha}} = 0.
\label{eq:geodesic}
\end{equation}
Here we are assuming that $\omega_\phi = \omega_0 +\sqrt{-1} \p \b
\p \phi
> 0$ in $M.\;$  In this section, we will simply write $g$ for the metric $g_\phi$ if there is no confusion.
Write the leaf vector field (cf. Definition \ref{def:leafvector1})
as
\begin{equation}
X = \sum_{\alpha=1}^n\; \eta^\alpha {\p \over \p w_\alpha} =
\sum_{\alpha=1}^n\; - g^{\alpha\b{\beta}} \; {{\partial^2
\phi}\over{\partial z {\partial \b w_{\beta}}}} {\p \over \p
w_\alpha}. \label{loca:leafvector}
\end{equation}

Denote the linearized operator by $\Delta_z$. There is  a natural
splitting of this operator since all disks are holomorphic:
\begin{equation}\Delta_z = \p_z \bar \partial_{z}, \qquad {\rm where}\;\;
\partial_z = {{\partial }\over {\partial z}} + \eta^{\alpha} {{\partial }\over{\partial w_{\alpha}}}.
\label{diskderiv}
\end{equation}


\begin{prop} \label{local:diskbrake} The leaf vector field $X$ (cf. Defi. \ref{def:leafvector} ) is holomorphic in
$z.\;$ In other words
\begin{equation} [\partial_z, \overline{\partial_z} ] = {\p_{\b z} X} = 0.
\label{diskbrake}
\end{equation}
\end{prop}
\begin{proof} Direct calculation.\end{proof}
\begin{prop} \label{local:commutator} The commutator of local differtial operator on $TM$   and the leaf derivative $\p_z$  is
\begin{equation}
[\partial_z, {{\partial} \over{\partial w_i}}] = - {{\partial
\eta^{\alpha}} \over {\partial w_i}} {{\partial }\over{\partial
w_{\alpha}}},\qquad [\partial_z, {{\partial} \over{ \partial
\overline{w_i}}} ] = - {{\partial \eta^{\alpha}} \over {\partial
\overline{w_i}}} {{\partial }\over{\partial w_{\alpha}}}.
\label{commutator}
\end{equation}
\end{prop}
\begin{proof}
\[
[\partial_z,  {{\partial} \over{\partial w_i}}]  =  [{{\partial
}\over {\partial z}} + \eta^{\alpha}
 {{\partial }\over{\partial w_{\alpha}}}, {{\partial} \over{ \partial w_i}}] = - {{\partial \eta^{\alpha}} \over {\partial w_i}}
  {{\partial }\over{\partial w_{\alpha}}},
\]
and
\[
[\partial_z,  {{\partial} \over{\partial \overline{w_i}}}]  =
[{{\partial }\over {\partial z}} + \eta^{\alpha} {{\partial
}\over{\partial w_{\alpha}}}, {{\partial} \over{ \partial w_i}}] =
- {{\partial \eta^{\alpha}} \over {\partial \overline{ w_i}}}
{{\partial }\over{\partial w_{\alpha}}}.
\]
\end{proof}
\begin{rem} Note that  that ${{\partial \eta^{\alpha}} \over
{\partial w_i}}$ is not a globally well defined tensor, while $
{{\partial \eta^{\alpha}} \over {\partial \overline{ w_i}}}$ is a
globally well defined tensor since
\[\begin{array}{lcl}
  {{\partial \eta^{\alpha}} \over
{\partial \overline{ w_i}}} & = & - g^{\alpha \b s} \left( {{\p
\phi} \over {\p z}}\right)_{,\b i \b s}.
\end{array}
\]

\end{rem}

In a local coordinate chart, suppose $g_{0,\alpha
\overline{\beta}} = {{\partial^2 \rho} \over{ \partial w_{\alpha}
\overline{\partial w_{\beta}}}}$ where $ \rho$ is independent of
the $z$ direction.
 Let $\Phi$ denote the local K\"ahler potenital for $g,$ then
\[
  \Phi = \phi + \rho.
\]

\begin{prop} \label{local:deltazpotential} For global K\"ahler distortion potentials, the following
is true:
\begin{equation}
  \Delta_z \phi =\p_z \bar \p_z (\phi) =- |X|^2_{g_0}.
  \label{local:deltazpotential1}
\end{equation}
\end{prop}
\begin{proof} By a straightforward calculation.
\end{proof}

\begin{lem}\label{local:bootstrap} (The Bootstrapping lemma): The following commutation formula
for the third transversal derivatives holds
\begin{equation}
{{\partial \eta^{\alpha}}\over {\partial w_i}} = - g^{\alpha {\b
\beta}} \;\partial_{z} g_{i {\b \beta}},\qquad {\rm and}\qquad
{{\partial \eta^{\alpha}}\over {{\partial \b w_i}}} = - g^{\alpha
{\b \beta}} \;\partial_{z} {{\partial^2 \Phi}\over{{\partial \b
w_i}{\partial \b w_\beta}}}. \label{eq:bootstrapping}
\end{equation}
\end{lem}
\begin{proof}In local coordinate, we have \[
\begin{array} {lcl} \partial_{z} g_{i{\b \beta}} & = &
\partial_{z} {{\partial^2 \Phi}\over {\partial w_i{\partial \b w_{\beta}}}}
= {{\partial }\over {\partial w_i}} \partial_{z} {{\partial
\Phi}\over {{\partial \b w_{\beta}}}} - {{\partial \eta^{\alpha}}
\over {\partial w_i}} {{\partial \Phi}\over{{\partial\b
w_{\beta}}\partial w_{\alpha}}} \\ & = & -{{\partial
\eta^{\alpha}} \over {\partial w_i}} g_{\alpha {\b \beta}}.
 \end{array}
\]
On the other hand,
\[
\begin{array} {lcl}
\partial_{z} {{\partial^2 \Phi}\over {\overline{\partial w_i} \overline{\partial w_{\beta}}}} & = & {{\partial }\over {\overline{\partial w_i}}} \partial_{z} {{\partial \Phi}\over { \overline{\partial w_{\beta}}}} -  {{\partial \eta^{\alpha}} \over {\overline{\partial w_i}}} {{\partial \Phi}\over{ \overline{\partial w_{\beta}}\partial w_{\alpha}}} \\ & = & - {{\partial \eta^{\alpha}} \over {\overline{\partial w_i}}} g_{\alpha \overline{\beta}}.
\end{array}
\]
\end{proof}
\begin{rem} The significance of this equation (\ref{eq:bootstrapping}) is that it changes
the 1st derivatives on the transversal direction into the 1st
derivatives along the disk of the 2nd order jet of local K\"ahler
potentials.
\end{rem}
\begin{lem}\label{local:regularity}(Regularity lemma). The following equation holds along the disk
(for any $\alpha, i =1,2,\cdots n$):
\begin{equation}
\overline{\partial_z} {{\partial \eta^{\alpha}}\over {\partial
w_i}} = - \overline{\left({\partial \eta^{\beta}}\over
{\overline{\partial w_i}}\right)}\;{{\partial \eta^{\alpha}}\over
{\overline{\partial w_{\beta}}}},\qquad {\rm and}\qquad
\overline{\partial_z} {{\partial \eta^{\alpha}}\over
{\overline{\partial w_i}}} = - \overline{\left({\partial
\eta^{\beta}}\over {{\partial w_i}}\right)}\;{{\partial
\eta^{\alpha}}\over {\overline{\partial w_{\beta}}}}.
\label{regularity}
\end{equation}
\end{lem}
\begin{proof} This lemma follow from the commutation formula
(\ref{commutator}) and Lemma 4.2.5 immediately.
\end{proof}

The major obstacle in establishing an {\it a priori} estimate for
(\ref{eq:hcma0}) is that its linearized operator $\Delta_z$ only
has
 rank $1.\;$ This is a severe restriction in deriving any
meaningful estimate directly. However,  we may restrict the
$k-(k=0,1,2\cdots)$th  jet of the potential function $\Phi$ along the
disk and study the geometric  equation it must satisfy.  Notice that the
restriction of $TM$ bundle to the disk is a trivial $C^n$ bundle.
At each point, $g_\phi$ is a Hermitian metric in this $TM$ bundle
over disk. Suppose that $\left(F^{\alpha}_{ \beta}\right)$ is the Hermitian
curvature of this Hermitian metric.
\begin{theo} \label{local:curvaturepositive}The curvature of the bundle metric is always non-positive.
\end{theo}
\begin{proof} \[\begin{array}{lcl} F^r_{\alpha} & =  & - \p_{\b z} \left(g^{r \overline{\delta}} \;\partial_{z} g_{\phi,\alpha
\overline{\delta}}\right) \\ & = &    \p_{\b z} {{\partial
\eta^{r}}\over {\partial w_{\b \alpha}}}
 =  - {{\partial \eta^{r}}\over
{\partial w_{\b i}}} \cdot {{\partial \eta^{\b i}}\over {\partial
w_\alpha}}.
\end{array}
\]
It is not difficult to see that the last term in the right hand
side is a Hermitian Symmetric non-positive 2-tensors. For any
holomorphic section $s: \Sigma \rightarrow T^{1,0}M,\;$ we have
\[
   \begin{array}{lcl} F(s,s) & = & F^{\alpha}_\beta s^\beta s^{\b \gamma}
   {g _\phi}_{\alpha \b \gamma} \\
   & = &  - {{\partial \eta^{\alpha}}\over
{\partial w_{\b i}}} \cdot {{\partial \eta^{\b i}}\over {\partial
w_\beta}} s^\beta s^{\b \gamma}
   {g _\phi}_{\alpha \b \gamma} \\ & = & -  g^{\alpha \b b}  g^{\b i a} \left( {{\partial \varphi}\over {\partial z}}\right)_{, \b i \b b} \cdot  \left( {{\partial  \varphi}\over {\partial \b z}}\right)_{, \beta  a}  s^\beta s^{\b \gamma}  {g}_{\alpha \b \gamma} \\
  & = &-  g^{\b i a} \left( {{\partial \varphi}\over {\partial z}}\right)_{, \b i \b r} \cdot  \left( {{\partial  \varphi}\over {\partial \b z}}\right)_{, \beta  a}  s^\beta s^{\b \gamma}
   \leq 0.
   \end{array}
\]
The last inequality holds since $\phi$ is real valued potential function and the $(2,0)$ or $(0,2)$ part of the Hessian
is symmetric.
\end{proof}
A quick corollary follows
\begin{theo} \label{local:vanishing curvature}This foliation by
holomorphic disks is a holomorphic foliation if and only if the
curvature of these Hermitian metrics vanishes.
\end{theo}
\begin{proof} This is quite evident from Theorem \ref{local:curvaturepositive}.
\end{proof}
\begin{prop} $g^{i \overline{j}} {{\partial \eta^{\alpha}}\over {\overline{\partial w_{j}}}} $ is holomorphic along
the disk.
\end{prop}
\begin{proof}
\[
\overline{\partial_z} g^{i \overline{j}}  = - g^{i \overline{k}}
\;\overline{\partial_z} g_{\overline{k} l } \;
 g^{l\overline{j}} =g^{i \overline{k}}  {{\partial \eta^{\overline{\beta}}}\over {\overline{\partial w_{k}}}}
   g_{l \overline{\beta}} g^{l\overline{j}}  = g^{i \overline{k}}  {{\partial \eta^{\overline{j}}}\over {\overline{\partial w_{k}}}}
  \;  \]
Now, we have
\[\begin{array}{lcl}
\overline{\partial_z} \left( g^{i \overline{j}} {{\partial
\eta^{\alpha}}\over {\overline{\partial w_{j}}}} \right) & =&
\left( \overline{\partial_z}  g^{i \overline{j}}\right) {{\partial
\eta^{\alpha}}\over {\overline{\partial w_{j}}}}  +
g^{i \overline{j}} \overline{\partial_z} {{\partial \eta^{\alpha}}\over {\overline{\partial w_{j}}}} \\
& = & g^{i \overline{k}}  {{\partial \eta^{\overline{j}}}\over
{\overline{\partial w_{k}}}}
  {{\partial \eta^{\alpha}}\over {\overline{\partial w_{j}}}}
   - g^{i \overline{j}} \overline{\left({\partial \eta^{\beta}}\over {{\partial w_j}}\right)}
   \;{{\partial \eta^{\alpha}}\over {\overline{\partial w_{\beta}}}}=0.
\end{array}
\]
\end{proof}
\begin{cor}\label{local:linecurvature} \footnote{This was first
derived in 1996 by the first author and S. Donaldson using some
different methods.}.The anti-canonical line bundle in $M$ equipped
with  $\omega_\phi^n$ as Hermitian metric, restricted to a
holomorphic disk,  has non-positive curvature. More precisely, we
have
\begin{equation}
 \Delta_z \log \omega_\phi^n = \partial_z \overline{\partial_z} \log \omega_\phi^n= {{\partial \eta^{\alpha}}\over {
\overline{\partial w_i}}} \overline {\left({{\partial
\eta^{i}}\over { \overline{\partial w_{\alpha}}}}\right)} .
\label{local:linecurvature1}
\end{equation}
Note that the right hand side measures whether the given
geodesic(or disk version) is holomorphic or not.
\end{cor}
\begin{proof} In a local
coordinate, we have
\[
   \partial_z \log \omega_\phi^n = g^{\alpha \overline{\beta}} \partial_z g_{\overline{\beta}\alpha}
                    = - {{\partial \eta^{\alpha}}\over {\partial w_{\alpha}}}.
\]
Thus,
\[
\begin{array}{lcl}
  \overline{\partial_z} \partial_z \log \omega_\phi^n & = & - \overline{\partial_z}\;{{\partial \eta^{\alpha}}\over {\partial w_{\alpha}}}
 =  - \left( {{\partial}\over{\partial w_{\alpha}}} \overline{\partial_z} - {{\partial \eta^{\overline{\beta}}}\over {\partial w_{\alpha}}} {{\partial}\over{\partial w_{\overline{\beta}}}}\right) \eta^{\alpha} \\
& = & {{\partial \eta^{\overline{\beta}}}\over {\partial
w_{\alpha}}} {{\partial \eta^{\alpha}}\over {\partial
\overline{w_{\beta}}}} \geq 0.
\end{array}\]
\end{proof}
Set $- S = {F^{\alpha}}_{\alpha} = {{\partial
\eta^{\overline{\beta}}}\over {\partial w_{\alpha}}} {{\partial
\eta^{\alpha}}\over {\partial \overline{w_{\beta}}}} =
{g_\phi}^{\alpha\overline{p}} \; \phi_{, z \overline{i}
\overline{p}}\;{g_\phi}^{\overline{i} q} \;
 \phi_{, \overline{z} \alpha\;q} \geq 0\;$ (this denotes the covariant derivatives w.r.t. $g_\phi$.) as  the trace of the curvature
 of this $TM$ bundle with Hermitian metric $g_\phi.\;$
\begin{theo} \label{local:3rdderivative}  The trace of the Hermitian curvature
satisfies the following
\begin{eqnarray}
\partial_z \overline{\partial_z} S & = & \left(\partial_z {{\partial \eta^{\alpha}}\over {\partial \overline{w_i}}} -
{{\partial \eta^{\beta}}\over {\partial \overline{w_{i}}}}
\;{{\partial \eta^{\alpha}}\over {\partial w_{\beta}}}\right)
\left( \overline{\partial_z}
 {{\partial \eta^{\overline{i}}}\over {\partial w_{\alpha}}} - {{\partial \eta^{\overline{i}}}\over {\partial \overline{w_p}}}\;
 {{\partial \eta^{\overline{p}}}\over {\partial w_{\alpha}}}\right) \nonumber  \\ & & \qquad \qquad \qquad + 2 {{\partial \eta^{\alpha}}\over {\partial \overline{w_{i}}}}
{{\partial \eta^{\overline{p}}}\over {\partial w_{\alpha}}}
{{\partial \eta^{\beta}}\over {\partial \overline{w_p}}}
{{\partial \eta^{\overline{i}}}\over {\partial w_{\beta}}} \label{eq:C3estimate} \\
& = & {{\partial}\over {\partial \overline{w_i}}} \partial_z
(\eta^{\alpha}) \; {{\partial}\over {\partial {w_{\alpha}}}}
\overline{\partial_z} (\eta^{\overline{i}})\nonumber
\\ & & \qquad
+  2 {g_\phi}^{\alpha \overline{i}} {g_\phi}^{ \overline{r}j }
{g_\phi}^{\beta \overline{k}}  {g_\phi}^{ \overline{\delta}l }
\phi_{,z \overline{i} \overline{r}} \phi_{,\overline{z} j \beta}
\phi_{,z \overline{k} \overline{\delta}} \phi_{,\overline{z} l
\alpha}
 \geq {2\over n} S^2. \nonumber
\end{eqnarray}
\end{theo}

\begin{proof} According to Lemma \ref{local:regularity}, we have
\[
 \overline{\partial_z} {{\partial \eta^{\alpha}}\over {\partial \overline{w_i}}}
= -{{\partial \eta^{\overline{\beta}}}\over {\partial
\overline{w_i}}}  {{\partial \eta^{\alpha}}\over {\partial
\overline{w_{\beta}}}}. \] Thus
\[
\begin{array}{lcl}
\partial_z \overline{\partial_z} {{\partial \eta^{\alpha}}\over {\partial \overline{w_i}}}
& = & - \partial_z \left({{\partial \eta^{\overline{\beta}}}\over
{\partial \overline{w_i}}} \right)  {{\partial \eta^{\alpha}}\over
{\partial \overline{w_{\beta}}}}
- {{\partial \eta^{\overline{\beta}}}\over {\partial \overline{w_i}}} \partial_z  {{\partial \eta^{\alpha}}\over {\partial \overline{w_{\beta}}}}\\
& = & {{\partial \eta^{p}}\over {\partial \overline{w_i}}}
{{\partial \eta^{\overline{\beta}}}\over {\partial w_p}}
{{\partial \eta^{\alpha}}\over {\partial \overline{w_{\beta}}}}
-{{\partial \eta^{\overline{\beta}}}\over {\partial
\overline{w_i}}} \partial_z  {{\partial \eta^{\alpha}}\over
{\partial \overline{w_{\beta}}}}.
\end{array}
\]
 We have
\[
\begin{array}{lcl}
 & & \Delta_z \left( {{\partial \eta^{\alpha}}\over {\partial \overline{w_i}}} {{\partial \eta^{\overline{i}}}\over {\partial
 w_{\alpha}}}\right)\\
& = & \partial_z \left( {{\partial \eta^{\alpha}}\over {\partial
\overline{w_i}}}\right) \overline{\partial_z} \left( {{\partial
\eta^{\overline{i}}}\over {\partial w_{\alpha}}}  \right) +
\overline{\partial_z }\left( {{\partial \eta^{\alpha}}\over
{\partial \overline{w_i}}}\right)
\partial_z \left( {{\partial \eta^{\overline{i}}}\over {\partial w_{\alpha}}}  \right)\\
& & \qquad \qquad \qquad + \triangle_z \left( {{\partial
\eta^{\alpha}}\over {\partial \overline{w_i}}}\right)
 \left( {{\partial \eta^{\overline{i}}}\over {\partial w_{\alpha}}}  \right)
+ \left( {{\partial \eta^{\alpha}}\over {\partial
\overline{w_i}}}\right)
\triangle_z \left( {{\partial \eta^{\overline{i}}}\over {\partial w_{\alpha}}}  \right)\\
& = & \partial_z \left( {{\partial \eta^{\alpha}}\over {\partial
\overline{w_i}}}\right) \overline{\partial_z} \left( {{\partial
\eta^{\overline{i}}}\over {\partial w_{\alpha}}}  \right) +
{{\partial \eta^{\overline{\beta}}}\over {\partial
\overline{w_i}}}  {{\partial \eta^{\alpha}}\over {\partial
\overline{w_{\beta}}}}
{{\partial \eta^{{p}}}\over {\partial w_{\alpha}}}  {{\partial \eta^{\overline{i}}}\over {\partial {w_{p}}}}\\
& & \qquad \qquad \qquad + \left( {{\partial \eta^{p}}\over {\partial \overline{w_i}}} {{\partial \eta^{\overline{\beta}}}\over {\partial w_p}}
 {{\partial \eta^{\alpha}}\over {\partial \overline{w_{\beta}}}} -{{\partial \eta^{\overline{\beta}}}\over {\partial \overline{w_i}}} \partial_z
  {{\partial \eta^{\alpha}}\over {\partial \overline{w_{\beta}}}}\right) {{\partial \eta^{\overline{i}}}\over {\partial w_{\alpha}}}\\
& &  \qquad \qquad \qquad +{{\partial \eta^{\alpha}}\over {\partial \overline{w_{i}}}}
\overline{\left({{\partial \eta^{p}}\over {\partial \overline{w_{\alpha}}}} {{\partial \eta^{\overline{\beta}}}\over {\partial w_p}}
{{\partial \eta^{i}}\over {\partial \overline{w_{\beta}}}} -{{\partial \eta^{\overline{\beta}}}\over {\partial \overline{w_{\alpha}}}} \partial_z
 {{\partial \eta^{i}}\over {\partial \overline{w_{\beta}}}} \right)}\\
& = & \partial_z \left( {{\partial \eta^{\alpha}}\over {\partial
\overline{w_i}}}\right) \overline{\partial_z} \left( {{\partial
\eta^{\overline{i}}}\over {\partial w_{\alpha}}}  \right) +
{{\partial \eta^{\overline{\beta}}}\over {\partial
\overline{w_i}}}  {{\partial \eta^{\alpha}}\over {\partial
\overline{w_{\beta}}}} {{\partial \eta^{{p}}}\over {\partial
w_{\alpha}}}  {{\partial \eta^{\overline{i}}}\over {\partial
{w_{p}}}} - \partial_z  \left({{\partial \eta^{\alpha}}\over
{\partial \overline{w_{\beta}}}}\right)\; {{\partial
\eta^{\overline{\beta}}}\over {\partial \overline{w_i}}}\;
 {{\partial \eta^{\overline{i}}}\over {\partial w_{\alpha}}}\\
& & - {{\partial \eta^{\alpha}}\over {\partial \overline{w_{i}}}}
\;{{\partial \eta^{\beta}}\over {\partial w_{\alpha}}}
\overline{\partial_z}  {{\partial \eta^{\overline{i}}}\over
{\partial w_{\beta}}} + {{\partial \eta^{p}}\over {\partial
\overline{w_i}}} {{\partial \eta^{\overline{\beta}}}\over
{\partial w_p}} {{\partial \eta^{\alpha}}\over {\partial
\overline{w_{\beta}}}} \;{{\partial \eta^{\overline{i}}}\over
{\partial w_{\alpha}}} + {{\partial \eta^{\alpha}}\over {\partial
\overline{w_{i}}}} {{\partial \eta^{\overline{p}}}\over {\partial
w_{\alpha}}} {{\partial \eta^{\beta}}\over {\partial
\overline{w_p}}}
{{\partial \eta^{\overline{i}}}\over {\partial w_{\beta}}}\\

& = & \partial_z \left( {{\partial \eta^{\alpha}}\over {\partial
\overline{w_i}}}\right) \overline{\partial_z} \left( {{\partial
\eta^{\overline{i}}}\over {\partial w_{\alpha}}}  \right) +
{{\partial \eta^{\overline{\beta}}}\over {\partial
\overline{w_i}}}  {{\partial \eta^{\alpha}}\over {\partial
\overline{w_{\beta}}}} {{\partial \eta^{{p}}}\over {\partial
w_{\alpha}}}  {{\partial \eta^{\overline{i}}}\over {\partial
{w_{p}}}} - \partial_z  \left({{\partial \eta^{\alpha}}\over
{\partial \overline{w_{i}}}}\right)\; {{\partial
\eta^{\overline{i}}}\over {\partial \overline{w_p}}}\;
 {{\partial \eta^{\overline{p}}}\over {\partial w_{\alpha}}}\\
& & - {{\partial \eta^{\beta}}\over {\partial \overline{w_{i}}}}
\;{{\partial \eta^{\alpha}}\over {\partial w_{\beta}}}
\overline{\partial_z}  {{\partial \eta^{\overline{i}}}\over
{\partial w_{\alpha}}} + {{\partial \eta^{p}}\over {\partial
\overline{w_i}}} {{\partial \eta^{\overline{\beta}}}\over
{\partial w_p}} {{\partial \eta^{\alpha}}\over {\partial
\overline{w_{\beta}}}} \;{{\partial \eta^{\overline{i}}}\over
{\partial w_{\alpha}}} + {{\partial \eta^{\alpha}}\over {\partial
\overline{w_{i}}}} {{\partial \eta^{\overline{p}}}\over {\partial
w_{\alpha}}} {{\partial \eta^{\beta}}\over {\partial
\overline{w_p}}}
{{\partial \eta^{\overline{i}}}\over {\partial w_{\beta}}}\\
& = & \left(\partial_z {{\partial \eta^{\alpha}}\over {\partial
\overline{w_i}}} - {{\partial \eta^{\beta}}\over {\partial
\overline{w_{i}}}} \;{{\partial \eta^{\alpha}}\over {\partial
w_{\beta}}}\right) \left( \overline{\partial_z}  {{\partial
\eta^{\overline{i}}}\over {\partial w_{\alpha}}} - {{\partial
\eta^{\overline{i}}}\over {\partial \overline{w_p}}}\;
 {{\partial \eta^{\overline{p}}}\over {\partial w_{\alpha}}}\right) - {{\partial \eta^{\beta}}\over {\partial \overline{w_{i}}}} \;{{\partial \eta^{\alpha}}\over {\partial w_{\beta}}}\;{{\partial \eta^{\overline{i}}}\over {\partial \overline{w_p}}}\;
 {{\partial \eta^{\overline{p}}}\over {\partial w_{\alpha}}} \\
& & + {{\partial \eta^{\overline{\beta}}}\over {\partial
\overline{w_i}}}  {{\partial \eta^{\alpha}}\over {\partial
\overline{w_{\beta}}}} {{\partial \eta^{{p}}}\over {\partial
w_{\alpha}}}  {{\partial \eta^{\overline{i}}}\over {\partial
{w_{p}}}} + {{\partial \eta^{p}}\over {\partial \overline{w_i}}}
{{\partial \eta^{\overline{\beta}}}\over {\partial w_p}}
{{\partial \eta^{\alpha}}\over {\partial \overline{w_{\beta}}}}
\;{{\partial \eta^{\overline{i}}}\over {\partial w_{\alpha}}} +
{{\partial \eta^{\alpha}}\over {\partial \overline{w_{i}}}}
{{\partial \eta^{\overline{p}}}\over {\partial w_{\alpha}}} {{\partial \eta^{\beta}}\over {\partial \overline{w_p}}} {{\partial \eta^{\overline{i}}}\over {\partial w_{\beta}}}\\

& = & \left(\partial_z {{\partial \eta^{\alpha}}\over {\partial
\overline{w_i}}} - {{\partial \eta^{\beta}}\over {\partial
\overline{w_{i}}}} \;{{\partial \eta^{\alpha}}\over {\partial
w_{\beta}}}\right) \left( \overline{\partial_z}  {{\partial
\eta^{\overline{i}}}\over {\partial w_{\alpha}}} - {{\partial
\eta^{\overline{i}}}\over {\partial \overline{w_p}}}\;
 {{\partial \eta^{\overline{p}}}\over {\partial w_{\alpha}}}\right)
- {{\partial \eta^{\beta}}\over {\partial \overline{w_{i}}}}
\;{{\partial \eta^{\alpha}}\over {\partial w_{\beta}}}\;
{{\partial \eta^{\overline{p}}}\over {\partial w_{\alpha}}}
\;{{\partial \eta^{\overline{i}}}\over {\partial
\overline{w_p}}}\;
\\
& & +   {{\partial \eta^{\alpha}}\over {\partial
\overline{w_{\beta}}}} {{\partial \eta^{{p}}}\over {\partial
w_{\alpha}}}  {{\partial \eta^{\overline{i}}}\over {\partial
{w_{p}}}} {{\partial \eta^{\overline{\beta}}}\over {\partial
\overline{w_i}}} + {{\partial \eta^{p}}\over {\partial
\overline{w_i}}} {{\partial \eta^{\overline{\beta}}}\over
{\partial w_p}} {{\partial \eta^{\alpha}}\over {\partial
\overline{w_{\beta}}}} \;{{\partial \eta^{\overline{i}}}\over
{\partial w_{\alpha}}} + {{\partial \eta^{\alpha}}\over {\partial
\overline{w_{i}}}}
{{\partial \eta^{\overline{p}}}\over {\partial w_{\alpha}}} {{\partial \eta^{\beta}}\over {\partial \overline{w_p}}} {{\partial \eta^{\overline{i}}}\over {\partial w_{\beta}}}\\

& = & \left(\partial_z {{\partial \eta^{\alpha}}\over {\partial
\overline{w_i}}} - {{\partial \eta^{\beta}}\over {\partial
\overline{w_{i}}}} \;{{\partial \eta^{\alpha}}\over {\partial
w_{\beta}}}\right) \left( \overline{\partial_z}  {{\partial
\eta^{\overline{i}}}\over {\partial w_{\alpha}}} - {{\partial
\eta^{\overline{i}}}\over {\partial \overline{w_p}}}\;
 {{\partial \eta^{\overline{p}}}\over {\partial w_{\alpha}}}\right) + 2 {{\partial \eta^{\alpha}}\over {\partial \overline{w_{i}}}}
{{\partial \eta^{\overline{p}}}\over {\partial w_{\alpha}}}
{{\partial \eta^{\beta}}\over {\partial \overline{w_p}}}
{{\partial \eta^{\overline{i}}}\over {\partial w_{\beta}}}.
 \end{array}
\]
The last equality holds because in the line above the last
equation, the 2rd term and the 3rd cancel each other, while the
3rd and 4th term are the same.\\

  Note that $\partial_z {{\partial \eta^{\alpha}}\over {\partial \overline{w_i}}} -
{{\partial \eta^{\beta}}\over {\partial \overline{w_{i}}}}
\;{{\partial \eta^{\alpha}}\over {\partial w_{\beta}}} $ is a
tensor since
\[
\begin{array}{lcl} \partial_z {{\partial \eta^{\alpha}}\over {\partial \overline{w_i}}} -
{{\partial \eta^{\beta}}\over {\partial \overline{w_{i}}}}
\;{{\partial \eta^{\alpha}}\over {\partial w_{\beta}}} & = & {{\p
}\over {\partial z}}  {{\partial \eta^{\alpha}}\over {\partial
\overline{w_i}}}  + \eta^r {{\p} \over {\p w_r}} {{\partial
\eta^{\alpha}}\over {\partial \overline{w_i}}} -{{\partial
\eta^{\beta}}\over {\partial \overline{w_{i}}}} \left({{\partial
\eta^{\alpha}}\over {\partial w_{\beta}}} +  \eta^r
\Gamma^{\alpha}_{r \;\beta}- \eta^r \Gamma^{\alpha}_{r
\;\beta}\right)\\
& = & {{\p }\over {\partial z}}  {{\partial \eta^{\alpha}}\over
{\partial \overline{w_i}}} + \eta^r \left( {{\p} \over {\p w_r}}
{{\partial \eta^{\alpha}}\over {\partial \overline{w_i}}} +
{{\partial \eta^{\beta}}\over {\partial \overline{w_{i}}}}
\Gamma^{\alpha}_{r \;\beta}\right)   -{{\partial
\eta^{\beta}}\over {\partial \overline{w_{i}}}}
{\eta^{\alpha}}_{,\beta} \\
& = & {{\p }\over {\partial z}}  {{\partial \eta^{\alpha}}\over
{\partial \overline{w_i}}} + \eta^r  {\eta^{\alpha}}_{,\b i r}
-{{\partial \eta^{\beta}}\over {\partial \overline{w_{i}}}}
{\eta^{\alpha}}_{,\beta}.
\end{array}
\]
Moreover, this is a (0,2) symmetric tensor since
\[
\begin{array}{lcl} \partial_z {{\partial \eta^{\alpha}}\over {\partial \overline{w_i}}} -
{{\partial \eta^{\beta}}\over {\partial \overline{w_{i}}}}
\;{{\partial \eta^{\alpha}}\over {\partial w_{\beta}}} & = &
\partial_z \left(- g^{\alpha \overline{\beta}} \;\partial_{z}
{{\partial^2 \Phi}\over{\overline{\partial w_i} \overline{\partial
w_\beta}}}\right) -
{{\partial \eta^{\beta}}\over {\partial \overline{w_{i}}}} \;{{\partial \eta^{\alpha}}\over {\partial w_{\beta}}} \\
& = & g^{\alpha \overline{p}} \partial_z g_{\overline{p} q} g^{q \overline{\beta}} {{\partial^2 \Phi}\over{\overline{\partial w_i} \overline{\partial w_\beta}}}  - g^{\alpha \overline{\beta}} \;\left(\partial_{z}\right)^2 {{\partial^2 \Phi}\over{\overline{\partial w_i} \overline{\partial w_\beta}}} - {{\partial \eta^{\beta}}\over {\partial \overline{w_{i}}}} \;{{\partial \eta^{\alpha}}\over {\partial w_{\beta}}} \\
& = & {{\partial \eta^{q}}\over {\partial \overline{w_{i}}}} \;{{\partial \eta^{\alpha}}\over {\partial w_{q}}} - g^{\alpha \overline{\beta}} \;\left(\partial_{z}\right)^2 {{\partial^2 \Phi}\over{\overline{\partial w_i} \overline{\partial w_\beta}}} - {{\partial \eta^{\beta}}\over {\partial \overline{w_{i}}}} \;{{\partial \eta^{\alpha}}\over {\partial w_{\beta}}} \\
& = & - g^{\alpha \overline{\beta}} \;\left(\partial_{z}\right)^2
{{\partial^2 \Phi}\over{\overline{\partial w_i} \overline{\partial
w_\beta}}} \\
& = & - g^{\alpha \overline{\beta}} \;\partial_{z} \left( {{\p
\phi}\over {\p z}} \right)_{,\b i \b \beta}.
\end{array}
\]
  Thus, the first term in the equation  (\ref{eq:C3estimate}) can be
changed into
\[
\begin{array}{lcl} & &
\left(\partial_z {{\partial \eta^{\alpha}}\over {\partial
\overline{w_i}}} - {{\partial \eta^{\beta}}\over {\partial
\overline{w_{i}}}} \;{{\partial \eta^{\alpha}}\over {\partial
w_{\beta}}}\right) \left( \overline{\partial_z}  {{\partial
\eta^{\overline{i}}}\over {\partial w_{\alpha}}} - {{\partial
\eta^{\overline{i}}}\over {\partial \overline{w_p}}}\;
 {{\partial \eta^{\overline{p}}}\over {\partial
 w_{\alpha}}}\right)\\ & & \qquad \qquad
= g^{\alpha \overline{\beta}} \;\left(\partial_{z}\right)^2
{{\partial^2 \Phi}\over{\overline{\partial w_i} \overline{\partial
w_\beta}}} \; g^{\overline{i} p }
\;\left(\overline{\partial_{z}}\right)^2 {{\partial^2
\Phi}\over{\partial w_{\alpha} \partial w_p}} \geq 0.
\end{array}
\]
The lemma is then proved.
\end{proof}
This theorem should be compared with Calabi's third derivative
estimate for the non-degenerate Monge-Ampere equation.
   Following a  result of R. Osserman\cite{Osserman57} (later generalized
   by E. Calabi\cite{Calabi57}), we have
   \begin{prop} \label{local:interiorboundofcurvature} Let $d$ denote the Euclidean distance to the boundary of
   $\Sigma.\;$ Then, there exists a uniform constant $C$ such that
   the trace of the curvature has the following interior estimate
   \[
      S(z\;,x) = - \displaystyle \sum_{\alpha=1}^n\;F^{\alpha}_\alpha  \leq {{C}\over {d(z,\p \Sigma)^2}}.
   \]
   \end{prop}
This estimate plays a crucial in obtaining a priori estimate for super regular discs.

\section{Compactness of holomorphic disks}
In this section, we continue our study of the HCMA equation from
the point of view initiated in the previous section. Namely, in
the foliation by holomorphic disks of $\Sigma\times M$,  we study
the family of restricted $TM$ bundles equipped with the varying
Hermitian metric $g_\phi.\;$ We introduce two geometrical
quantities which can be used to study the compactness problem of
holomorphic disks. For any ${\phi_0}:\p \Sigma \rightarrow \cal H$
and any holomorphic disk $f \in {\cal M}_{\phi_0},$ we define its
area as
\begin{equation}
   A(f) = \displaystyle\;\int_\Sigma\; \left( {\sqrt{-1} \over 2 } d\,z \wedge d\,\b z + f^* \omega_0\right).
   \label{eq:area}
\end{equation}
Note that this is the area of $\pi\circ f(\Sigma)$ in
$\Sigma\times M, $ not the area of $f(\Sigma)\in \Sigma\times
\cW_M.\;$  When no confusion is possible, we will not distinguish
between $f,\;\pi\circ f,\;$ or even $\pi_2\circ \pi\circ f.\;$
Similarly, we also define the {\bf Capacity} for any super regular
disk $f$ in ${\cal M}_{\phi_0}$ by
\begin{equation}Cap(f) =\displaystyle\; \int_\Sigma\;  {{f^*\omega_{\phi_0(z_0)}^n}\over {f^*\omega_\phi^n}}\; {\sqrt{-1}\over 2}
\;d\,z \;d\,\b z.\label{eq:capacity}
\end{equation}
for some $z_0\in \p \Sigma.\;$ For simplicity, we fix $z_0$ in
this section and assume without loss of generality that
$\phi_0(z_0,\cdot ) = 0.\;$ Under this assumption, $\omega_{\phi_0(z_0,\cdot)}
=\omega_0.\;$
Obviously, a non-super regular disk has infinite capacity.   \\

Let us set up some notations first. We will fix a positive number
$\alpha \in (0,1)$ in this section.  For any $\epsilon,\delta,
\Lambda, $  let  ${\phi_0}$ be  any map from $\p\Sigma $
to $ \cal H$ which satisfies
\begin{equation}
\omega_{\phi_0} \geq \delta \omega_0,\qquad
\|{\phi_0}\|_{C^{2,\alpha}(\p \Sigma\times M)} \leq \Lambda.
\label{eq:boundpsi}
\end{equation}
We further define ${\cal C}(\delta,\Lambda)$ as the space of
all embedded holomorphic disks in $\Sigma\times \cW_M$ with
vanishing normal {\it Maslov} indice, whose boundary lies in $\bar
\Lambda_{\phi_0}. \;$
The space ${\cal C}(\delta,\Lambda, L_0)$ is a subset of ${\cal
C}(\delta,\Lambda)$ such that for each disk  $ f\in  {\cal
C}(\delta,\Lambda, L_0), $ the corresponding compatible solution $\phi$ of HCMA (1.1)
is smooth locally and
\begin{equation}
\mid\phi\mid_{C^{1,1}} \leq L_0  \label{eq:boundvarphi}
\end{equation}
hold in a small neighborhood of $f(\Sigma).\;$
 Define
\[
 {\cal C}(\delta, \Lambda,L_0, L_1) = \{ f\in {\cal C}(\delta,\Lambda,L_0)\mid A(f) \leq L_1\}\]
and
\[
 {\cal C}(\delta, \Lambda,L_0, L_1,L_2) = \{ f\in {\cal C}(\delta,\Lambda,L_0)\mid A(f) \leq L_1, Cap(f) \leq L_2\}.\]

 In this
Section, we will prove
\begin{theo} \label{comp:nobubble} The space ${\cal C}(\delta, \Lambda,L_0, L_1)$  is compact in   ${\cal C}(\delta,\Lambda,L_0).\;$
\end{theo}

\begin{theo} \label{comp:superregularlimit} The space ${\cal C}(\delta, \Lambda,L_0, L_1, L_2)$ is compact.
\end{theo}

Suppose $\cF_{\phi_0}$ is an almost super regular foliation, we
want to identify $\cU_{\phi_0}$ with an open and dense subset of
$M$ via evaluation map:
\[
\begin{array}{llll} \sharp: & \Sigma\times {\cal M_{\phi_0}} &  \rightarrow & \Sigma \times M 
 \\ & (z, f) & \rightarrow  & \pi\circ ev(z, f)  
   \end{array}
\]
 Then, $\sharp$ is
invertible on $\Sigma\times \cU_{\phi_0}.\;$  Define $F =\pi \circ
ev \circ \sharp^{-1}.\;$ We can identify $\cU_{\phi_0}$ with some
open dense subset $\pi\circ ev (z_0, \cU_{\phi_0})$ of $M.\;$ We
will use this point of view from time to time in this section.

\begin{theo} \label{comp:boundpsi}  If ${\phi_0}$ satisfies inequality (\ref{eq:boundpsi}) and if $\cF_{\phi_0}$
is an almost super regular foliation, then there exist two
constants $L_0, L_1$ which depend only on $\delta,\Lambda$ such
that
\[
\cU_{\phi_0} \subset {\cal C} (\delta, \Lambda, L_0, L_1).
\]
Moreover,
\[
 \displaystyle\; \int_{ f\in \cU_{\phi_0}} \; Cap(f)\;\omega_{\phi_0(z_0,\cdot)}^n = \displaystyle\; \int_{ f\in \cU_{\phi_0}} \; Cap(f)\;\omega_0^n \leq C.
\]
\end{theo}
These three theorems can be used to derive a compactness theorem of sequences of
almost super regular foliations.

\begin{theo} \label{comp:foliationlimit} Suppose $\{{\phi_0}^{(m)}, m\in \cN\}$ is a sequence of functions
in $C^{\infty}(\p\Sigma, \cH)$ which satisfies the uniform bound
of (\ref{eq:boundpsi}) and converges
 to ${\phi_0}^{(\infty)} \in C^{\infty}(\p\Sigma, \cH)\;$ in the $C^{2,\alpha}(\p \Sigma\times M)\;$
 norm.
 Suppose  that $\{ \cF_{\phi_0^{(m)}}, m\in \cN\}$ is a sequence of
almost super regular foliations, while $\{\phi^{(m)}\}$ is the
corresponding sequence of almost smooth solutions with Drichelet
boundary value $\{{\phi_0}^{(m)}\}.\;$ Passing to a subsequence if
necessary, $\cF_{{\phi_0}^{(m)}}$ converges to a partially smooth
foliation $\cF_{{\phi_0}^{(\infty)}}.\;$ In particular, at least
one component of ${\cal M}_{{\phi_0}^{(\infty)}}$ contains at
least one super regular disk. Moreover, $\phi^{(m)}$ converges in
the  weak $C^{1,1}$ norm to $\phi^{(\infty)}$ such that
$\omega_{\phi^{(\infty)}}^n$ is a continuous volume form on
$\Sigma^0 \times M.$
\end{theo}

\subsection{Proof of Theorem \ref{comp:boundpsi}}
\begin{proof} Suppose that $\phi$ is the corresponding almost smooth solution of the HCMA equation
\ref{eq:hcma0} with the prescribed boundary map ${\phi_0}:  \p
\Sigma \rightarrow \cal H.\;$ By Theroem 1.2 \cite{chen991}, there
is a uniform constant $C(\delta, \Lambda)$ such that
\[
   \mid \p \b \p \phi\mid \leq C(\delta, \Lambda).
\]
It is clear that for any super regular disk $f \in \cU_{\phi_0},$
any local compatible solution to HCMA equation (1.1)  must agree with $\phi$ in any small open and
saturated neighborhood of $\pi\circ f(\Sigma) \subset \Sigma\times
M.\;$ Next, we want to show that each disk must also have a
uniform upper bound  on its area.

\begin{lem} For any regular disk $f: (\Sigma,\p \Sigma)\rightarrow (\Sigma\times \cW_M,\b \Lambda_{\phi_0})$
 where $\phi$ is the corresponding ``almost smooth" solution of equation \ref{eq:hcma0}. There is a uniform constant $L_1$ which depends
 only on $\delta, \Lambda$ such that $A(f) \leq L_1\;$ holds
uniformly.
\end{lem}
\begin{proof}  Recalled that the leaf vector field $X$ (cf. eq. (\ref{loca:leafvector})) along this disk in $\Sigma\times M$
can be expressed as
\begin{equation}\label{comp:diskvectorfield}
X = \displaystyle \sum_{\alpha=1}^n\; \eta^{\alpha} {\p \over {\p
w_\alpha} }= - \displaystyle \sum_{\alpha=1}^n\; g_\phi^{\alpha
\overline{\beta}}{{\partial^2 \phi}\over{\partial z \partial w_{\b
\beta}}} {\p \over {\p w_\alpha} }.
\end{equation}

According to Corollary \ref{local:deltazpotential}, we have
\[ {\p\over {\p z}} {\p \over {\p \b z}} (\phi \circ f) =  -
{g_{0}}_{\alpha \overline{\beta}} \eta^{\alpha}
  \eta^{\overline{\beta}}. \]
 Using this holomorphic map $f$, the pull back of the fixed product metric
metric on $\Sigma\times M$ to  $\Sigma$ is:
\[\begin{array}{lcl}
  f^{*}(g_0 + |d\, z|^2) & = & |{\partial  \over {\p z}} + X |_{g_0}^2  |d\,
  z|^2\\
  & = & \left(1 + {g_{0}}_{\alpha \overline{\beta}} \eta^{\alpha}
  \eta^{\overline{\beta}} \right) |d\, z|^2.
  \end{array}
\]
  Thus, the area of any disk is :
  \[\begin{array} {lcl}
    \int_{f(\Sigma)} 1 & = & \int_{\Sigma} f^{*}(g_0 + |d\,
    z|^2)\\ & = & \displaystyle \int_\Sigma \; \left(1 + {g_{0}}_{\alpha \overline{\beta}} \eta^{\alpha}
  \eta^{\overline{\beta}}  \right) |d\,z|^2\\
    & = &  |\Sigma| - \int_{\Sigma}\; \partial_z
    \partial_{\b z} \; \phi \; |d\, z|^2\\ & =&  |\Sigma| - \int_{\Sigma}\; {{\partial^2 }\over {\partial
    z \partial \overline{z}}}
     \;  (\phi\circ f)  \; |d\, z|^2  \\
    & = & |\Sigma| + \int_{\partial \Sigma}
    {\partial \over {\partial z}} (\phi \circ f) \; \cdot {\bf n}_{\Sigma}
    \\ & = &  |\Sigma| + \int_{\partial \Sigma}
     \partial_z \phi  \mid_{\p \Sigma} \cdot {\bf n}_{\Sigma}  \; |d\, z|^2\\
     & = & |\tilde{\Sigma}| + \int_{\partial {\Sigma}}
   \left({\partial \over {\partial
  z}} + \eta^{\alpha} {\partial \over {\partial w_{\alpha}}} \right)\phi \mid_{\p \Sigma} \cdot
  {\bf n}_{{\Sigma}}  \; |d\, z|^2\\
  & = &|{\Sigma}| + \int_{\partial {\Sigma}}
   \left({ {\partial \phi} \over {\partial
  z}}- g_\phi^{\alpha \b \beta} {{\partial^2 \phi}\over {\partial z \partial w_{\b \beta}}} {\partial \phi \over {\partial w_{\alpha}}} \right) \mid_{\p \Sigma} \cdot
 {\bf n}_{{\Sigma}}  \; |d\, z|^2,
    \end{array}
  \]
  where ${\bf n}_{{\Sigma}}$ represents the normal direction on
  the boundary of ${\Sigma}.\;$  On $\p \Sigma, $ we have
  \[
\begin{array}{lcl}
      g_{\phi, \alpha \b \beta}(z,\cdot)&  =&  {g_{0}}_{\alpha \b \beta}(\;z\;,\cdot) + \sqrt{-1}
      \p \b \p \phi(\;z\;,\cdot)  \\
       & = &  {g_{0}}_{\alpha \b \beta}(\;z\;,\cdot) + \sqrt{-1}
      \p \b \p {\phi_0}(\;z\;,\cdot) \geq
      \delta \;{g_{0}}_{\alpha \b \beta}.\end{array}
  \]
    Thus,
  \[
  \begin{array}{lcl}
\int_{\Sigma} 1 & = & |{\Sigma}| + \int_{\partial {\Sigma}}
   \left({ {\partial \phi} \over {\partial
  z}}- g_\phi^{ \alpha \b \beta} {{\partial^2 \phi }\over {\partial z \partial \b \beta}}
   {\partial \phi \over {\partial w_{\alpha}}} \right) \mid_{\p \Sigma} \cdot
  {\bf n}_{{\Sigma}}  \; |d\, z|^2 \\
  &\leq & L_1(\delta, \Lambda).
  \end{array}
  \]
\end{proof}
Now we return to the proof of Theorem \ref{comp:boundpsi}.   Note
that
   \[  F^{*} \omega_{\phi^{}}^n = \omega_{\phi_0(z_0)}^n =  \omega_0^n,
   \]
   and
   \[\begin{array}{lcl}
      F^{*} \omega_{\phi_0(z_0)}^n &  = &  F^*  \left({{\omega_{0}^n} \over {\omega_{\phi^{}}^n} }\; \omega_{\phi^{}}^n \right)
                     \\ & = & \left({{\omega_{0}^n} \over {\omega_{\phi}^n}
                     }\right)\circ F \cdot \omega_{0}^n.
                     \end{array}
   \]
   Thus
   \[
   \begin{array}{lcl} \int_M  \int_\Sigma \left({{\omega_{0}^n} \over
   {\omega_{\phi^{}}^n}}\right)\circ F^{} \;{\sqrt{-1}\over 2} d\,z\wedge d
   \b z z\; \omega_{0}^n & = & \int_M  \int_\Sigma \left({{\omega_{0}^n}
\over {\omega_{\phi^{(m)}}^n}
                     }\right)\circ F^{} \; \omega_{0}^n\;\;{\sqrt{-1}\over 2}\;  d\,\;z\wedge d\,\b z\;\;\\
                     & = & \int_M  \int_\Sigma F^{*}\omega_{0}^n \;{\sqrt{-1}\over 2}\;  d\,\;z\wedge d\,\b z\;\\
                     & = & \int_\Sigma \int_M \omega_{0}^n \;{\sqrt{-1}\over 2}\; d\,z\wedge d\,\b z = C.
   \end{array}
   \]
  In other words, we have
\[
\displaystyle \int_M\;{\rm Cap}(\Sigma_x)\; d\,x = \displaystyle
\int_M\; \int_\Sigma \left({{\omega_{0}^n} \over {\omega_{\phi}^n}
                     }\right)\circ F \;{\sqrt{-1}\over 2}\;  d\,\;z\;
                     d\,\wedge {\b z\;}\;d\,x \leq C.
\]
  where $C$ is a topological constant.  This concludes our proof of Theorem \ref{comp:boundpsi}.
  \end{proof}

\subsection{No vertical bubble--Proof of Theorem \ref{comp:nobubble}}
We want to re-phrase Theorem \ref{comp:nobubble}.
\begin{theo} \label{comp:nobubble1} For any sequence of super regular disks (in a sequence of almost
super regular foliations $\{({\phi_0}^{(m)},\cU_{{\phi_0}^{(m)}}), m \in \NN\})$,
there is no vertical bubble in the limit provided that the
corresponding sequence of boundary maps ${\phi_0}^{(m)}$ converges
in the $C^{2,\alpha}(\p\Sigma\times M)$ norm to some potential function $\phi_0^{(\infty)} $ in the completion
of $\cH$ by $ C^{2,\alpha}(\p\Sigma\times M)$ norm.
\end{theo}
\begin{proof}  Suppose
  \[
\begin{array}{llll}
f^{(m)}:& \Sigma & \mapsto &  \Sigma \times \cW_M\\
  &   z            & \mapsto  &  (\;z\;, ,f^{(m)}(z),\zeta^{(m)}(f^{(m)} (z) ) ), \qquad m=1,2,\cdots\infty
  \end{array}
\]
is a sequence  of super regular disks in $\{({\phi_0}^{(m)},\cU_{{\phi_0}^{(m)}}), m \in \NN\}.\;$ According
to  Theorem 5.0.16, there exists two constants $L_0, L_1$ such that all $\{f^{(m)}, m \in \NN\} $
lies in ${\cal C}(\delta,\Lambda,
L_0,L_1.\;$ Recalled that
\[
{{\partial f^{(m) \alpha} }\over{\partial z}} = - g^{(m)\; \alpha
\overline{\beta}}{{\partial^2 \phi^{(m)}}\over{\partial z
\overline{\partial w_\beta}}}.
\]
Here $\zeta$ is the corresponding fibre component of
$f^{(m)}(\Sigma)$ in $\cW_M.\;$ In a local coordinate, we write
\[
 \zeta^{(m)\; i} (z, x) = {{\p (\phi^{(m)}+\rho)} \over {\p w_{
 i}}}(z, x),\qquad \forall \;i = 1, 2,\cdots n
\]
where $\rho$ is a local K\"ahler potential for the given form
$\omega.\;$ Note that in a uniform size neighborhood of $\pi\circ
f^{(m)}(\Sigma)\subset \Sigma\times M$,  we have
\[
  \mid\phi^{(m)}\mid_{C^{1,1}} \leq L_0.
\]
In particular, there exists a uniform constant $C (\delta,
\Lambda)$ such that
\begin{equation}
   \displaystyle \max_{\p \Sigma} \mid {{\p f^{(m)}}\over {\p z}}\mid_{g_0} \leq C. \label{comp: uniform bound in the boundary}
\end{equation}
  Consequently, all
bubble points\footnote{A point $\{(z_m, x_m), m\in \NN\}$  is
called bubble point if a) $\mid {{\p f^{(m)}}\over {\p
z}}(z_m, x_m)\mid_{g_0} \rightarrow \infty$ and b) it is a global maximal of
$\mid {{\p f^{(m)}}\over {\p z}}\mid_{g_0}$ in $\Sigma\times M.\;$} that occur must
occur in the interior of $\Sigma \times M$ (although the bubble
may travel to some boundary
point in the limit.).  We want to show that no such bubble point exists; which in turns implies the present theorem.\\

Suppose bubbling does occur and
  there is a sequence of points  $ (z_m, x_m)$ such that
\[
\epsilon_m^{-1} =  \displaystyle \max_{\Sigma\times M} \mid {{\p
f^{(m)}}\over {\p z}}\mid_{g_0}  =\mid {{\p f^{(m)}}\over {\p
z}}(z_m, x_m)\mid_{g_0} \rightarrow \infty. \]

W.l.o.g,   we may assume that $\displaystyle \lim_{m\rightarrow
\infty} \;z_m = \;z_\infty \in  \Sigma\;$ and $\displaystyle
\lim_{m \rightarrow \infty}\; x_m = x_\infty \in M.\;$  We want to
argue that  there exists a point  $z_m' \in B_{2
\sqrt{\epsilon_m}}(z_m) \bigcap \Sigma$ such that the area of
$B_{\epsilon_m}(z_m')\geq c_0 $ for some uniform constant
$c_0>0.\;$  If the area functional were the area of a holomorphic
disk in $\Sigma\times \cW_M$, then this follows from standard
literatures in this
 direction.  In our setting, this is still true.  This in turn implies that there are at most finite number
 of bubbles.  We give a brief explanation here and leave interested readers to fill in the details.\\

 Set $d_m = d(z_m, \p\Sigma).\;$ If $d_m > {1\over 2}\epsilon_m$, then an easy calculation implies

 \begin{equation}
 A\left(f^{(m)}, B_{\epsilon_m} (z_m) \bigcap \Sigma\right)\geq c_0.
 \label{comp:nontrivialarea}
 \end{equation}
   However, we need to establish inequality (\ref{comp:nontrivialarea}) even when ${d_m\over \epsilon_m}\rightarrow 0.\;$  In such a case,
   there must exists another point $z_m'\in \p B_{1-d_m-\epsilon_m}(O) \bigcap B_{\sqrt{\epsilon_m}}(z_m)$
 such that
 \[
  \mid {{\p f^{(m)}}\over {\p z}}(z_m', x_m')\mid_{g_0} > {\epsilon_m\over 2}, \qquad x_m' = f^{(m)}(z_m').
 \]
 This can be proved by using inequality (\ref{comp: uniform bound in the boundary}) and the maximum
 principle for holomorphic function along a long strip.  The main point is that,  for this point $z_m'$,  the
 inequality (\ref{comp:nontrivialarea}) holds.  Consequently, there exists at most a finite number of bubble
 points. Next we want to argue that there is no bubble point at all. For this purpose, we
 consider
 two cases:   $z_\infty \in \p \Sigma $ or
$z_\infty \in \Sigma^0.\; $ In both cases, we want to show that
the existence of a non trivial bubble  must lead to contradiction.


{\bf Part I: No bubble in the boundary}.

 We choose a number $\varsigma \in B_1, $ to be fixed throughout
 the following argument. Set
\[
\phi^{(m)} (\zeta) = z_m +  {\zeta \cdot \epsilon_m}, \qquad
\forall \zeta \in {\bf C}
\] and
\[ l^{(m)} = f^{(m)} \circ \phi^{(m)},\qquad {\rm and}\;\; \tilde{\;z\;}_m =
\phi^{(m)}(\varsigma).\]
 By definition, we have
\[
   |\tilde{\;z\;}_m - \;z_m| \leq \epsilon_m \rightarrow 0, \qquad \forall |\varsigma|\leq 1.
\]
Thus $\displaystyle \lim_{m \rightarrow \infty} \; \tilde{\;z\;}_m
= \;z_\infty.\;$   First set
\[
  \tilde{x}_m =\pi \circ f^{(m)}(\phi^{(m)}(\varsigma)).
\]
 Note that
\[
\begin{array}{lcl}
  l^{(m)} (\varsigma)
  & = & (\tilde{\;z\;}_{m}, \tilde{x}_m, \zeta^{(m)}(\tilde{\;z\;}_{m},
  \tilde{x}_m))\\
   & = & (\phi^{(m)}(\varsigma), \pi\circ f^{(m)}(\phi^{(m)}(\varsigma)), \zeta^{(m)}
  (\phi^{(m)}(\varsigma), \pi\circ f^{(m)}(\phi^{(m)}(\varsigma)))).
  \end{array}
\]
Since $\displaystyle \lim_{m \rightarrow \infty} l^{(m)} =
l^{\infty},$ we may assume
\[
 \displaystyle \lim_{m \rightarrow \infty} \tilde{x}_m
 =\pi\circ l^\infty(\varsigma).
\]
Set
\[
  \displaystyle \lim_{m \rightarrow \infty} \;
  \zeta^{(m)}(\tilde{\;z\;}_{m}, \tilde{x}_m) = \zeta^{\infty}
  (\varsigma).
\]
for some function $\zeta^{\infty}.\;$ We want to show
\[
  \zeta^{\infty} (\varsigma) = \b \p \phi^\infty(\;z_\infty,
  \pi\circ l^{\infty}(\varsigma)).
\]
Now for any $ 0 < \alpha < 1$, there exists a uniform constant $C$
such that
\[
  {{|\zeta^{(m)}(\tilde{\;z\;}_{m}, \tilde{x}_m)
  - \zeta^{(m)}(\;z_\infty, \pi\circ l^{\infty}(\varsigma))| }\over {(|\tilde{z}_m - \;z_\infty|
  + |\tilde{x}_m - \pi\circ l^\infty(\varsigma)|)^{\alpha}
  }} < C.
\]
Since $\displaystyle \lim_{m \rightarrow \infty} \;
(|\tilde{\;z\;}_m - \;z_\infty| + |\tilde{x}_m - \pi\circ
l^\infty(\varsigma)|) = 0, $ we have
\[
 \displaystyle \lim_{m \rightarrow \infty}\;
\left(\zeta^{(m)}(\tilde{z}_{m}, \tilde{x}_m)
  - \zeta^{(m)}(\;z_\infty, \pi\circ l^{\infty}(\varsigma)) \right) = 0.
\]
On the other hand, we have
\[
  \begin{array}{lcl}& & \displaystyle \lim_{m \rightarrow \infty}\;
\left(\zeta^{(m)}(\;z_\infty, \pi\circ l^{\infty}(\varsigma) )
  - \b \p \phi^{\infty}(\;z_\infty, \pi\circ l^{\infty}(\varsigma)) \right)\\
   & & \qquad \qquad = \displaystyle \lim_{m \rightarrow \infty}\;
\left(\b \p \phi^{(m)}(\;z_\infty, \pi\circ l^{\infty}(\varsigma)
)
  - \b \p \phi^{\infty}(\;z_\infty, \pi\circ l^{\infty}(\varsigma))
  \right)\\ &  & \qquad\qquad = \; 0.
  \end{array}
\]
Thus,
\[
\displaystyle \lim_{m \rightarrow \infty}\;
\left(\zeta^{(m)}(\tilde{\;z\;}_{m}, \tilde{x}_m)
  - \b \p \phi^{\infty}(\;z_\infty, \pi\circ l^{\infty}(\varsigma)) \right)
  = 0.
\]
Consequently,
\[
\zeta^{\infty} (\varsigma) = \b \p \phi^\infty(\;z_\infty,
  \pi\circ l^{\infty}(\varsigma)).
\]
Thus
\[
  l^{\infty}(\varsigma) = (\;z_\infty, \pi\circ l^\infty(\varsigma), \b \p \phi^{\infty}(\;z_\infty,
  \pi\circ l^\infty(\varsigma))) \in \{z_\infty\} \times \Lambda_{z_\infty,
  \phi^{\infty}(z_\infty)}.
\]
Since $\varsigma \in B_1$ is chosen randomly, we have
\[
  l^\infty(B_1) \subset  \{z_\infty\} \times
  \Lambda_{z_\infty,\phi^{\infty}(z_\infty)}\subset \{\;z_\infty\}
  \times \cW_M.
\]
Next note that $\;z_\infty \in \p \Sigma, $ we have
\[
 \phi^{\infty}(z_\infty, \cdot) = {\phi_0}^\infty(z_\infty,\cdot).
\]
Thus
\[
\Lambda_{z_\infty,\phi^{\infty}(z_\infty)} =
\Lambda_{z_\infty,{\phi_0}^{\infty}(z_\infty)} \subset
\{z_\infty\}
  \times \cW_M
\]
is a totally real sub-manifold. This contradicts with the fact
that $l^{\infty}(B_1) \subset \{z_\infty\} \times \cW_M $ is a
holomorphic disk! Consequently, there is no bubble sphere/disk
developing at the boundary of $\Sigma.\;$\\

Now we proceed to part 2.\\

\noindent {\bf Part 2: Non existence of interior bubbles}. Suppose
that $\;z_\infty \in \Sigma^0 = \Sigma \setminus \p \Sigma.\;$ Set
\[
  \phi^{(m)}(\varsigma) = z_m + {\epsilon_m\over \delta_m}
  \varsigma
\]
and $l^{(m)} = f^{(m)} \circ \phi^{(m)}$ for $\epsilon_m \ll
\delta_m \rightarrow 0.\;$ For any $k$ fixed, the following map
 has a non-trivial limit:
 \[
   l^\infty = \displaystyle \lim_{m \rightarrow \infty} l^{(m)}: B_k \subset \displaystyle \bigcup_{m=1}^\infty\;
   B_{\delta_m} \rightarrow \Sigma \times \cW_M.
 \]
 Let $k \rightarrow \infty. \;$ Then $l^\infty$
 defines a holomorphic map from $R^2$ to $\cW_M$ with bounded
 area in $\Sigma\times M.\;$ Therefore, the image must be a holomorphic $S^2.\;$

 Here
 \[
   {{\p l^{(m)\alpha}}\over {\p \varsigma}} = {{\p f^{(m)\alpha}}\over {\p
   z}}  {\epsilon_m\over \delta_m} = - g^{(m)\; \alpha
\overline{\beta}}{{\partial^2 \phi^{(m)}}\over{\partial z
\overline{\partial w_\beta}}} {\epsilon_m\over \delta_m}.
 \]
 Set
\[\left\{\begin{array}{lcl}
 \eta^{(m)\alpha} & = & - g^{(m)\; \alpha
\overline{\beta}}{{\partial^2 \phi^{(m)}}\over{\partial z
\overline{\partial w_\beta}}},\nonumber\\
\tilde{\eta}^{(m)\alpha} &  = &- g^{(m)\; \alpha
\overline{\beta}}{{\partial^2 \phi^{(m)}}\over{\partial z
\overline{\partial w_\beta}}} {\epsilon_m\over \delta_m}.
 \end{array}\right.\]
 By assumption, $\tilde{\eta}^{(m)}$ has a non-trivial limit
 $\eta^{\infty}$ such that
 \[
{{\p l^{(m)\alpha}}\over {\p \varsigma}} = \eta^{(m)
\alpha},\qquad {\rm and}\;\;\;
 {{\p l^{ \infty\alpha}}\over {\p \varsigma}}  =  \eta^{\infty
 \alpha}. \]
 The above equations imply
 \[
 \begin{array}{lcl}  \p_\varsigma
({{\p \phi^{(m)}} \over {\p w_{\b \beta}}}  \circ l^{(m)}) & = &
\left({\p \over \p \varsigma} + \tilde{\eta}^{(m) \alpha} {\p
\over {\p w_\alpha}}\right){{\p \phi^{(m)}} \over {\p w_{\b
\beta}}}\\ & = & {{\epsilon_n}\over {\delta_n}} {{\p^2 \phi^{(m)}}
\over {\p z \p w_{\b \beta}}} + \tilde{\eta}^{(m) \alpha} {{\p^2
\phi^{(m)}} \over {\p w_\alpha \p w_{\b \beta}}}\\ & = &
{{\epsilon_n}\over {\delta_n}} {{\p^2 \phi^{(m)}} \over {\p z \p
w_{\b \beta}}} + \tilde{\eta}^{(m)
\alpha} (g^{(m)}_{\alpha \b \beta} - g_{0,\alpha \b \beta})\\
& = & {{\epsilon_n}\over {\delta_n}} {{\p^2 \phi^{(m)}} \over {\p
z \p w_{\b \beta}}} +  (- g^{(m)\; \alpha
\overline{r}}{{\partial^2 \phi^{(m)}}\over{\partial z
\overline{\partial w_r}}} {\epsilon_m\over \delta_m})
g^{(m)}_{\alpha \b \beta} -\tilde{\eta}^{(m) \alpha} g_{0,\alpha
\b \beta}\\ & = & - g_{0,\alpha \b \beta} \tilde{\eta}^{(m)
\alpha}.
\end{array}
 \]
Thus
\[
\begin{array}{lcl} - g_{0,\alpha \b \beta} \tilde{\eta}^{(m)
\alpha} \tilde{\eta}^{(m)\b \beta} & = & \p_\varsigma ({{\p
\phi^{(m)}} \over {\p w_{\b \beta}}}  \circ l^{(m)})
\tilde{\eta}^{(m)\b \beta} \\
& = & \p_\varsigma \left( \tilde{\eta}^{(m)\b \beta} {{\p
\phi^{(m)}} \over {\p w_{\b \beta}}} \circ l^{(m)} \right)\\
& = & \p_\varsigma \left( \p_{\b \varsigma} (\phi^{(m)} \circ
l^{(m)}) - {\epsilon_n \over \delta_n}  {{\p \phi^{(m)}} \over
{\p \b z}} \circ l^{(m)} \right) \\
& = & \p_\varsigma  \p_{\b \varsigma} (\phi^{(m)} \circ l^{(m)}) -
{\epsilon_n \over \delta_n} \p_\varsigma ({{\p \phi^{(m)}} \over
{\p \b z}} \circ l^{(m)}).
\end{array}
\]
Let $\chi(\varsigma) $ be any smooth test function which vanish
outside a compact domain of $R^2.\;$ For $m$ large enough, the
domain of $\chi$ is contained inside in the  domain of
$l^{(m)}.\;$ Then,
\[
\begin{array}{lcl} & & - \displaystyle \int_{R^2} \; \chi\; g_{0,\alpha \b \beta} \tilde{\eta}^{(m)
\alpha} \tilde{\eta}^{(m)\b \beta}\; | d\,\varsigma|^2\\ &  &
\qquad \qquad =  \displaystyle \int_{R^2} \; \chi\; \p_\varsigma
\p_{\b \varsigma} (\phi^{(m)} \circ l^{(m)}) \; | d\,\varsigma|^2
- \displaystyle \int_{R^2} \; \chi\;  {\epsilon_m \over \delta_m}
\p_\varsigma ({{\p \phi^{(m)}} \over {\p \b z}} \circ l^{(m)} )
\; | d\,\varsigma|^2\\
&  & \qquad \qquad = - \displaystyle \int_{R^2} \;{ {\p \chi}\over
{ \p \varsigma}} \; \p_{\b \varsigma} (\phi^{(m)} \circ l^{(m)})
\; | d\,\varsigma|^2  +  {\epsilon_m \over \delta_m} \displaystyle
\int_{R^2} \; {\p \chi\over {\p \varsigma}} {{\p \phi^{(m)}} \over
{\p \b z}} \circ l^{(m)} \; | d\,\varsigma|^2
\end{array}
\]
Taking limit as $m \rightarrow \infty$, we have
\[
- \displaystyle \int_{R^2} \; \chi\; g_{0,\alpha \b \beta}
\tilde{\eta}^{\infty\; \alpha} \tilde{\eta}^{\infty\; \b \beta}\;
| d\,\varsigma|^2 = - \displaystyle \int_{R^2} \;{ {\p \chi}\over
{ \p \varsigma}} \; \p_{\b \varsigma} (\phi^{\infty} \circ
l^{\infty}) \; | d\,\varsigma|^2.
\]

 This holds for any test function in $\RR^2 .\;$ Now the image of
 $l^\infty$ is a smooth $S^2$ in $\{z_\infty\} \times M.\;$
 Therefore, $l^{\infty*} \omega_0$ is an induced (smooth)
 K\"ahler form in this class. We may as well assume $l^{\infty*} \omega_0$ is
 cohomologous to the standard K\"ahler form in $S^2.\;$ Then there exists a bounded
 smooth function $\lambda$ in this $S^2$ such that
 \[
 \begin{array}{lcl}
l^{\infty*} \omega_0 & = & g_{0,\alpha \b \beta}
\tilde{\eta}^{\infty\; \alpha} \tilde{\eta}^{\infty\; \b \beta}\;
| d\,\varsigma|^2   \\ & = & {{\p^2 }\over {\p \varsigma \p \b
\varsigma}}
 \left( - 2 \ln (1 + |\varsigma|^2) + \lambda \circ l^\infty \right)\; d\,\varsigma\; d\,\b \varsigma.
 \end{array}
 \]
 Then
 \[
\displaystyle \int_{R^2} \; {{\p \chi(\varsigma)}\over {\p
\varsigma}} \cdot {{\p } \over {\b \p \varsigma}} \left(- 2 \ln (1
+ |\varsigma|^2) + \lambda \circ l^\infty + \phi^{\infty} \circ
l^{\infty}\right)\; |d\,\varsigma|^2 = 0.
 \]
 Set
 \[
 \Phi = \left( - 2 \ln (1 + |\varsigma|^2) + \lambda \circ l^\infty +
\phi^{\infty} \circ l^{\infty}\right).
 \]
Since $\chi$ is an arbitrary test function, this implies that
$\Phi$ is a weakly $C^{1,1},\;$ one side bounded harmonic function
in ${\bf R}^2.\;$ Then $\Phi$ must be a constant function $c.\;$
Then
 \[
\phi^{\infty} \circ l^{\infty} = c + 2 \ln (1 + |\varsigma|^2) +
\lambda \circ l^\infty
 \]
 is not a bounded function as $\varsigma\rightarrow \infty.\;$
 This contradicts with the fact that $\phi^\infty$ is uniformly
 bounded, which in turn implies that there is no bubble in the
 interior. \\

 The proof of this theorem is then completed.
 \end{proof}

\subsection{Proof of Theorem \ref{comp:superregularlimit}}
\subsubsection{Uniform $C^1$ transversal derivatives of almost super regular foliations}
In this subsection, we continue to derive the {\it a priori}
regularity
estimate for the disk with uniform upper bound on area and capacity.\\

   \begin{theo} \label{comp:volumeinteriorbound} Let $\Omega $ be any compact sub-domain in $\Sigma^0.\;$  For any holomorphic disk
  $f\in {\cal C}(\delta, \Lambda, L_0,L_1, L_2),$   there exists a constant $C> 1$ such that:
\[
{1\over C} \leq  \left({{\omega_0^n} \over {\omega_{\phi^{}}^n}
                     }\right)\circ f \leq C.
                     \]
Here $C$ depends on $\delta, \Lambda, L_0, L_2 $ and $d(\p \Omega,
\p \Sigma).\;$ Moreover, this constant approaches to $\infty$ if $
d(\p \Omega, \p  \Sigma) \rightarrow 0.\;$

   \end{theo}
\begin{proof} Proposition \ref{local:interiorboundofcurvature} implies that the trace of curvature of the $TM$ bundle over the
disk $\Sigma$ has interior estimates:

\[   0 \leq S^{}\circ f^{}(\;z\;, \cdot) \leq  {C_1\over {d(\;z\;,\p \Sigma)^2}}.
\]
Corollary \ref{local:linecurvature} then implies
\[
   0\; \leq\; \p_z \bar \p_z\; (\ln  {{\omega_{\phi^{}}^n} \over {\omega_{0}^n}})  \leq  C_2, \; \;\forall \;z\; \in
   \Omega,\]
for some constants $C_2$ which depends on $C_1.\;$  On the other
hand, finite capacity implies
   \[\
 \int_{\Omega} \;\left({{\omega_0^n} \over
{\omega_{\phi}^n}
                     }\right)\circ f \; {\sqrt{-1}\over 2} \;d\,\;z\;
                     d\,\bar z\;  \leq  L_2.
\]
 This in
turn implies that in a slightly smaller sub-domain $\Omega_1
\subsetneq \Omega$, we have
\[
|\ln {{\omega_{\phi^{}}^n} \over {\omega_{0}^n}} - C_3| \leq C_4
\]
for some constant $C_3$ which might depends on $\phi.\;$ The Key
observation is that $C_4$ depends only on $C_2,  C_3 $ and
$L_2.\;$ Consequently, there exists a  constant $C_5$ such that
\[ {1\over C_5} \leq
\left({{\omega_0^n} \over {\omega_{\phi^{}}^n}
                     }\right)\circ f \leq C_5.
\]
Here $C_5$ depends only on $C_1, L_2.\;$
\end{proof}

Note that $\nabla_\p, \nabla_{\b \p}$ induces naturally a map from
$TM $ to $TM \bigotimes T^*M.\;$ The image of  the leaf vector field
$X$  under these two operators are of particular interest. In a
local coordinate chart, we have
\[
\nabla_{\p} X = \left( {{\p \eta^\alpha}\over {\p w^\beta}}
 + \eta^\nu \Gamma^{\alpha}_{\nu \beta}(g_0)
\right) \cdot {\p \over {\p w^\alpha}} \bigotimes d w^\beta =
{\eta^{\alpha}}_{,\beta} {\p \over {\p w^\alpha}} \bigotimes d
w^\beta
\]
and
\[
\nabla_{\b \p} X = {{\p \eta^\alpha}\over {\p w^{\b \beta}}}  \cdot
{\p \over {\p w^\alpha}} \bigotimes d w^{\b \beta} =
{\eta^{\alpha}}_{, \b \beta}  {\p \over {\p w^\alpha}} \bigotimes d
w^{\b \beta}.
\]
\begin{theo}   \label{comp:1stderivboundeta} For any super regular disk $f\in {\cal C}(\delta,\Lambda, L_0, L_1, L_2), $
 the  $T^{*(1,0)}M$ component of the first transversal derivatives of the
 leaf
vector field is bounded in any compact sub-domain $\Omega\subset
\Sigma.\;$ Namely, there exists a constant $C$ depends on
$\delta,\Lambda, L_0, L_1, L_2$ and $d(\p\Omega, \p \Sigma),$ such
that (c.f. equation (\ref{comp:diskvectorfield}))
    \[
      \| \nabla_{\p} X\|_{g_0} < C, \qquad {\rm in}\;\; \Omega.
    \]
  Moreover, this constant $C$ blows up if $d(\p \Omega, \p \Sigma) \rightarrow 0. $
       \end{theo}
\begin{proof} For any $z\in \Omega$,
consider function $d(z,
\partial \Omega) \cdot \|\nabla_{\p} X\|_{g_0}(z) .\;$ This is a non-negative
 function in $\Omega$ which vanishes on $\partial \Omega.\;$
The maximum value must be attained in $\Omega^0.\;$ If this
theorem is false, then there exists a sequence of super regular
holomorphic disks $\{f^{(m)}\} \subset {\cal C}(\delta,\Lambda,
L_0,L_1,L_2)$ such that
\[
    \displaystyle \lim_{m\rightarrow \infty}\;\displaystyle \max_{\Omega} d(\;z\;,
    \partial \Omega) \cdot
    \|\nabla_{\p} X^{(m)}\|_{g_0}(\;z\;) =  \infty.
\]
Without loss of generality, one may assume that the maximum is
attained at the point $\;z_m.$  Set
\[\displaystyle \lim_{m \rightarrow \infty}
\;z_m = \;z_\infty \in \b \Omega.
\]

On the other hand, Theorem \ref{comp:nobubble} implies that there
exists a subsequence of $f^{(m)}$ which converges in $ {\cal
C}(\delta,\Lambda, L_0,L_1,L_2)$ as an embedded holomorphic disk.
W.l.o.g., we may assume that $f^{(m)}$ is fixed but the restricted
TM bundle varies. Denote the sequence of re-scaling factors as
\[ {1\over \epsilon_m} =  \|\nabla_{\p} X\|_{g_0}
(\;z_m) \rightarrow \infty.
\]

Write this sequence of disks as
  \[
\begin{array}{lllllll}
f^{(m)}:& \Sigma & \mapsto &  \Sigma \times M & \hookrightarrow & \Sigma \times \cW_M\\
  &   \;z\;            & \mapsto  &  (z, f^{(m)}(\;z\;) ) & \hookrightarrow
  & (z, f^{(m)}(\;z\;), \xi^{(m)} )
  \end{array}
\]
where
\[
\xi^{(m)\alpha}(z) = {{\p (\rho + \phi)}\over {\p w_\alpha}} \circ
f^{(m)}(z), \qquad \forall \alpha =1, 2,\cdots n.
\]
Moreover
\[
  X^{(m)} = \displaystyle \sum_{i=1}^n\; \eta^{(m)i} {\p \over {\p
  w_i}}
\]
and
\[
 \eta^{(m)\alpha } = {{\partial f^{(m) \alpha} }\over{\partial z}} =
- g^{(m)\; \alpha \overline{\beta}}{{\partial^2
\phi^{(m)}}\over{\partial z \overline{\partial w_\beta}}}.
\]

 Theorem \ref{comp:volumeinteriorbound}
implies that  there is a positive number $ C_3 > 0$ such that \[
 {1\over C} \leq
\left({{\omega_0^n} \over {\omega_{\phi^{(m)}}^n}
                     }\right)\circ f^{(m)} \leq C_3
\]
hold uniformly in $\Omega $ since $d(\Omega, \p \Sigma) > 0.\;$
Combining this with the $C^{1,1}$ estimate in \cite{chen991},
there exists a small positive constant $\epsilon_0
> 0,$ such that
\[
 \epsilon_0 \left(g_{0, i\b j}\right)_{n \times n} \leq  \left(g^{(m)}_{ i\b j}\right)_{n \times
 n} \leq C\left(g_{0, i\b j}\right)_{n \times n}
\]
holds for these disks on $\Omega.\;$ Set
 \[ \tilde{\phi}^{(m)}( \;z\;, w) =
\phi^{(m)}(\;z_m + \epsilon_m \cdot \;z\;, w), \qquad
\tilde{X}^{(m)} = \displaystyle \sum_{\alpha=1}^n\;
\tilde{\eta}^{(m)\alpha} {\p \over {\p w_\alpha}}
\]
where
\[\begin{array}{lcl}
  \tilde{\eta}^{(m)\alpha} ( \;z\;, w) & =&  g^{(m) \alpha \b \beta} {{\p^2
  \tilde{\phi}^{(m)}}\over {\p  \;z\; \p w_{\b \beta}}}\\
  & = & \epsilon_m\;
  \eta^{(m)\alpha}(\epsilon_m \cdot \;z\; + \;z_m, w).
  \end{array}
\]
Moreover,
\[
\|\nabla_\p \tilde{X}^{(m)}\|_{g_0} (0) = \epsilon_m \cdot
\|\nabla_\p {X}^{(m)}\|_{g_0} (\;z_m) = 1.
\]
Set
\[
  \Omega^{(m)} =  \{\;z\; | \epsilon_m \cdot \;z\; + \;z_m \in
   \Omega\}.
\]

 Clearly
\[ \begin{array}{lcl} \displaystyle \lim_{m\rightarrow \infty}\; d(0, \p
\Omega^{(m)})& = &
   \displaystyle \lim_{m\rightarrow \infty}\; d(0, \p \Omega^{(m)})
   \cdot |\nabla_\p \tilde{X}^{(m)} |_{g_0} (0) \\
   & = & \displaystyle
   \lim_{m\rightarrow \infty}\; d(z_m,\p \Omega) \cdot  \|\nabla_\p X^{(m)}\|_{g_0}
(\;z_m)  = \infty.\;\end{array}
\]
In other words,
 \[ \Omega^{(m)} \rightarrow \RR^2. \]
 Here
\[\begin{array}{lcl}
\tilde{S}^{(m)}(\;z\;, w ) & = & {{\partial
\tilde{\eta}^{(m)\overline{\beta}}}\over {\partial w_{\alpha}}}
{{\partial \tilde{\eta}^{(m)\alpha}}\over {\partial
\overline{w_{\beta}}}} \\
& = & \epsilon_m^{2} {{\partial {\eta}^{(m)\overline{\beta}}}\over
{\partial w_{\alpha}}} {{\partial {\eta}^{(m)\alpha}}\over
{\partial \overline{w_{\beta}}}}\\ & = & \epsilon_m^2 \cdot
S^{(m)}( \;z_m + \epsilon_m \;z\;, w).
\end{array}
\]
Then $\tilde{S}^{(m)}$ still satisfies the inequality
\[
\p_z \p_{\b z}\;\tilde{S}^{(m)} \geq {2\over n} \tilde{S}^{(m)\;
2}
\]
in $\Omega^{(m)}.\;$  Consequently, we have (cf. Proposition.
\ref{local:interiorboundofcurvature})
\[
 0 \leq  \tilde{S}^{(m)}(\;z\;) \leq {{C}\over {d(\;z\;, \partial \Omega^{(m)})^2}} \rightarrow
  0,\qquad \forall\; \;z\; \in \Omega^{(k)}
\]
for any fixed $k$ and $m \rightarrow \infty.\;$ Recall that \[
\begin{array}{lcl} \tilde{S}^{(m)} & = & \epsilon_m^2 {{\p
\eta^{(m)\alpha}} \over {\p w_{\b\beta}}} {{\p \eta^{(m) \b
\beta}} \over {\p w_{\alpha}}}\\
& = & \epsilon_m^2 g^{(m)\alpha \b a} g^{(m) \b \beta b}
\left({{\p \phi^{(m)}} \over {\p z}}\right)_{, \b a \b \beta}
\left({{\p \phi^{(m)}} \over {\p \b z}}\right)_{, \alpha a}\\
& \geq & C^{-2} \epsilon_m^2 g_0^{\alpha \b a} g_0^{ \b \beta b}
\left({{\p \phi^{(m)}} \over {\p z}}\right)_{, \b a \b \beta}
\left({{\p \phi^{(m)}} \over {\p \b z}}\right)_{, \alpha a}.
\end{array}
 \]
 The last inequality holds because $g^{(m)}_{\alpha\b \beta}$ has a uniform upper bound.
 Thus,
 \[
\displaystyle \lim_{m \rightarrow \infty}\; \epsilon_m \cdot
\left({{\p \phi^{(m)}} \over {\p \b z}}\right)_{, \alpha a} = 0.
 \]
Consequently, we have
\[
\begin{array}{lcl} \displaystyle
\lim_{m \rightarrow \infty}\;{{\p \tilde{\eta}^{(m)\alpha}} \over
{\p w_{\b\beta}}} & = & \displaystyle \lim_{m \rightarrow
\infty}\;\epsilon_m \cdot {{\p \eta^{(m)\alpha}} \over {\p
w_{\b\beta}}} \\& = &  \displaystyle \lim_{m \rightarrow \infty}\;
\epsilon_m \cdot g^{(m) \alpha \b a} \left({{\p \phi^{(m)}} \over
{\p z}}\right)_{, \b \beta \b a} =0.
\end{array}
\]
The last inequality holds since $g^{(m)}_{\alpha \b\beta}$ has a
uniform positive lower bound on $\Omega.\;$

Moreover, for any fixed $\;z\;$, we have
\[
\|\nabla_\p \tilde{X}^{(m)}\|_{g_0} (\;z\;) d(\;z\;, \p
\Omega^{(m)}) \leq \|\nabla_\p {X}^{(m)}\|_{g_0} (0) d(0, \p
\Omega^{(m)}).
\]
Therefore \footnote{The second inequality holds since $
\displaystyle \lim_{m \rightarrow \infty} \;
 d(0, \p \Omega^{(m)}) =\infty$ while $ d(0,\;z\;) = |\;z\;|$ is fixed.}
 \[\begin{array}{lcl} \|\nabla_\p \tilde{X}^{(m)}\|_{g_0}
(\;z\;) & \leq  &\|\nabla_\p \tilde{X}^{(m)}\|_{g_0} (0) {d(0,
\p \Omega^{(m)}) \over {d(z\, \p \Omega^{(m)})}} \\
& \leq &  2 \cdot \|\nabla_\p \tilde{X}^{(m)}\|_{g_0} (0)
\\ & \leq & 4.
\end{array}
\]
Consequently, $\nabla_\p \tilde{X}^{(m)}$ \footnote{Note that
\[
       \tilde{\eta}^{(m)\;\alpha}_{,i}  =  {{\p \tilde{\eta}^{(m) \alpha}} \over {\p  w_i}} +
      \tilde{ \eta} ^{(m)\beta} \Gamma^{\alpha}_{\beta i} (g_0).
    \]
    where
\[ \displaystyle \lim_{m \rightarrow \infty} \tilde{ \eta}
^{\beta} \Gamma^{\alpha}_{\beta i} (g_0) = \displaystyle \lim_{m
\rightarrow \infty}  \epsilon_m \cdot \eta^{(m)\;\beta}
\Gamma^{\alpha}_{\beta i} (g_0) = 0. \] This is because the disk,
when restricted in $\Omega$, uniformly converge to a smooth limit
surface. Hence $\Gamma^{\alpha}_{\beta i} (g_0) $ is uniformly
bounded. Thus, \[{{\p \tilde{\eta}^{(m)}} \over {\p w_i}}\approx
\tilde{\eta}^{(m)\;\alpha}_{,i}.
\]} is uniformly bounded and $|\nabla_\p \tilde{X}^{(m)}|_{g_0} \approx 1 $ at the origin. By
Lemma \ref{local:regularity}, both ${{\p \tilde{\eta}^{(m)}} \over
{\p w_i}}$ and $ {{\p \tilde{\eta}^{(m)}} \over {\b \p w_i}} $ are
uniformly $ C^{\alpha} (\forall \alpha < 1)$ bounded. Since
$g^{(m)}_{i \b j}$ has a uniform upper and lower bound in $\Omega,
$ we have \[ \displaystyle \lim_{m \rightarrow \infty}{{\p
\tilde{\eta}^{(m)}} \over {\b \p w_i}} (z, \cdot) = 0.
\]
On the other hand,  ${{\p \tilde{\eta}^{(m)}} \over {\p w_i}}$ is
a bounded holomorphic function in the limit since ${{\p
\tilde{\eta}^{(m)}} \over {\p \b w_i}} =0$ in the limit (cf. Lemma
\ref{local:regularity}). Therefore, ${{\p \tilde{\eta}^{(m)}}
\over {\p w_i}}$ is a constant matrix everywhere in the limit!
 Set
\[
   g^\infty_{i \b j} = \displaystyle \lim_{m \rightarrow \infty}
   g^{(m)}_{i \b j}, \qquad \tilde{\eta}^\infty =  \displaystyle \lim_{m \rightarrow
   \infty}\tilde{\eta}^{(m)}.
\]
Then
\[
   {1\over C}\; I_{n \times n} \leq \left( g^{(m)}_{i \bar j}\right) \leq
   C\; I_{n \times n}.
\]
Theorem \ref{local:curvaturepositive} takes the form

\[
\partial_z   \overline{\partial_z}\; g^\infty_{,i\overline{j}}  = {a^{\alpha}}_{i}  \overline{{a^{\beta}}_j}
g^\infty_{,\overline{\beta}\alpha}.
\]
where
\[
{{\partial \tilde{\eta}^{\infty \; \alpha}}\over {\partial w_i}} =
{a ^{\alpha}}_i, \qquad {\rm and}\qquad  {{\partial
\tilde{\eta}^{\infty\;\overline{\beta}}}\over {\overline{\partial
w_{j}}}} = \overline{{a^{\beta}}_j}
\]
are constant matrixes. This clearly violates the maximum
principle. Note that $g^\infty_{i \bar j}$ and its derivatives are
uniformly bounded in the entire plan.  The contradiction implies
the constant matrix $ \left( {a^{\alpha}}_i\right)$ must vanish
identically which contradicts the blowingup assumption.
\end{proof}
Following  Prop. \ref{local:interiorboundofcurvature}, we can easily
derive the following
\begin{cor}  \label{comp:1stderivboundeta1} For any super regular
 disk $f\in {\cal C}(\delta,\Lambda, L_0, L_1, L_2), $ the $T^*(0,1)M$ component of the
  first transversal derivatives of the leaf
vector field is bounded in any compact sub-domain $\Omega\subset
\Sigma.\;$ Namely, there exists a constant $C$ depending only on
$\delta,\Lambda, L_0, L_1, L_2$ and $d(\p\Omega, \p \Sigma),$ such
that (c.f. equation (\ref{comp:diskvectorfield}))
    \[
      \| \nabla_{\b \p} X\|_{g_0} < C, \qquad {\rm in}\;\; \Omega.
    \]
  Moreover, this constant $C$ blows up if $d(\p \Omega, \p \Sigma) \rightarrow
  0.$
\end{cor}
\subsubsection{The Limit of Super-regular disks with finite capacity is Super regular}

\begin{theo}  \label{comp:limit  super regulardisks}Suppose that
$\{{\phi_0}^{(m)}, \;m\in \NN\}$ is a sequence of loops
in $C^{2,\alpha}(\p\Sigma, {\cal H})$ which satisfies the uniform
bound (\ref{eq:boundpsi}) and $ f^{(m)}: (\Sigma,\p \Sigma)
\rightarrow \left(\Sigma \times \cW_M,
 \bar \Lambda_{{\phi_0}^{(m)}}\right) $ is a sequence of super-regular
 holomorphic disks in ${\cal C}(\delta, \Lambda, L_0,L_1, L_2).\;$
 There exists a subsequence of $\{f_m\}$ which
converges to an embedded, super-regular holomorphic disk.
\end{theo}
Note that this in fact implies Theorem
\ref{comp:superregularlimit}.

\begin{proof} We use $1\leq i,j,k, \alpha,\beta, \gamma \leq n$ to denote the labels
on K\"ahler manifold $M,\;$  and use $1\leq p,q,\cdots \leq 2n $
to denote the labels in the {\it moduli} space $\cU_{\phi_0}.\;$
By definition, we have
\[ {{\p F^{(m)i}} \over { \p \b z\;}} = 0, \qquad \forall i = 1, 2, \cdots.
 \]
 Let $\nu_p^{(m)} = \|{{\p F^{(m)}} \over {\p
 t_p}}\|_{W^{1,2}(\Sigma)}.\;$  Following the standard theory on elliptic problems, there exists a uniform constant $C(\delta, \Lambda)$ (which
 depends on the totally real boundary sub-manifold) such that
 \begin{equation}
     \|{{\p F^{(m)}}\over {\p t_p}}\|_{L^\infty} < C \cdot \nu_p^{(m)}.
 \label{eq:regu00}
 \end{equation}

 On the other hand, we have
 \begin{eqnarray}{{\p}\over {\p \;z\;}} { {\p F^{(m)i}} \over {\p
 t_p}}  
 & = & {{\p \eta^{(m)i}}\over {\p w_\alpha}}  { {\p
F^{(m)\alpha}} \over {\p
 t_p}} + {{\p \eta^{(m)i}}\over { \p w_{\b \beta}}}  { {\p F^{(m)\b \beta}} \over {\p
 t_p}} \label{eq:regu1}\\
{{\p}\over {\p \b z\;}} {{\p F^{(m)i}}\over {\p t_p}} &  = & 0.
\nonumber
\end{eqnarray}
By Theorem \ref{comp:1stderivboundeta}, we have
\begin{equation}
|{{\p \eta^{(m)i}}\over {\p w_\alpha}}|_{g_0} + |{{\p
\eta^{(m)i}}\over { \p w_{\b \beta}}}|_{g_0} \leq C
\label{eq:regu0}
\end{equation}
uniformly in any sub-domain $\Omega \subsetneq \Sigma$ such that
$d(\p \Omega, \p \Sigma) > 0.\;$ Fix $\Omega$ now.\\

For any $\;z_0 \in  \Omega$, and  we re-parametrize the family of
disks  such that for
\[
   F^{(m)}(z_0, p) = (z_0, x(p)),
\]
for any $(z_0,t) $ in the domain of $F^{(m)}.\;$  This is possible
since $f^{(m)}$ is a super regular disk.  By definition, we have
\[
   \left( { {\p F^{(m)i}} \over {\p
 t_p}} \;{ {\p F^{(m)\b i}} \over {\p
 t_p}}    \right)_{2n \times (2n)}
\]
is non-singular at $z=z_0.\;$ From the equation (\ref{eq:regu1})
and {\it a priori} estimate (\ref{eq:regu0}), we have the
following important Harnack type inequality
\[
  C^{-1}\cdot  |{{\p F^{(m)}}\over {\p t_p}} |_{g_0} (\;z_1) \leq |{{\p F^{(m)}}\over {\p t_p}} |_{g_0} (\;z_2) \leq  C \cdot
   |{{\p F^{(m)}}\over {\p t_p}} |_{g_0}
  (\;z_1),\qquad \forall \;z_1,\;\;z_2 \in \Omega.
\]
Here $C$ depends on $\Omega$ only.  Consequently, we have shown
that
\[
  \left( { {\p F^{(m)i}} \over {\p
 t_p}} \;{ {\p F^{(m)\b i}} \over {\p
 t_p}}    \right)_{2n \times (2n)}
\]
  is a bounded and non-singular matrix in $\Omega\;$  In particular,
   there exists a small constant $c$ such that

  \begin{equation}
 \mid \det \left( { {\p F^{(m)i}} \over {\p
 t_p}} \;{ {\p F^{(m)\b i}} \over {\p
 t_p}}    \right)_{2n \times (2n)}\mid  > \;c.
 \label{eq:regu2}
  \end{equation}
Here  $c$  depend only on  (\ref{eq:regu1}). This implies that
$\nu_p ( 1\leq p \leq 2n) $ has a uniform positive low bound. Now
we claim that they all have uniform upper bound:
\begin{equation}
  \displaystyle \sup_{m\rightarrow \infty} \nu^{(m)}_p <
  \infty.
  \end{equation}
  Otherwise, there exists a subsequence (use the same notation for convenience) such that
\[
 \displaystyle \lim_{m\rightarrow \infty} \nu^{(m)}_p = \infty.
\]
   The matrix
   \[
  \left( { {\p F^{(m)i}} \over {\p
 t_p}} \;{ {\p F^{(m)\b i}} \over {\p
 t_p}}    \right)_{2n \times (2n)}
\]
 is uniformly bounded from above and below on $\Omega. \;$
 Set
 \[ u_p^{(m)} = \left(u_p^{(m)1}, u_p^{(m)2},\cdots u_p^{(m)n}\right)\]
 where
 \[ u_p^{(m)\;i} ={ { {\p F^{(m)\;i}} \over {\p
 t_p}} \over {\nu_p^{(m)} }}.\]
  Then,
 $u_p^{(m)} \rightarrow 0$ in $\Omega $ uniformly. Since
 $u_p^{(m)}$ is holomorphic on $\Sigma,$ then
 $u_p^{(m)}$ converges to $0$ at least in $\Sigma^0.\;$ This contradicts  the fact that $u_p^{(m)}$ has a non-zero
 limit.  Consequently, our claim holds
 and $\nu_p$ has a uniform upper bound and
   \[
  \left( { {\p F^{(m)i}} \over {\p
 t_p}} \;{ {\p F^{(m)\b i}} \over {\p
 t_p}}    \right)_{2n \times (2n)}
\]
 is uniformly bound from above on $\Sigma.\; $  Recall that $F^* \omega_{\phi^{(m)}}$ is an invariant form along disc direction.  
 In particular, the pulled back of volume form is constant along the disc.  Thus, we have
  \[
  \begin{array}{lcl} \det\; (g^{(m)}_{pq}) (z_0) &  = &\det \; ( g^{(m)}_{i \b j})(\;z\;) \det \left({ {\p F^{(m)i}} \over {\p
 t_p}}  { {\p F^{(m) \b j}} \over {\p
 t_p}}\right)(\;z\;) \\ & \leq & C \;\det  ( g^{(m)}_{i \b j})(\;z\;), \qquad\forall\; \;z\;\in \Sigma.
 \end{array}
  \]
  Thus $g^{(m)}_{i \b j}$ has a uniform positive lower bound on $\Sigma.\;$
  According to Donaldson \cite{Dona01}, the limiting disk
  must be super-regular.
\end{proof}

\subsection{Proof of Theorem \ref{comp:foliationlimit}:  Compactness of almost super regular foliations}
We give a proof directly based on our work in the previous two
subsections.
 \begin{proof} As before, identify $\cU_{{\phi_0}^{(m)}}$ with an open dense set $\tilde{M}$ of $M.\;$
 For every point $x$ in $ \;\tilde{M},\;$ consider $g_x$ as the holomorphic disk in
$\cU_{{\phi_0}^{(m)}}$ which passes through the point $(z_0,x).\;$
Set
\[
  h^{(m)}(x) = {\rm Cap}(g^{(m)}_x(\Sigma)).
\]  Theorem \ref{comp:boundpsi} then implies that
\[
\int_{\tilde{M}}\; h^{(m)}(x) d\,x =  \int_{\tilde{M}} \;{\rm
Cap}(g^{(m)}_x(\Sigma)) \; d\,x \leq C.
\]
Therefore, for generic points $x \in \tilde{M}, $ there exists a
subsequence of $g^{(m)}_x$ such that
\begin{equation}
  h^{(m)}(x) =  {\rm Cap}(g^{(m)}_{x}) \leq C(x),\qquad \forall \;
  m=1,2,\cdots
    \infty.\label{eq:cap1}
    \end{equation}
  According to Theorem \ref{comp:superregularlimit}, after passing to a subsequence if
  necessary, this sequence $\{g^{(m)}_{x},\;m \in \NN\}$ of disks has a uniform
   limit $g^\infty_{x}.\;$ Moreover,  the limit disc  is super-regular, in particular, regular.
   Therefore,  there exists a small open subset $B_{r(x)}(x) $ in the moduli space such that the foliation
$F^{(m)}\mid_{\Sigma \times B_{r(x)}(x)}: \Sigma \times
B_{r(x)}(x) \rightarrow \Sigma \times \cW_M$ has a unique smooth
limit $F^{(\infty)}: \Sigma \times B_{r(x)}(x)
\rightarrow \Sigma \times \cW_M.\;$ Moreover, $\pi\circ F^{(\infty)}$ induces a foliation in a small
open, tublar neighborhood of $g_x^{(\infty)} \subset \Sigma\times M.\;$
  Consequently, $h^{(m)}$ is
uniformly continuous in $B_{r(x)}(x) \subset \tilde{M}.\;$
Following Lemma \ref{comp:technical3} below, there exists a set $E$ of
measure $0$, and a subsequence of $h^{(m)}$ (ultimately of
$g^{(m)}_x$) such that $h^{(m)}$ is uniformly continuous function
in any compact subset of $\tilde{M}\setminus E.\; $ Theorem
\ref{comp:superregularlimit} again implies that the foliation has
a uniform limit in this compact subset. We define
  $\cU_{{\phi_0}^{(\infty)}} = \displaystyle \lim_{m \rightarrow
\infty}\;F^{(m)}:\Sigma\times \tilde{M}\setminus E
\rightarrow \Sigma
\times \cW_M.\;$ \\

For every $x\in E$, any disk in $\{g^{(m)}_x, \; m\in \NN \}$ has a uniform upper bound on its area. Following
Theorem \ref{comp:nobubble},  there exists a subsequence of  disks
$\{g^{(m)}_x, \; m=1,2, \cdots \infty \}$ which converges to an
embedded holomorphic disk $S_x^{(\infty)}.\;$  This limit depends on the selection of subsequence
apriori and  may not be
unique in general.  However, the image of each limiting disk $S_x^{(\infty)}$
does not intersect
the image of any disk of $\cU_{{\phi_0}^\infty}$  on $\Sigma^0\times M.\;$ \\

Let us expand the open set $F^{(\infty)} \subset \cM_{\phi_0^{(\infty)}}$  by including $E_\infty
= \displaystyle \bigcup_{x \in E} S_x^{(\infty)}.\;$
Then,  $F^{(\infty)}$ is a partially smooth foliation (c.f.
defi. \ref{defo:inducedfoliation}), where
$\cU_{{\phi_0}^{(\infty)}}$ is the set of super regular disks
in $F^{(\infty)}.\;$

The only thing remain is to show that
$\omega_{{\phi_0}^{(\infty)}}^n$ is a continuous form on
$\Sigma^0\times M.\; $ As before, denote the image of super regular discs
under the evaluation map $\pi\circ ev$ as
$\cV_{{\phi_0}}.\;$ Let $\cS_{\phi_0} = \Sigma\times M\setminus
\cV_{{\phi_0}}$ be the union of all of the singular points.
Clearly, $\omega_{{\phi_0}^{(\infty)}}^n$ is a smooth $(n,n)$ form
in $\cV_{{\phi_0}}$ and vanishes completely in its complement $\cS_{\phi_0}..\;$  To show
that this $(n,n)$ form is continuous, we just need to verify that,
for any sequence $(z_i, x_i)\in \cV_{{\phi_0}}$ which converges
to $(\b z,\b x) \in \cS_{\phi_0}, $ we have
\[
\displaystyle \lim_{i\rightarrow \infty}\;
\omega_{{\phi_0}^{(\infty)}}^n (z_i,x_i) = 0.
\]
Set $f_i$ to be the unique super regular disk in
$\cU_{{\phi_0}^{(\infty)}} $ passing through the point
$(z_i,x_i).\;$ Without loss of generality, we may assume that this
sequence of disks converges to
some other disk $f.\;$\\

By definition, $\displaystyle \lim_{i\rightarrow \infty} Cap(f_i)
= \infty.\;$ In other words,
\begin{equation}
 \lim_{i \rightarrow \infty} \int_\Sigma \;
 \left({{\omega_{{\phi_0}^{(\infty)}}^n}\over {\omega^n}}\right)^{-1} (z_i,x_i)\;
 |d\,z|^2 = \infty.
\end{equation}
However, $\log \left({{\omega_{{\phi_0}^{(\infty)}}^n}\over
{\omega^n}}\right)$ is a sub-harmonic function with a uniform
upper bound.  Theorem \ref{local:deltazpotential} and
\ref{local:interiorboundofcurvature} imply that, for any compact
sub-domain $\Omega\subset \Sigma^0,\;$ we have
\[
|\triangle_z\;\log \left({{\omega_{{\phi_0}^{(\infty)}}^n}\over
{\omega^n}}\right)| \leq C_\Omega
\]
where $C_\Omega$ depends on $d(\p\Omega, \p\Sigma).\;$ Choose
$\Omega$ so that $z_i, z\in \Omega.\;$ Harnack inequality for
negative harmonic functions implies that, either $\log
\left({{\omega_{{\phi_0}^{(\infty)}}^n}\over {\omega^n}}\right)$
tends to $-\infty$ simultaneously in $\Omega$ or $\log
\left({{\omega_{{\phi_0}^{(\infty)}}^n}\over {\omega^n}}\right)$
are uniformly bounded from above and below. This dichotomy holds
for any compact sub-domain of $\Sigma^0.\;$ For any $z\in \Sigma^0$ fixed, there is
 a unique map $w_i: \Sigma\rightarrow M$
such that $(z,w_i)$ lies in the image of $\pi\circ f_i.\;$ If (passing to a
subsequence if necessary)  \[ \displaystyle \lim_{i
\rightarrow \infty} \left({{\omega_{{\phi_0}^{(\infty)}}^n}\over
{\omega^n}}\right)(z_i,x_i) = c > 0,
\]
then
\[
  \displaystyle \lim_{i\rightarrow \infty} \;\left({{\omega_{{\phi_0}^{(\infty)}}^n}\over
{\omega^n}}\right)(z, w_i(z)) > 0.
\]
Following the proof of Theorem \ref{comp:volumeinteriorbound} and
Proposition \ref{local:interiorboundofcurvature}, we can show that
in fact,
\[
\displaystyle \lim_{i\rightarrow \infty}
\;\left({{\omega_{{\phi_0}^{(\infty)}}^n}\over
{\omega^n}}\right)(z, w_i(z)) > 0, \qquad \forall z\in \Sigma.
\]
 The fact contradicts the assumption that the capacity of this sequence of disks blows up. Consequently,
   the varying volume form ratio $\left({{\omega_{{\phi_0}^{(\infty)}}^n}\over {\omega^n}}\right) $
   must converges to $0.\;$ In other words, the volume form must be continuous in the
interior of $\Sigma \times M.\;$
\end{proof}

\begin{lem} \label{comp:technical3} Suppose $\{h^{(m)}, m \in \NN\}$ is a sequence of continuous,
positive functions in a fixed domain $\Omega$ which satisfies
the following two conditions
\begin{enumerate}
\item The $L^1$ norm of $h^{(m)}$ is uniformly bounded; \item For
any $x\in \Omega$, if $\displaystyle \sup_{1\leq m \leq \infty}\;
h^{(m)}(x) < \infty,\;$   there exists a small neighborhood ${\cal
O}_x$ of $x$ such that this sequence (pass to a a subsequence if
necessary) of functions $\{h^{(m)}\}_{m=1}^\infty$ is uniformly
continuous on ${\cal O}_x.\;$
\end{enumerate}
Then,  there exists a set $E$ of measure at most $0$ and a
subsequence of $\{h^{(m)}, m \in \NN\}$ such that this
subsequence is uniformly continuous on any compact subset of
$\Omega\setminus E.\;$ Moreover,  there exists a limit function
$h^\infty$ such that $\displaystyle \lim_{m \rightarrow \infty}
h^{(m)}=  h^\infty$ on $M\setminus E.\;$ Moreover, ${1\over
{1+h^\infty}}$ is a continuous function.
\end{lem}
\begin{proof} By an elementary and straightforward argument. \end{proof}

\section{The (modified) K energy along almost solutions }


In this section, we want to prove first that  the K energy
function is sub-harmonic along any almost smooth solution of HCMA
equation \ref{eq:hcma0}.  One can view this as a generalization of
the fact that the K energy functional is
convex along a smooth geodesic.  Secondly, we want to use this
property of subharmonicity  to prove that the (modified) K energy
has a lower bound in any K\"ahler class, provided that there exists  a constant
scalar curvature metric (or an extremal K\"ahler metric) in this class.
\subsection{The Subharmonicity of the K energy}

Suppose that $\cF_{\phi_0}$ is an almost super regular foliation
and $\phi_0: \Sigma \rightarrow \overline{\cal H}$ is an
 almost smooth solution corresponding to it.  Note that  the
evaluation map $ev: \Sigma \times {\cal M}_{\phi_0} \rightarrow
\Sigma \times \cW_M$ is smooth everywhere. The set of holomorphic
disks which are not super regular has codimension at least $1$ in
the {\it moduli} space. Recalled that $\cV_{\phi_0} = \pi\circ ev
(\Sigma\times \cU_{\phi_0}).\;$  As before, set ${\cal S}_{\phi_0}
= \Sigma\times M \setminus \cV_{\phi_0}.\;$  Clearly, ${\cal
S}_{\phi_0}$ is a smooth sub-manifold and ${\cal
S}_{\phi_0}\bigcap (\p\Sigma\times M)$ has codimension at least
$1$ at $\p \Sigma\times M.\;$
We follow  notations in Sections \ref{local:defomation} in
general. For convenience
of the readers, let us re-state Theorem \ref{th:hessianofK-energy} here\\

\begin{theo}\label{appl:kenergyissubharmonic} Suppose that $\phi: \Sigma \rightarrow \overline{\cal H}$
is an almost smooth solution described as in Definition
\ref{def:almostsmooth}. Then the induced K energy function $\bE:
\Sigma\rightarrow \RR$ (by $\bE(z) = \bE(\phi(z,\cdot))$) is
weakly sub-harmonic and $C^1$ continuous (up to the boundary).
More precisely,
\[
{{\p^2}\over {\p z\p\b z}}  \bE(\phi(z,\cdot)) = \int_{\pi\circ
ev(z,\cU_{\phi_0})} | {\cal D} {{\p \phi}\over {\p \b
z}}|_{{\omega_{\phi}}}^2\, {\omega_{\phi}}^n\;\geq 0, \qquad
\forall\; z \in \Sigma^0
\] holds in $\Sigma^0$ in the weak sense. On $\p \Sigma$, we have
\[
 \displaystyle \int_{\p
\Sigma}{{\p \bE}\over{\p\,{\bf n}}} (\phi) d\,s
 =  \displaystyle \int_{\pi\circ
ev(z,\cU_{\phi_0})} | {\cal D} {{\p \phi}\over {\p \b
z}}|_{{\omega_{\phi}}}^2\, {\omega_{\phi}}^n\, d\,s,
\] where  $d\,s$
is the length element of $\p \Sigma,\;$ and $\bf n$ is the outward
pointing unit normal direction at $\p \Sigma.\;$

\end{theo}
Before we proceed to the proof, let us note the following
calculation scheme. For any function $f \in C^\infty(\Sigma\times
M), $ we have
\[
\begin{array}{lcl} {\p\over {\p z}} \;\int_M\; f\;
\omega_{\phi}^n & = & \int_M\; {{\p f}\over{\p z}} \omega_\phi^n +
\int_M\; f\; \triangle_\phi {{\p
\phi}\over{\p z}}\; \omega_\phi^n\\
& = &  \int_M\; {{\p f}\over{\p z}} \omega_\phi^n - \int_M\; {{\p
f}\over {\p w_\alpha}} \; g^{\alpha\b\beta}_\phi\;  {{\p^2
\phi}\over{\p z \p w_{\b \beta}}}\; \omega_\phi^n \\
& = & \int_M\; {{\p f}\over{\p z}} \omega_\phi^n + \int_M\;
{{\p f}\over {\p w_\alpha}} \; \eta^{\alpha}\omega_\phi^n\\
& = & \int_M\;\p_z(f)\;\omega_\phi^n.
\end{array}
\]

Similarly, we have
\[
 {\p\over {\p  \b z}} \;\int_M\; f\;
\omega_{\phi}^n = \int_M\;\p_{\b z}(f)\;\omega_\phi^n.
\]
We will use these schemes throughout the proof below.  We also
need to use the decomposition formula of the K energy given in \cite{chen993}
(c.f.\cite{tian98}). For any $\phi \in \cal H$, we have
\begin{eqnarray} \bE(\phi) & = &  \displaystyle \int_M\;
\ln {{\omega_\phi^n}\over {\omega_0^n}}\; \omega^{[n]} + J(\phi)
+ \underline{R}\; I(\phi),  \label{appl:kenergy}\\
 I(\phi) & = &
\displaystyle \sum_{p=0}^n  {1 \over { p+1}} \displaystyle  \int_M
\; \phi\; \omega_0^{[n-p]} \wedge (\sqrt{-1} \partial
\overline{\partial} \phi )^{[p]}, \label{appl:Ifunc}\\
  J(\phi) & = & - \displaystyle \sum_{p=0}^{n-1}  {1 \over { p+1}} \displaystyle
  \int_M \; \phi\; Ric(\omega_0) \wedge \omega_0^{[n-p-1]} \;
\wedge (\sqrt{-1} \partial \overline{\partial} \phi )^{[p]}.
\label{appl:Jfunc}
\end{eqnarray}
Here $\underline{R}$ is the average of the scalar curvature
function of any K\"ahler metric in $[\omega].\;$ Note that  the
component $I$ makes no contribution to the calculation of the 2nd
derivatives of  the K energy $\bE.\;$ Thus, we can basically leave
it aside as we calculate the second derivatives of the $K\;$ energy.

\begin{proof}


Let $\chi$ be any non-negative function whose support lies inside
of the set $\Sigma\times M\setminus {\cal S}_{\phi_0}.\;$ Set
\[
   K_\chi = \int_M\;\chi \log {\omega_\phi^n\over \omega_0^n}
\; \omega_\phi^{[n]}.\] Then
\[
\begin{array}{lcl} {{\p^2 K_\chi}\over {\p z\p\b z}} & = & \int_M\;\p_{\b z}\;\p_z \left(\chi \log {\omega_\phi^n\over
\omega_0^n}\right) \; \omega_\phi^{[n]}\\
&= & \int_M\;\chi \left(\triangle_z\;\log {\omega_\phi^n\over
\omega_0^n} + \p_{\b z}\;\left( (\p_z \chi) \;
\log{\omega_\phi^n\over \omega^n}\right)\right)\;\omega_\phi^{[n]}
+\int_M \left(\p_z \chi\right) \left( \log {\omega_\phi^n\over
\omega_0^n} \right) \;\omega_\phi^{[n]}\\ & = & \int_M\;\chi \;
\mid {\cal D} {{\p \phi}\over {\p z}}\mid^2_\phi\;
\omega_\phi^{[n]} + \int_M\;\chi \; Ric(\omega_0)_{\alpha\b\beta}
\;\eta^\alpha\eta^{\b\beta}\;\omega_\phi^{[n]}\\ && \qquad  +
{\p\over {\p\b z}} \int_M\; (\p_z \chi) \log {\omega_\phi^n\over
\omega_0^n} \; \omega_\phi^{[n]}  + \int_M\; (\p_z \chi)(\p_{\b
z}\;\log {\omega_\phi^n\over \omega_0^n})\;\omega_\phi^{[n]}.
\end{array}
\]
 Let $v(z)$ be any non-negative cut off function in $\Sigma^0.\;$
\[
\begin{array}{lcl} \int_\Sigma\; K_\chi\; \triangle_z\;v & = &
\int_\Sigma\; v(z)\; \triangle_z K_\chi\\
& = & \int_\Sigma\; v(z)\;\int_M\;\chi \; \mid {\cal D} {{\p
\phi}\over {\p z}}\mid^2_\phi\; \omega_\phi^{[n]} + \int_\Sigma\;
v(z)\; \int_M\;\chi \; Ric(\omega_0)_{\alpha\b\beta}
\;\eta^\alpha\eta^{\b\beta}\;\omega_\phi^{[n]}\\
& & \qquad- \int_\Sigma\; {{\p v}\over {\p\b z}}\cdot  \int_M\;
(\p_z \chi) \log {\omega_\phi^n\over \omega_0^n} \;
\omega_\phi^{[n]}   +  \int_\Sigma\; v\;\int_M\; (\p_z \chi)
\p_{\b z}\left({\omega_\phi^n\over \omega_0^n}\right)\;
\omega_0^{[n]}.
\end{array}
\]

 Consider the
evaluation map:
\[
\begin{array}{llll} \sharp: & \Sigma\times {\cal M_{\phi_0}} &  \rightarrow & \Sigma \times
   M\\ & (z, f) & \rightarrow  & (z, \pi(f(z))).
   \end{array}
\]
 Then, $\sharp$ is
invertible on $\Sigma\times \cU_{\phi_0}.\;$  Consider any
$C^\infty$ function $\phi \in C_0^\infty(\cU_{\phi_0}) \subset
C^\infty({\cal M}_{\phi_0}) \;$ which vanishes on the boundary of
${\cal U}_{\phi_0}.\;$  Set
\begin{equation}
\chi(z, w) = \phi (\sharp^{-1}(z,w)),\qquad \forall\; (z,w)\in
\Sigma\times M. \label{eq:Kenergyconvex2}
\end{equation}
Then $\chi(z,w)$ is a smooth function in $\Sigma\times M$ whose
support  lies completely inside $\Sigma\times M\setminus {\cal
S}_{\phi_0}.\;$    By definition, the disk derivative $\p_z\chi$
vanishes completely along super regular disks. Consequently, for
any cut off function defined via formula
(\ref{eq:Kenergyconvex2}), we have
\[
\int_\Sigma\; K_\chi\; \triangle_z\;v = \int_\Sigma\;
v(z)\;\int_M\;\chi \; \mid {\cal D} {{\p \phi}\over {\p
z}}\mid^2_\phi\; \omega_\phi^{[n]} + \int_\Sigma\; v(z)\;
\int_M\;\chi \; Ric(\omega)_{\alpha\b\beta}
\;\eta^\alpha\eta^{\b\beta}\;\omega_\phi^{[n]}.
\]
Now let $\phi$ tend to Characteristic function of $\cU_{\phi_0}$
inside ${\cal M}_{\phi_0}.\;$ Then, we have:

\[
\begin{array}{lcl} \int_\Sigma \triangle_z v(z)\;\int_M \;\log {\omega_\phi^n\over\omega_0^n}
\;\omega_\phi^{[n]} & = & \int_{\Sigma\times M \setminus {\cal
S}_{\phi_0}}\; \triangle_z v(z)\; \log {\omega_\phi^n\over\omega_0^n}
\;\omega_\phi^{[n]}\\
& = & \int_{\Sigma\times M \setminus {\cal S}_{\phi_0}}\; v(z)\; \; \mid
{\cal D} {{\p \phi}\over {\p z}}\mid^2_\phi\; \omega_\phi^{[n]} +
\int_{\Sigma\times M \setminus {\cal S}_{\phi_0}}\; v(z)\; \;
Ric(\omega)_{\alpha\b\beta}
\;\eta^\alpha\eta^{\b\beta}\;\omega_\phi^{[n]}\\
& = & \int_{\Sigma\times M \setminus {\cal S}_{\phi_0}}\; v(z)\; \; \mid
{\cal D} {{\p \phi}\over {\p z}}\mid^2_\phi\; \omega_\phi^{[n]} +
\int_{\Sigma\times M }\;  v(z)\;\; Ric(\omega)_{\alpha\b\beta}
\;\eta^\alpha\eta^{\b\beta}\;\omega_\phi^{[n]}.
\end{array}
\]
The first and the last equality holds because that \[
 \;\log {\omega_\phi^n\over\omega_0^n}
\;{\omega_\phi^{n}\over \omega_0^n} \qquad {\rm and}\qquad
Ric(\omega_0)_{\alpha\b\beta}
\;\eta^\alpha\eta^{\b\beta}\;{\omega_\phi^{n}\over \omega_0^n}
 \] both vanish on $\cS_{\phi_0}.\;$ On the
other hand,
\[
\begin{array}{lcl}
{{\p^2 J(\phi)}\over {\p z \p \b z}}  &  = & -\int_M\; {{\p^2
\phi}\over {\p z\p \b z}}\; Ric(\omega)\wedge \omega_\phi^{[n-1]}
- \int_M\; {{\p \phi}\over {\p z}}\; Ric(\omega_0)\wedge \sqrt{-1}
\p\b\p {{\p \phi}\over {\p \b
z}}\wedge \omega_\phi^{[n-2]}\\
& = &  -\int_M\; {{\p^2 \phi}\over {\p z\p \b z}}\;
Ric(\omega_0)\wedge \omega_\phi^{[n-1]} + \int_M\;
Ric(\omega_0)\wedge \sqrt{-1} \p {{\p \phi}\over {\p \b z}}\;
\wedge \b\p {{\p \phi}\over {\p z}}\wedge
\omega_\phi^{[n-2]}\\
& = &  -\int_M\; {{\p^2 \phi}\over {\p z\p \b z}}\;
Ric(\omega_0)\wedge \omega_\phi^{[n-1]}\\
&&\qquad\qquad + \int_M\;\left(g^{\alpha\b\beta}_\phi
Ric(\omega_0)_{\alpha\b\beta} \cdot g^{r\b\delta}_\phi {{\p^2
\phi}\over {\p \b z \p w_r}}\; {{\p^2 \phi}\over {\p z \p
w_{\b\delta}}} -  Ric(\omega_0)_{\alpha\b\beta} \eta^\alpha
\eta^{\b\beta} \right) \omega^{[n]}\\
& = &  -\int_M\; \left({{\p^2 \phi}\over {\p z\p \b z}} -
g^{r\b\delta}_\phi \eta^\alpha\;\eta^{\b\beta}\right) \;
Ric(\omega_0)\wedge \omega_\phi^{[n-1]} - \int_M\;
Ric(\omega_0)_{\alpha\b\beta} \eta^\alpha \eta^{\b\beta}
\omega^{[n]}\\
& = &  - \int_M\; Ric(\omega_0)_{\alpha\b\beta} \eta^\alpha
\eta^{\b\beta} \omega^{[n]}.
\end{array}
\]
The last equality holds since $\phi$ is a solution to the
Homogenous complex Monge-Ampere equation. Therefore,
\[
\begin{array}{lcl} \int_\Sigma \triangle_z v\; J(\phi) & = & -
\int_\Sigma\triangle_z v\; \int_M\; Ric(\omega_0)_{\alpha\b\beta}
\eta^\alpha \eta^{\b\beta} \omega^{[n]}
\end{array}
\]
Using the decomposition formula for K energy (\ref{appl:kenergy}),
we have
\begin{equation}
\int_\Sigma (\triangle_z v(z))\; \bE(\phi(z,\cdot)) =
\int_{\Sigma\times M \setminus {\cal S}_{\phi_0}}\; v(z)\; \; \mid {\cal D}
{{\p \phi}\over {\p z}}\mid^2_\phi\; \omega_\phi^{[n]}\geq 0.
\label{eq:Kenergyconvex0}
\end{equation}

This implies that the K energy functional is sub-harmonic in
$\Sigma^0.\;$  Next we want to derive a formula for the first
derivative of the K energy. For any $v\in C_0^\infty (\Sigma),\;$
we have
\[
\begin{array}{lcl}\int_\Sigma  {{\p v(z)}\over {\p z}}\cdot  K_\chi & = &
 \int_\Sigma\; v(z) \cdot {{\p }\over {\p z}} K_\chi \\
& = & \int_\Sigma\; v(z) \left(\int_M\;\;\p_z (\chi)  \log
{\omega_\phi^n\over \omega_0^n} \; \omega_\phi^{[n]} +
\int_M\;\;\chi \p_z  \left({\omega_\phi^n\over
\omega_0^n}\right)  \omega^{[n]}\right).\\
\end{array}
\]
For any small $\delta>0, $ let $\Sigma_\delta = \{z\in \Sigma: |z|
\leq 1 - \delta\}.\;$ Since  $v(z)$ is  an arbitrary compactly
supported function in $\Sigma^0, $ we obtain
\[
\int_{\p \Sigma_\delta} \;\zeta(z)\;  {\p\over {\p z}} K_\chi =
\int_{\p \Sigma_\delta}\; \zeta(z)\;\left(\int_M\;\;\p_z (\chi)
\log {\omega_\phi^n\over \omega^n} \; \omega_\phi^{[n]} +
\int_M\;\; \chi \p_z \left({\omega_\phi^n\over \omega_0^n}\right)
\; \omega_0^{[n]}\right),\] where $\zeta$ is any
smooth function.  \\

Now let $\chi$ tend to the Characteristic function of
$\Sigma\times M\setminus {\cal S}_{\phi_0}.\;$ As before, the first term in
the right hand side vanishes, we have
\[
\int_{\p \Sigma_\delta} \;\zeta(z)  {\p\over {\p z}} \int_M\; \log
{\omega_\phi^n\over \omega_0^n}\; \omega_\phi^{[n]} = \int_{\p
\Sigma_\delta}\; \zeta(z)\;\int_{M\setminus {\cal S}_{\phi_0}}\;\p_z
\left({\omega_\phi^n\over \omega_0^n}\right) \; \omega_0^{[n]}.
\]
For the K energy, we have
\[\begin{array}{lcl}
\int_{\p \Sigma_\delta} \;\zeta(z)\;  {\p\over {\p z}}
\bE(\phi(z,\cdot)) & = & \int_{\p \Sigma_\delta}\;\zeta(z)\;
\int_{M\setminus {\cal S}_{\phi_0}}\;\p_z
\left({\omega_\phi^n\over \omega_0^n}\right) \; \omega_0^n\\
& &\qquad\qquad -\int_{\p \Sigma_\delta}\;\zeta(z)\; \int_M\; {{\p
\phi}\over {\p z}}\; (Ric(\omega_0) - \underline{R}
\omega_\phi)\wedge \omega_\phi^{[n-1]}
\end{array}
\]
Integrating by parts on the left hand side of equation
(\ref{eq:Kenergyconvex0}) and letting $v$ approach  the
characterstic function of $\Sigma_\delta, $  we obtain
\[
\begin{array}{lcl} \int_{\Sigma_\delta \times M \setminus {\cal S}_{\phi_0}}\; \; \mid {\cal D}
{{\p \phi}\over {\p z}}\mid^2_\phi\; \omega_\phi^{[n]} & = &
\int_{\p \Sigma_\delta} \;{{\b z}\over {|z|}}\;  {\p\over {\p z}}
\bE(\phi(z,\cdot))\\
& = & \int_{\p \Sigma_\delta}\;\;{{\b z}\over {|z|}}\;
\int_{M\setminus {\cal S}}\;\p_z \left({\omega_\phi^n\over
\omega_0^n}\right) \omega_0^n -\int_{\Sigma_\delta}\;\;{{\b
z}\over {|z|}}\; \int_M\; {{\p \phi}\over {\p z}}\; (Ric(\omega_0)
- \underline{R} \omega_\phi)\wedge \omega_\phi^{[n-1]}.

\end{array}
\]
The first and the third terms in the formula above,  are both
integration on $\pi\circ ev (\Sigma\times \cU_{\phi_0})$ where
K\"ahler metric is smooth. Therefore, taking limit as $\delta
\rightarrow 0, $ we arrive
\[
\int_{\Sigma\times M \setminus {\cal S}_{\phi_0}}\; \; \mid {\cal
D} {{\p \phi}\over {\p z}}\mid^2_\phi\; \omega_\phi^{[n]} =
\int_{\p \Sigma}\; \;{{\b z}\over {|z|}}\; \int_{M\setminus {\cal
S}_{\phi_0}}\;\p_z \left({\omega_\phi^n\over \omega_0^n}\right)
\omega_0^n -\int_{\p \Sigma}\;\;{{\b z}\over {|z|}}\; \int_M\;
{{\p \phi}\over {\p z}}\; (Ric(\omega_0) - \underline{R}
\omega_\phi)\wedge \omega_\phi^{[n-1]}.
\]
In the above process of taking the limit, the only term which
needs special attentions is:
\[
 \displaystyle \lim_{z\rightarrow z_0\in \p \Sigma}\; \int_{\{z\}\times M\setminus {\cal
S}_{\phi_0}}\;\p_z \left({\omega_\phi^n\over \omega_0^n}\right)
\omega_0^n = \int_{\{z_0\}\times M\setminus {\cal
S}_{\phi_0}}\;\p_z \left({\omega_\phi^n\over \omega_0^n}\right)
\omega_0^n.
\]
This is equivalent to say that the $z-$ derivative of the K energy
is continuous as $z\rightarrow z_0 \in \p \Sigma \;(|z| < 1 =
|z_0|, \;\forall \;z_0 \in \p \Sigma).\; $ For any $\delta>0$
fixed, choose any $\delta$ neighborhood of the set of non-super
regular disks (Denoted by $E_\delta \subset {\cal M}_{\phi_0} $)
such that
\[
\lim_{\delta\rightarrow 0} \; mes (E_\delta) = 0.\;
 \]
Set $${\cal S}_\delta = \sharp(\Sigma\times E_\delta).\;$$ Let
$(t_1,t_2,\cdots t_{2n})$ be the coordinate variables in
$\cM_{\phi_0}$ and $w_1,w_2,\cdots w_n$  be the complex coordinate
variables in $M.\;$ Set
\[
  J = \left({{\p w_\alpha}\over {\p t_i}} {{\p w_{\b \beta}}\over {\p t_i}} \right)
\]
as the Jacobi matrix.    Then, $J$ is a smooth complex  matrix
valued function in $\cal M,$ and invertible at ${\cal M}_{\phi_0}.\;$ Denote
by $\Gamma(\omega_0)$  the connection form of the K\"ahler metric
which corresponds to the K\"ahler form $\omega_0.\;$ In the
following calculation, we take covariant derivatives with respect
to $\omega_0.\;$ Clearly,
\[
 \displaystyle \lim_{z\rightarrow z_0\in \p \Sigma}\;
\int_{\{z\}\times M\setminus {\cal S}_\delta}\;\p_z
\left({\omega_\phi^n\over \omega_0^n}\right) \omega_0^n =
\int_{\{z_0\}\times M\setminus {\cal S}_\delta}\;\p_z
\left({\omega_\phi^n\over \omega_0^n}\right) \omega_0^n.\]
 Now we need to show that the remaining portion is
$o(\delta).\;$

 \[\begin{array}{lcl}
 \displaystyle \lim_{z\rightarrow z_0\in \p \Sigma}\;
\int_{\{z\}\times {\cal S}_\delta \setminus {\cal
S}_{\phi_0}}\;\p_z \left({\omega_\phi^n\over \omega_0^n}\right)
\omega_0^n & = & \int_{\{z_0\}\times {\cal S}_\delta \setminus
{\cal S}}\; {\eta^{\alpha}}_{, w_\alpha}(\omega_0)
\left({\omega_\phi^n\over \omega_0^n}\right) \omega_0^n\\
& = & \int_{\{z_0\}\times {\cal S}_\delta \setminus {\cal
S}}\;\left({{\p \eta^{\alpha}}\over {\p w_\alpha}} -
\Gamma^{\alpha}_{\beta\alpha}(\omega_0) \eta^\beta\right)
\left({\omega_\phi^n\over \omega_0^n}\right) \omega^n
\\
& = & \int_{\{z_0\}\times ({E}_\delta \bigcap {\cU_{\phi_0}}
)}\;\left({{\p \eta^{\alpha}}\over {\p w_\alpha}} -
\Gamma^{\alpha}_{\beta\alpha}(\omega_0) \eta^\beta\right)
\left({\omega_\phi^n\over \omega_0^n}\right) \det
\;g_{\alpha\b\beta} \det (J) \; d \;t
\\& = &  \int_{\{z_0\}\times ({E}_\delta \bigcap {\cU_{\phi_0}}
)}\;{{\p \eta^{\alpha}}\over {\p x_k}} \left( {{\p x_k}\over {\p
w_\alpha}} \; \det (J) \right) \det \;g_{\alpha\b\beta}  d
\;t\\
&& \qquad - \int_{\{z_0\}\times ({E}_\delta \bigcap {\cU_{\phi_0}}
)}\; \Gamma^{\alpha}_{\beta\alpha}(\omega_0) \eta^\beta
\left({\omega_\phi^n\over \omega_0^n}\right) \det
\;g_{\alpha\b\beta} \det (J) \; d\,t  \rightarrow 0
\end{array}
\] as $\delta \rightarrow 0.\;$ This is because all terms in the
last formula are uniformly bounded and the measure of $E_\delta$
tends to $0$ as $\delta \rightarrow 0.\;$ Here
\[d\,t = d \;t^1 d\;t^2\cdots d\,t^{2n}.\]

Consequently, we have shown that
\[
\int_{\Sigma\times M \setminus {\cS_{\phi_0}}}\; \; \mid {\cal D}
{{\p \phi}\over {\p z}}\mid^2_\phi\; \omega_\phi^{[n]} = \int_{\p
\Sigma}\; \;{{\b z}\over {|z|}}\; \int_{M\setminus
{\cS_{\phi_0}}}\;\p_z \left({\omega_\phi^n\over \omega_0^n}\right)
\omega_0^n -\int_{\p \Sigma}\;\;{{\b z}\over {|z|}}\; \int_M\;
{{\p \phi}\over {\p z}}\; (Ric(\omega_0) - \underline{R}
\omega_\phi)\wedge \omega_\phi^{[n-1]}.
\]
In other words, we have

\[\begin{array}{lcl}
\int_{\Sigma\times M \setminus {\cS_{\phi_0}}}\; \; \mid {\cal D}
{{\p \phi}\over {\p z}}\mid^2_\phi\; \omega_\phi^{[n]} & = &
\int_{\p \Sigma} \;\;{{\b z}\over {|z|}}\;  {\p\over {\p z}}
\bE(\phi(z,\cdot))\\
&= & \int_{\p \Sigma} \;  {\p E\over {\p {\bf n}}}. \end{array}.
\]
The theorem is then proved.
\end{proof}

If we replace the almost smooth solution by a partially smooth
solution, then
\begin{cor} \label{appl:kenergyissubharmonic1}Suppose that $\phi: \Sigma \rightarrow \overline{\cal H}$ is a
partially smooth solution described as in Definition
\ref{def:partiallysmooth}. Then the induced K energy function
$\bE: \Sigma\rightarrow {\bf R}$ (by $\bE(z) =
\bE(\phi(z,\cdot))$)  is a bounded weakly sub-harmonic function in
$\Sigma$ such that
\[
{{\p^2}\over {\p z\p\b z}}  \bE(\phi(z,\cdot)) \geq
\int_{M\setminus {\cal S}_{\phi_0}} | {\cal D} {{\p \phi}\over {\p
\b z}}|_{{\omega_{\phi}}}^2\, {\omega_{\phi}}^n\;\geq 0
\] holds in $\Sigma$ in the weak sense. Moreover,
\[
 \displaystyle \int_{\p
\Sigma}{{\p \bE}\over{\p\,{\bf n}}} (\phi(z, \cdot)) d\,s
 \geq  \int_{\Sigma\times M\setminus {\cal S}_{\phi_0}} | {\cal D}
{{\p \phi}\over {\p \b z}}|_{{\omega_{\phi}}}^2\,
{\omega_{\phi}}^n\, d\,s,
\] where  $d\,s$
is the length element of $\p \Sigma,\;$ and $\bf n$ is the outward
pointing unit normal direction at $\p \Sigma.\;$
\end{cor}

 \noindent In the case of {\bf extremal K\"ahler
metrics}, E. Calabi showed \cite{calabi85} that any extremal
K\"ahler metric must be invariant under some maximal compact
subgroup of the automorphism group. Consider all the metrics which
are symmetric under the same maximal connected compact subgroup of
the automorphism group. According to \cite{futaki92}, there exists
a unique  extremal holomorphic vector field \[ Y = Y^\alpha
{\p\over {\p w_\alpha}},
\]
 which is the gradient vector field of scalar curvature if the
metric is extremal. Note that this vector filed is unique in each K\"ahler class.  Consider
\[
  {\cal L}_Y \omega_\phi = \sqrt{-1}\p\b\p \theta(\phi).
\]
Here $\theta(\phi)$ is a real valued potential function for
$\phi.\;$ It is well known that one can modify the definition of
the K energy by this potential function such that the critical
point of the new functional is the extremal K\"ahler metric. Set
\[
\begin{array}{lcl}
  {{d\, \tilde{\bE}} \over {d\,t}} (\phi(t)) & = & - \displaystyle
  \int_M\; (R(\phi) - \underline{R} - \theta(\phi))
  {{\p \phi}\over {\p t}} \omega_\phi^{[n]}\\
  & = & {{d\, \bE} \over {d\,t}} (\phi(t))  +  \displaystyle
  \int_M\;  \theta(\phi) {{\p \phi}\over {\p t}} \omega_\phi^{[n]}.
  \end{array}
\]
Here $\theta(\phi)$ is the real valued potential function for this
holomorphic vector field.  This is true when all K\"ahler
potentials are invariant under the maximal compact subgroup.
\[
   {{\p \theta(\phi)}\over {\p w_\alpha}} = Y^{\b \beta}
   g_{\phi, \alpha \b \beta}.
\]
It is easy to see that \[ {{\p \theta} \over {\p t}} = Y({{\p
\phi}\over {\p t}}) = g^{\alpha\b\beta}_\phi {{\p \theta}\over {\p
w_\alpha}} {{\p^2 \phi}\over {\p t\p w_{\b\beta}}}.
\]
Then,\[
\begin{array}{lcl} {d\over {d\, t}}  \displaystyle
  \int_M\;  \theta(\phi) {{\p \phi}\over {\p t}} \omega_\phi^{[n]} & = &
\int_M\; \left({{\p \theta}\over {\p t}} {{\p \phi}\over {\p t}} +
\theta \; {{\p^2 \phi}\over {\p t^2}}  + \theta\; {{\p \phi}\over
{\p t}} \triangle_\phi ({{\p \phi}\over {\p
t}})\right) \omega_\phi^{[n]} \\
& = & \int_M\; \left(g^{\alpha\b\beta}_\phi {{\p \theta}\over {\p
w_\alpha}} {{\p^2 \phi}\over {\p t\p w_{\b\beta}}} {{\p \phi}\over
{\p t}} + \theta \; {{\p^2 \phi}\over {\p t^2}} + \theta\; {{\p
\phi}\over {\p t}} \triangle_\phi ({{\p
\phi}\over {\p t}})\right) \omega_\phi^{[n]} \\
& = & \int_M\; \left( - \theta\; \triangle_\phi({{\p \phi}\over
{\p t}}) - \theta \;{1\over 2} |\nabla {{\p\phi}\over {\p
t}}|^2_\phi + \theta \; {{\p^2 \phi}\over {\p t^2}} + \theta\;
{{\p \phi}\over {\p t}} \triangle_\phi ({{\p \phi}\over {\p
t}})\right) \omega_\phi^{[n]}\\ & = &
\int_M\;\theta(\phi)\left({{\p^2 \phi}\over {\p t^2}} - {1\over 2}
|\nabla {{\p\phi}\over {\p t}}|^2_\phi \right)\omega_\phi^{[n]}.
\end{array}
\]
By a similar calculation, we obtain
\[
{{\p^2\,\tilde{\bE}}\over {\p z\,\p \b z}} = {{\p^2\,{\bE}}\over
{\p z\,\p \b z}} + \displaystyle \int_M \; ({{\p^2\phi}\over {\p
z\p \b z}} - {1\over 2} |\nabla {{\p \phi}\over {\p z}}|^2_\phi)
\; \theta(\phi(\cdot, z)) \; \omega_\phi^{[n]}.
\]

 For an
almost smooth solution to HCMA equation (\ref{eq:hcma0}), the second
term vanishes completely. Note that if $\omega_\phi$ is uniformly
bounded from above, then $\theta(\phi)$ is uniformly Lipschitz.
This is a key technical step in  generalizing Theorem 6.1.1 to the case of extremal K\"ahler metrics.   Similarly, the same approximating proof will will work in this
case as well.  Thus \footnote{The crucial point is that
$\theta(\phi(z,\cdot))$ is uniformly Lipschitz.\; }.
\begin{cor}\label{appl:kenergyissubharmonic2} For any partially smooth
solution $\phi \in \b {\cal H} $ (c.f. Definition
\ref{def:partiallysmooth}) which is invariant under the maximal
compact subgroup, we have
\[
\int_{M\times \Sigma \setminus {\cal S}_{\phi_0}} | {\cal D} {{\p
\phi}\over {\p z}}|_{{\omega_{\phi}}}^2\, {\omega_{\phi}}^n\,
d\,\tau d\,\b z  \leq  \displaystyle \int_{\p \Sigma}{{\p
\tilde{\bE}}\over{\p\,{\bf n}}} (\phi) \;d\,s,
\]
where the left hand side is evaluated at points where K\"ahler
metric is smooth. Equality holds for any partially  smooth
solution. Moreover, $\tilde{E}$ is a bounded weakly sub-harmonic
function on $\Sigma.\;$
\end{cor}

\subsection{\label{appl:lowerboundkenergy}The  lower bound of the (modified) K energy}
In this subsection, we want to use Theorem
\ref{appl:kenergyissubharmonic} and Theorem
\ref{appl:kenergyissubharmonic1} to establish a lower bound for the
(modified) K energy.  In Sections 6-8, for simplicity, we always
deal with constant scalar curvature metric and K energy functional.
The corresponding treatment of extremal K\"ahler metric is parallel
and we leave for interested readers.
Now we set up some notations here first.\\

For any two K\"ahler potentials $\phi_0, \phi_1\in \cal H,$ we
want to use almost smooth solution to approximate the $C^{1,1}$
geodesic between $\phi_0$ and $\phi_1.\;$ For any integer $l, $
consider Drichelet problem for HCMA equation \ref{eq:hcma0} on the
rectangle domain $\Sigma_l = [-l,l] \times [0,1]$ with boundary
value as
\begin{equation}
\phi(s,0) = \phi_0, \phi(s,1) = \phi_1;\qquad \phi(\pm l,t) =
(1-t)\phi_0 + (1-t) \phi_1, \qquad (s,t) \in \Sigma_l.
\label{appl:boundarymap}
\end{equation}
 We
may modify this boundary map in the four corners so that the domain
is smooth without corner. Denote the almost smooth solution by
$\phi^{(l)}: \Sigma_l \rightarrow \cal H\;$  which corresponds to
this boundary map\footnote{We may need to alter the boundary value
slightly. }.  According to \cite{chen991}, $\phi^{(l)}$ has a
uniform $C^{1,1}$ upper bound which is independent of $l.\;$ Set
\begin{equation}
  \bE^{(l)}(s,t) = \bE(\phi^{(l)}(s,t)), \qquad\forall\; (s,t) \in
  \Sigma_l.
  \label{appl:kenergyondisk}
\end{equation}
Then, $E^{(l)}$ is a sequence of weakly sub-harmonic function with
uniform bound \footnote{because the decomposition formula for K
energy (\ref{appl:kenergy}) and the uniform $C^{1,1}$ bound on
$\phi^{(l)}.$} such that
\begin{equation}\bE^{(l)}(s,0) =
\bE(\phi_0) = A,  \qquad {\rm and} \;\; \bE^{(l)}(s,1) =
\bE(\phi_1) = B. \label{appl:kenergyondisk1}
 \end{equation}
  Now we are ready to prove (cf. Theorem \ref{intro:th:minimum})
\begin{theo}\label{appl:kenergyboundbelow} Any extremal K\"ahler metric realizes the
 absolute minimum of the modified K energy. Furthermore, the K\"ahler class is K-semistable if it admits a
constant scalar curvature metric.
\end{theo}
\begin{proof} We give a detailed proof for the case of constant scalar curvature metric.
The more general extremal K\"ahler metric case is similar and we
leave it for interested readers.  Suppose that $\phi_0$ is a
constant scalar curvature metric. Then ${{\p \bE^{(l)}} \over {\p
t}} = {{\p \bE^{(l)}} \over {\p s}} = 0$
when $t=0.\;$ Our theorem is reduced to the following claim:\\

\noindent {\bf Claim}:  $B = \bE(\phi_1) \geq \bE(\phi_0) = A.\;$ \\

Let ${\kappa}:(-\infty,\infty) \rightarrow {\bf R}$ be a smooth
non-negative function such that ${\kappa} \equiv 1 $ on $[ -
{1\over 2},{1\over 2}]$ and vanishes outside  of $[-{3\over
4},{3\over 4}].\;$  Set
\[
{\kappa}^{(l)} (s) = {1\over v} {\kappa}({s\over l}), \qquad {\rm
where} \;v = \int_{-\infty}^\infty\; {\kappa}(s)\;d\,s.
\]
Set
\[
f^{(l)}(t) = \displaystyle \int_{-\infty}^\infty\;
{\kappa}^{(l)}(s) \bE^{(l)}(s,t)\;d\,s.
\]
Then
\[\begin{array}{lcl}
f^{(l)}(0) & = & \displaystyle \int_{-\infty}^\infty\;
{\kappa}^{(l)}(s) \bE^{(l)}(s,0)\;d\,s = \displaystyle
\int_{-\infty}^\infty\; {\kappa}^{(l)}(s) A \;d\,s = A,
\\
f^{(l)}(1) & = &\displaystyle \int_{-\infty}^\infty\;
{\kappa}^{(l)}(s) \bE^{(l)}(s,1)\;d\,s = \displaystyle
\int_{-\infty}^\infty\; {\kappa}^{(l)}(s) B \;d\,s = B,
\\ {{d f^{(l)}}\over {d\,t}}\mid_{t=0} & = & 0.
\end{array}
\]
Now
\[
\begin{array}{lcl} \bE(\phi_1) - \bE(\phi_0) &= &  f^{(l)}(t)\mid_0^1
\\& = & \displaystyle \int_0^1\; \displaystyle\int_0^\theta \; {{d^2
f^{(l)}}\over
{d\,t^2}}\;d\,t\; d\,\theta \\
& = & \displaystyle \int_0^1\; \displaystyle\int_0^\theta
\;\displaystyle \int_{-\infty}^\infty {\kappa}^{(l)}(s)\; {{\p^2
\bE^{(l)}}\over {\p\,t^2}}\;d\,s\;d\,t\; d\,\theta\\ & = &
\displaystyle \int_0^1\; \displaystyle\int_0^\theta
\;\displaystyle \int_{-\infty}^\infty {\kappa}^{(l)}(s)\;
\triangle_{s,t}  \bE^{(l)}\;d\,s\;d\,t\; d\,\theta - \displaystyle
\int_0^1\; \displaystyle\int_0^\theta \;\displaystyle
\int_{-\infty}^\infty {\kappa}^{(l)}(s)\; {{\p^2 \bE^{(l)}}\over
{\p\,s^2}}\;d\,s\;d\,t\; d\,\theta\\ & \geq & - \displaystyle
\int_0^1\; \displaystyle\int_0^\theta \;\displaystyle
\int_{-\infty}^\infty {\kappa}^{(l)}(s)\; {{\p^2 \bE^{(l)}}\over
{\p\,s^2}}\;d\,s\;d\,t\; d\,\theta\\ & = & - \displaystyle
\int_0^1\; \displaystyle\int_0^\theta \;\displaystyle
\int_{-\infty}^\infty {{d^2 {\kappa}^{(l)}(s)}\over {d\,s^2}}\;
\bE^{(l)}(s,t)\;d\,s\;d\,t\; d\,\theta\\ & = & -{1\over l^2}
{1\over v} \displaystyle \int_0^1\; \displaystyle\int_0^\theta
\;\displaystyle \int_{-\infty}^\infty {{d^2 {\kappa}^{(l)}}\over
{d\,s^2}}\mid_{s\over l}\; \bE^{(l)}(s,t)\;d\,s\;d\,t\; d\,\theta

\end{array}
\]
Note that $|\bE^{(l)}(s,t)| $ has a unform bound $C.\;$ Thus,
\[
\begin{array}{lcl} {1\over l^2} {1\over
v} \mid \displaystyle \int_{-\infty}^\infty {{d^2
{\kappa}^{(l)}}\over {d\,s^2}}\mid_{s\over l}\;
\bE^{(l)}(s,t)\;d\,s  \mid & \leq & {1\over l^2} {1\over v}
\displaystyle \int_{-\infty}^\infty |{{d^2 {\kappa}^{(l)}}\over
{d\,s^2}}\mid_{s\over l}|\; \bE^{(l)}(s,t)\;d\,s
\\& \leq & {C\over {v\; l^2}}\;\displaystyle \int_{-\infty}^\infty |{{d^2
{\kappa}^{(l)}}\over {d\,s^2}}\mid_{s\over l}|\; \;d\,s \\
& = & {C\over {v\; l}}\;\displaystyle \int_{-\infty}^\infty |{{d^2
{\kappa}^{(l)}}\over {d\,s^2}}\mid_{s}|\; \;d\,s = {C\over {v\;
l}}\;\displaystyle \int_{-1}^1  |{{d^2 {\kappa}^{(l)}}\over
{d\,s^2}}\mid_{s}|\; \;d\,s\\
& \leq & {C\over l}
\end{array}
\]
for some uniform constant $C.\;$ Therefore, we have
\[
\begin{array}{lcl}
\bE(\phi_1) - \bE(\phi_0) & \geq & -{1\over l^2} {1\over v}
\displaystyle \int_0^1\; \displaystyle\int_0^\theta
\;\displaystyle \int_{-\infty}^\infty {{d^2 \kappa}\over
{d\,s^2}}\mid_{s\over l}\; \bE^{(l)}(s,t)\;d\,s\;d\,t\; d\,\theta
\\
&\geq & - \displaystyle \int_0^1\; \displaystyle\int_0^\theta
\;{1\over l^2} {1\over v}\;\mid \displaystyle
\int_{-\infty}^\infty {{d^2 {\kappa}}\over {d\,s^2}}\mid_{s\over
l}\; \bE^{(l)}(s,t)\;d\,s\mid \;d\,t\; d\,\theta
\\&\geq &- \displaystyle \int_0^1\; \displaystyle\int_0^\theta {C\over
l} \;d\,t\; d\,\theta = -{C\over {2 l}}.
\end{array}
\]
As $l \rightarrow \infty,$ we have
\[
\bE(\phi_1) \geq \bE(\phi_0).
\]
Since $\phi_1$ is an arbitrary K\"ahler potential,  the theorem is
then proved.
\end{proof}

\section{Partial regularity of the K energy minimzer}

\subsection{Strong convergence lemma for volume form}
In this subsection, we want to  prove 
\begin{theo} \label{appl:volumestrongconvergence}Suppose that $\{\phi_m, m \in \NN\}$ is a sequence of K\"ahler potentials in
$\cal H$ with uniform $C^{1,1}$ bound and suppose that $\phi_m
\rightarrow {\phi_0}  \in\b {\cal H}$ strongly in
$C^{1,\alpha}(\forall\;\alpha< 1)$ and weakly in $W^{2,p}$ for $p$
large enough.  If the corresponding sequence of K energies
$\{\bE(\phi_m), m \in \mathbb{N}\}$ is a cauchy sequence and \begin{equation}
\displaystyle \lim_{l\rightarrow \infty} \; \bE(\phi_l) \leq
\bE({\phi_0}),\label{appl:lowersemilimit}
\end{equation}
 then ${\omega_{\phi_m}^n\over
\omega^n} $ converges strongly to ${{\omega_{\phi_0}^n}\over
\omega^n}$  in $L^2(M,\omega).\;$
\end{theo}
\begin{proof} Set $f_m = {{\omega_{\phi_m}^n}\over \omega^n} \leq C $ and $g
= {{\omega_{{\phi_0}}^n}\over \omega^n}.\;$   Applying the
decomposition formula of the K energy
(\ref{appl:kenergy}), we obtain that $\{ \int_M f_m \log f_m, m \in \mathbb{N} \}$ is a Cauchy sequence. \\

Since $\{\phi_l, l \in \mathbb{N}\}$ weakly converges to ${\phi_0}$ in  the $W^{2,p}$
norm for $p$ large enough,  the lower order part of the K energy
converges to the corresponding lower order part of the K energy of
${\phi_0}.\;$ Thus
\begin{equation}
 \int_M\; f_l \log f_l \;\omega^n - \int_M \;g \log g
=  \bE(\phi_l) - \bE({\phi_0}) + o({1\over l}) \leq
 o( {1\over l}). \label{eq:strongvolumeconvergence}
\end{equation}
 Define
$F(u) = u \log u.\;$ For any $l$ large enough and for any
$\epsilon>0$, set $F(t) = F(t f_l + (1-t) (g+\epsilon)) = F(a\,t +
b)\;$ where
\[
  a =  f_l - g -\epsilon, \qquad {\rm and}\qquad b = g + \epsilon.
\]
Note that $a,b $ are both functions in $M.\;$ Clearly, we have
\[ |a| + |b| \leq C.\;\]
Note that
\[
F'(t) = a \log (a\,t + b) + a,
\]
and
\[
F''(t) = {{a^2}\over {a t + b}} \geq {a^2\over C}, \qquad \forall
\; t \in [0,1].
\]
Thus,
\[
\begin{array}{lcl} \int_M F'(0)\omega^n& =& \int_M\; (a \log b + a)\omega^n\\ & =  &
\int_M\; (f_l-g-\epsilon) \log (g_m +\epsilon) \omega^n   +
\int_M\; (f_l-g-\epsilon) \omega^n\\ & = & \int_M\;
(f_l-g-\epsilon) \log (g +\epsilon) \omega^n   + \int_M\;
\omega_{\phi_l}^n - \int_M\;
\omega_{{\phi_0}}^n -\epsilon \;vol(M) \\
& = & \int_M\; (f_l-g-\epsilon)\; \log (g+\epsilon) \;\omega^n
-o(\epsilon).
\end{array}
\]
Taking the following double limits
\[
\begin{array}{lcl}
\displaystyle\;\lim_{\epsilon \rightarrow 0} \displaystyle\;
\lim_{l\rightarrow \infty} \displaystyle\; \int_M\; F'(0) \omega^n
& = & \displaystyle\;\lim_{\epsilon\rightarrow 0}
\displaystyle\;\lim_{l\rightarrow \infty} \left(\int_M\;
(f_l-g-\epsilon)\; \log (g+\epsilon) \;\omega^n -o(\epsilon)
\right)\\
& = & \displaystyle\; \lim_{\epsilon \rightarrow 0} \left(\int_M\;
(g-g-\epsilon)\; \log (g+\epsilon) \;\omega^n -o(\epsilon) \right)
\\&
= & \displaystyle\; \lim_{\epsilon\rightarrow 0}\; o(\epsilon) =0.
\end{array}
\]
The second equality used the fact that $f_l \rightharpoonup g$
weakly in $L^p(M,\omega)$ and the 3rd equality used the fact that
$|g|$ is bounded.  Thus,
\[
\begin{array}{lcl} F(1)-F(0) & =& F(f_l) - F(g+\epsilon) \\
& = & \int_0^1\; F'(t) d\,t = F'(0) + \int_0^1 \int_0^t
\;F''(s)\;d\,s\,d\,t \\
& = & F'(0) +  \int_0^1 \int_0^t \;{{a^2}\over {a s + b}}
\;d\,s\,d\,t\\
&\geq & F'(0) +   \int_0^1 \int_0^t \;{{a^2}\over { C }}
\;d\,s\,d\,t =  F'(0)  + {a^2\over {2 C}}.
\end{array}.
\]
Integrating this over $M,\;$ we have,
\[\begin{array}{lcl} &&
\int_M\; (f_l-g-\epsilon)^2 \omega^n  =  \int_M \;a^2 \omega^n
\\ \qquad \qquad & \leq & 2 C \int_M\; (F(1) - F(0))\omega^n -  2 C \int_M F'(0)
\omega^n \\ \qquad \qquad & = & 2 C \left(\int_M f_l \log f_l \;
\omega^n - \int_M \;(g +\epsilon) \log(g+\epsilon) \omega^n\right)
-  2 C \int_M F'(0) \omega^n .
\end{array}
\]
Using inequality (\ref{eq:strongvolumeconvergence}), we have
\[
\begin{array}{lcl} &&
\int_M\; (f_l-g-\epsilon)^2 \omega^n \\
&\leq & C\left(\int_M g \log g \; \omega^n - \int_M \;(g
+\epsilon) \log(g+\epsilon) \omega^n + o({1\over l})\right) - C
\int_M F'(0) \omega^n
\\ &\leq & o(\epsilon + {1\over l}) -  C \int_M F'(0) \omega^n. \end{array}
\]
Consequently, we have
\[
\begin{array}{lcl} && 2 \displaystyle \lim_{l\rightarrow 0} \int_M\;(f_l-g)^2
\omega^n = 2 \lim_{\epsilon \rightarrow 0} \displaystyle
\lim_{l\rightarrow 0} \int_M\;(f_l-g-\epsilon +\epsilon)^2
\omega^n
\\&\leq &  \displaystyle \lim_{\epsilon \rightarrow 0} \displaystyle \lim_{l\rightarrow
\infty}\int_M\; (f_l-g-\epsilon)^2 \omega^n
  +  \displaystyle \lim_{\epsilon \rightarrow 0} \displaystyle \lim_{l\rightarrow \infty} \int_M\;\epsilon^2 \;\omega^n
\\& = &  \displaystyle \lim_{\epsilon\rightarrow 0}
\displaystyle \lim_{l\rightarrow \infty} \int_M\;
(f_l-g-\epsilon)^2 \omega^n \\
& \leq & \displaystyle \lim_{\epsilon\rightarrow 0} \displaystyle
\lim_{l\rightarrow \infty} \displaystyle  o({1\over l} +\epsilon)
+ C \displaystyle\;\lim_{\epsilon \rightarrow 0} \displaystyle\;
\lim_{l\rightarrow \infty} \displaystyle\; \int_M\; F'(0) \;\omega^n\\
& = & 0.
\end{array}
\]
Therefore,$ {{\omega_{\phi_l}^n}\over {\omega^n}} $ converges
strongly to ${{\omega_{{\phi_0}}^n}\over {\omega^n}}\; $ in
$L^2(M,\omega).\;$
\end{proof}

\subsection{Special properties of the K energy minimizer}
 We follow the setup in  Subsection \ref{appl:lowerboundkenergy}.
 Passing to a subsequence if necessary,
there exists a $C^{1,1}$ map $\underline{\phi_0}: \Sigma_1
\rightarrow \cal H$ such that
\begin{enumerate}
\item  $\phi^{(l)}$ converges to $\underline{\phi_0}$  weakly in
$W^{2,p} (\Sigma^1 \times M) $ for $p$ sufficiently large, with
respect to a fixed K\"ahler metric $\pi_2^*\omega + |d\,z|^2.$
 \item $\phi^{(l)}$ converges to $\underline {\phi_0}$ strongly in $C^{1,\alpha}$
for any $0< \alpha < 1.\;$ \item $\underline {\phi_0}(s,0)=\phi_0$
and $\underline {\phi_0}(s,1)=\phi_1.\;$
\end{enumerate}

The first key step in this subsection is to improve weak $L^p
(p>1) $ convergence  to a strong $L^2$ convergence for the volume ratio
 ${\omega_{\phi^{(l)}}^n\over \omega^n}.\;$\\

Recalled  the notation (\ref{appl:kenergyondisk})
\[
  \bE^{(l)}(s,t) = \bE(\phi^{(l)}(s,t)), \qquad\forall\; (s,t) \in
  \Sigma_l.
\]
As before, $E^{(l)}$ is a uniformly bounded, weakly sub-harmonic
function in $\Sigma^{(l)}$ with boundary condition
(\ref{appl:kenergyondisk1}).

 In this subsection, we assume that both of
$\phi_0$ and $ \phi_1 $ are  K\"ahler metrics with constant scalar
curvature in the fixed K\"ahler class.   Then,  ${{\p \bE^{(l)}} \over {\p t}} = {{\p \bE^{(l)}}
\over {\p s}} = 0$ when $t=0, 1.\;$ Theorem
\ref{appl:kenergyboundbelow} implies the following:
 \begin{equation}
A = B = \displaystyle \inf_{\phi \in {\cal H}}\; \bE(\phi)
\label{eq:Kenergyminimum1}
\end{equation}
and
\begin{equation}
\bE^{(l)}(s,t) \geq A, \qquad \forall\; (s,t)\in \Sigma_l.
\label{eq:Kenergyminimum2}
\end{equation}

\begin{lem} \label{appl:kenergyisapproharmonic}As $l\rightarrow \infty,$ the $L^1$ measure of the
Laplacian $\triangle_{s,t} \bE^{(l)}$ tends to $0$ in any fixed
compact subdomain.
\end{lem}
When there is no confusion, we will drop the superscript $(l)$.
\begin{proof}
 Let $\xi:(-\infty,\infty) \rightarrow \RR$ be a smooth
non-negative cut-off function such that $\xi \equiv 1 $ on $[ -
{1\over 2},{1\over 2}]$ and vanishes outside $[-{3\over 4},{3\over
4}].\;$
\[
  \begin{array}{lcl} \int_{t=0}^1 \;\int_{s=-{l\over 2}}^{l\over
2} \; |\triangle_{s,t}\; \bE(s,t)|\;d\,s\;d\,t & \leq &
\int_{t=0}^1 \;\int_{s=-{l}}^{l} \xi({s\over l}) \triangle_{s,t}\;
\bE(s,t)\;d\,s\;d\,t\\
&=& \int_{s=-l}^l {{\p \bE(s,t)}\over {\p t}}\mid_0^1\;
\xi({s\over l})\;d\,s - {1\over l} \int_{t=0}^1
\;\int_{s=-{l}}^{l}
\xi'({s\over l}) {{\p \bE(s,t)}\over {\p s}} \;d\,s\;d\,t\\
&=& 0 + {1\over l^2}  \int_{t=0}^1 \;\int_{s=-{l}}^{l}
\xi''({s\over l}) \bE(s,t) \;d\,s\;d\,t \\
&\leq & {1\over l^2}  \int_{t=0}^1 \;\int_{s=-{l}}^{l}
|\xi''({s\over l})|\cdot |\bE(s,t)| \;d\,s\;d\,t \leq {1\over l^2}
\int_{t=0}^1 \;\int_{s=-{l}}^{l} C \;d\,s\;d\,t\\
& = & {C\over l} \rightarrow 0.
\end{array}
\]
\end{proof}

\begin{lem}  \label{appl:kenergyisalmostconstant}For any point $(s,t)$ in a fixed compact domain in
$\Sigma^{(l)}$, except perhaps a set of measure $0$, we have $
\displaystyle\lim_{l\rightarrow \infty} \;\bE(s,t) = \displaystyle
\lim_{l\rightarrow \infty}\; \bE(\phi^{(l)}(s,t)) = A.\; $\\
\end{lem}

\begin{proof} Set $f^{(l)} = \triangle_{s,t} \; \bE^{(l)}(s,t)\geq 0.\;$ In $\Sigma_{l\over 2} \subset \Sigma_l$,
we have $ \displaystyle \lim_{l\rightarrow \infty}\displaystyle
\int_{\Sigma_{l\over 2}} \;f^{(l)} = 0. $ Next, we decompose
$\bE^{(l)}$ into two parts:
\[
    \bE^{(l)} = u^{(l)} + v^{(l)}, \qquad \;{\rm in}\qquad \Sigma_{l\over 2}
\]
such that
\[
\left\{ \begin{array}{ll} \triangle_{s,t} \; u^{(l)} \; = \; 0, & {\rm where}\;  u^{(l)}\mid_{\p \Sigma_{l\over 2}} =E^{(l)},\\
 \triangle_{s,t}\; v^{(l)} \; =\;f^{(l)} \geq 0, &\; {\rm where}\; v^{(l)}\mid_{\p \Sigma_{l\over 2}} =0.\end{array}
\right.
\]
It is clear that $v^{(l)}\leq 0.\;$   Since $E^{(l)}$ is uniformly bounded, then $u^{(l)}$ is a bounded
harmonic function in $\Sigma_{l\over 2}$ such that
\[
  u^{(l)}(s,0) =u^{(l)}(s,1) = A,\qquad \forall \; s\in [-{l\over
2},{l\over 2}].
\]
Taking limit as $l\rightarrow \infty, $ in any compact
subdomain $\Omega,$ we have $\displaystyle\lim_{l\rightarrow
\infty}\; u^{(l)} =A.\;$ Consequently,
\[\begin{array}{lcl} A &\leq &
\displaystyle\limsup_{l\rightarrow \infty} \;\bE^{(l)} \\
& = & \displaystyle\limsup_{l\rightarrow \infty}\;(u^{(l)} +
v^{(l)})
\\&\leq & \displaystyle\lim_{l\rightarrow \infty}\; u^{(l)} = A.
\end{array}\]
Therefore, for every point in $\Omega$(fixed), we have
\[
\displaystyle\lim_{l\rightarrow\infty} \bE^{(l)} =
\displaystyle\limsup_{l\rightarrow\infty} \bE^{(l)}= A.
\]

\end{proof}

Combining this with Theorem \ref{appl:volumestrongconvergence}, we
have
\begin{cor}\label{appl:volumestrongconvergence1}
 For any
$(s,t)\in \Sigma_{1}$, except perhaps  a set of measure $0$ in
$\Sigma_{1}$, the volume ratio ${{\omega_{\phi^{(l)}}^n}\over
{\omega^n}}$ converges strongly in $L^2(M,\omega)$ sense.
\end{cor} This corollary is
crucial in the following arguments.\\

For notational simplicity, set $\phi_l = \phi^{(l)} \in \cal H.\;$
 Let
$\zeta(s) = {1\over {1+s^2}}.\;$ Then
\[
\displaystyle \int_{z\in \Sigma_l}  \displaystyle \int_M\;
\omega^n \;\zeta(z) \;d\,z \wedge d \b z \leq C, \forall \;l\in
(1,\infty).
\]
Lemma \ref{appl:kenergyisapproharmonic} and Theorem
\ref{appl:kenergyissubharmonic} imply
\[
\int_{\Sigma_{l\over 2}}\;\int_M\; |{{\p \eta^{l,\alpha}}\over {\p
w_{\b \beta}}}|_{\phi_l}^2 \omega_{\phi_l}^n  d\,z\wedge d\,\b z\;
\leq \int_{z \in \Sigma_{l\over 2}} \triangle_{z} E(z,\cdot)
d\,z\wedge d \b z\; \rightarrow 0
\]
where
\[
\eta^{l,\alpha} = -\;g^{\alpha \b \beta}_{\phi_l}\; {{\p^2
\phi_l}\over {\p \;z\;\p w_{\b \beta}}},\qquad \qquad {\rm
and}\qquad g_{\phi_l,\alpha\b\beta} = g_{\alpha\b \beta} + {{\p^2
\phi_l}\over {\p w_\alpha \p w_{\b \beta}}}.
\]
When there is no confusion arised, we will drop the dependence on
$l.\;$\\

Now we adopte a sympletic point of view now: For any $l > 1,$  the product manifold
$\Sigma_{l}\times M$ is foliated by smooth holomorphic discs
which are transversal to $M$ direction, dictated by the structure of almost smooth solutions $\phi_l$ of HCMA
equation 1.1. In other words,  for any $z_0= t_0 + \sqrt{-1} s_0 \in \Sigma_1$ fixed,  we may use
$\{z_0\}\times M$ as the parametrization space of the set of
holomorphic discs which are  transversal to $\{z_0\} \times M.\;$ Note that this parametrization is effective
except a set of codimension 2.   Along each holomphic disc,
the $(n,n)$ form $\omega_{\phi}^n$ is invariant. The above two
inequalities can be re-stated as:
\begin{equation}
\int_{M}\; \int_{\Sigma_{l}}\; {{\omega^n}\over
{\omega_{\phi_l}^n}}\; \zeta(z) \;d\,z\wedge d\,\b z\wedge
\omega_{\phi_l}^n \leq C, \label{eq:Kenergyminimum3}
\end{equation}
and
\begin{equation}
\int_{M}\; \int_{\Sigma_{l\over 2}}\;|{{\p \eta^{l,\alpha}}\over
{\p w_{\b \beta}}}|_{\phi_l}^2 \;d\,\;z\wedge
d\,\b z\wedge \omega_{\phi_l}^n \rightarrow 0.
\label{eq:Kenergyminimum4}
\end{equation}

 Choosing any point $\;z_0$ in the interior of
$\Sigma_1$ such that ${\omega_{\phi_l(\;z_0,\cdot)}^n \over
\omega^n}$ converges strongly to
${\omega_{{\underline{\phi}_0}(\;z_0,\cdot)}^n \over \omega^n} $ in
$L^2(M,\omega)$.
For any $L^2$ function $h$ in $M$, we can normalize it by the following
\[
h(x) = \left\{\begin{array} {lcl} \displaystyle \lim_{\epsilon\rightarrow 0} \oint_{B_{\epsilon}(x)}\; h, & {\rm if \;limit\; exists,}
\\
0, & {\rm otherwise}.  \end{array} \right.\]
Then, $h(x)$
differs from the original function at most at a set of measure $0.\;$
Now, we decompose $\{z_0\}\times M$ into the union of two subsets $F_1
$ and $F_2$ such that
\[
{{\omega_{\underline{\phi_0}}}^n\over \omega^n} (\;z_0, x) >
0,\forall\;x \in\; F_1
\]
and
\[
{{\omega_{\underline{\phi_0}}}^n\over \omega^n} (\;z_0, x)  =
0,\forall\;x \in\; F_2.
\]
Clearly, $mes(F_1) > 0 $  since
\[
\begin{array}{lcl} mes(F_1) &  = &
\int_{F_1} \omega^n \\ & \geq &  c \int_{F_1}\;
{{\omega_{\underline{\phi_0}}}^n\over \omega^n} (z_0, x) \omega^n\\
& = & c \int_M\; {{\omega_{\underline{\phi_0}}}^n\over \omega^n} (\;z_0, x)
\omega^n = c\; vol(M)>0,\end{array}
\]
where ${1\over c} $ is the upper bound of the volume form ratio
${{\omega_{\underline{\phi_0}}}^n\over \omega^n} (z_0, \cdot) \;$
in $M.\;$\\

Ultimately, we
want to show that $mes(F_2)=0. \;$  This is not attainable at this point.
However, we can prove the following strong statement about the
volume form ratio in the limit.
\begin{theo}\label{appl:volumepositive} There exists a uniform constant $\epsilon_0$, which depends only on $\varphi_0, \varphi_1 \in \cal H$ (in particular,  independent of $x$),  such that, excluding at most
a set of measure $0$ from $F_1$, we have
\[ {{\omega_{\underline{\phi_0}}}^n\over \omega^n} (\;z_0, x) >
\epsilon_0, \qquad \forall\; x \in F_1.
\]
\end{theo}
\begin{proof} By our choice of $z_0$, ${{\omega_{\phi_l}}^n\over \omega^n} (z_0, x) \rightarrow
{{\omega_{\underline{\phi_0}}}^n\over \omega^n} (\;z_0, x) $
strongly in $L^2(M,\omega).\;$ For any $\delta>0,$ there exists an
open set $E_\delta$ with measure $E_\delta <{\delta\over 2}$ such
that ${{\omega_{\phi_l}}^n\over \omega^n} (z_0, x) \rightarrow
{{\omega_{\underline{\phi_0}}}^n\over \omega^n} (z_0, x) $
pointwisely in
$\{\;z_0\}\times (M\setminus E_\delta).\;$
Set $\ell_l(z): \{z_o\} \times M \rightarrow \{z\}\times M$ as ths syndromy map such that $(z_0,x) $ and
$(z, \ell_l(z) (x))$ lies in the same holomorphic disc.   Then, $\ell_l(z)$ is well defined for generic point
of $x\in M.\;$    Now set
\[S_l(z,x) = |{{\p \eta^{l,\alpha}}\over {\p w_{\b
\beta}}}|_{\phi_l}^2 (z, \ell_l(z)(x)) \qquad {\rm and}\qquad
f_l(z,x) = {\omega^n\over \omega_{\phi_l}^n}(z,\ell_l(z)(x)) \zeta(z).
\]
Denote
\[
S_l(x) = \int_{\Sigma_{l\over 2}}\; S_l(z,x)\;d\,z\wedge d\,\b z\;
\qquad {\rm and}\qquad f_l(x) = \int_\Omega f_l(z,x)\;d\,z\wedge d\,\b z.
\]

Then,
equations (\ref{eq:Kenergyminimum3}) and
(\ref{eq:Kenergyminimum4}) imply
\[
\int_M\; S_l(x) \cdot {\omega_{\phi_l}^n\over \omega^n}\;\omega^n
\rightarrow 0, \qquad {\rm and}\; \;
\int_M f_l(x) \cdot {\omega_{\phi_l}^n\over
\omega^n}\;\omega^n \leq C.
\]
The first assertion implies that $ \sqrt{S_l(x)} \cdot
\sqrt{{\omega_{\phi_l}^n\over \omega^n}} \rightarrow 0\;$ strongly
in $L^2(M)$. Consequently, $\{S_l(x)\}$ uniformly converges to
$0$ in $(M\setminus E_\delta)\bigcap F_1.\;$ On the other hand,
there exists a set $E'_\delta$ of measure at most ${\delta\over
2}$ such that
\[
\displaystyle\liminf_{l \rightarrow \infty}\;\left(f_l(x) \cdot
{\omega_{\phi_l}^n\over \omega^n}\right) < C(\delta),\qquad {\rm
whenever}\; x \in M \setminus E'_\delta.
\]
Let $F_\delta = F \subset
(E_\delta\bigcup E'_\delta).\;$ Then
\[
mes(F_\delta)\geq mes (F) -\delta.
\]
We proceed to prove that our theorem holds in $F_\delta.\;$\\

Let us pick any point $x_0\in F_\delta$ and fix it for now. Passing to a
subsequence if necessary, we may assume (w.l.o.g.) that
\[
{\omega_{\phi_l}^n\over \omega^n}(\;z_0,x_0) =
{\omega_{{\underline{\phi_0}}}^n\over \omega^n}(z_0,x_0) =
\epsilon
> 0.
\]
Here $\epsilon>0$ may be very small.  The goal is to show that a
uniform positive lower bound of the volume form ratio exists.  Clearly,
$S_l(x_0) \rightarrow 0$ for any fixed compact subset\footnote{As
a matter of fact, we may choose $\Omega=[ - {l\over 2},{l\over 2}]
$ and our claim still holds.} $\Omega.\;$  Next,
\[
  \displaystyle \liminf_{l\rightarrow \infty}\; \left( f_{l}(x_0) \cdot
{\omega_{\phi_l}^n\over \omega^n} (\;z_0,x_0)\right) < C(\delta).
\]
Passing to a subsequence if necessary, we have
\[
f_{l}(x_0) \leq C(\delta,\epsilon, x_0),
\]
in this subsequence. The main point is that it has a uniform bound
in terms of this subsequence.\\

 For any $l$, consider any holomorphic disc which passes through $(\;z_0,x_0)$ and
denote it by: \[ \ell_l: (\Sigma_{(l)},\p \Sigma) \rightarrow
(\Sigma_{(l)}\times M, (\p\Sigma_{(l)})\times M)\] such that
$l(z_0)=x_0.\;$  It is easy to see that this holomorphic disc
has uniformly bounded area in any fixed compact sub-domain (c.f.
Section 4). Let
\[
g_l(z) = \log {{\omega_{\phi_l}^n}\over
\omega^n}(\;z\;,\ell_l(z)),\qquad \forall \;z\; \in \Sigma_{(l)}.
\]
By definition,  $g_l$ is uniformly bounded in the boundary $\p \Sigma_l.\;$
Then
\[
g_l(z_0) = \log {{\omega_{\phi_l}^n}\over
\omega^n}(z_0,\ell_l(z_0)) =  \log {{\omega_{\phi_l}^n}\over
\omega^n}(z_0,x_0) = \ln \epsilon, \qquad \forall l.
\]
Moreover, there exists a uniform constant $C$ such that $|g_l(z)|
\leq C$ in the boundary $\p \Sigma_{(l)}.\;$  Since $\ell_l(\Sigma)$
is a measure $0$ set in $\Sigma\times M, $
 the limit of $g_l(z) (\;z\;\neq 0) $ most likely have no bearing
on  $\log {\omega_{\phi_0}^n\over \omega^n} $
in $\Sigma\times M.\;$ However,  we are interested in  obtaining a uniform positive
lower bound on $\epsilon$ through this procedure.
\\

Recall that Corollary \ref{local:linecurvature1} implies that
\[
\triangle_z\;\; g_l(z) = S(z\;,\ell_l(z)) + R_{0,\alpha\b\beta}
\eta^{l,\alpha} \eta^{l,\b\beta}\mid_{(z,\ell_l(z))}.
\]
Split  $g_l = u_l + v_l$ such that \[ \triangle_z\; v_l =
R_{0,\alpha\b\beta} \eta^{l,\alpha}
\eta^{l,\b\beta}\mid_{(z,\ell_l(z))}
\]
and \[ v_l(z) = g_l(z),\; \qquad\forall\;\;z\; \in \p \Sigma.
\]
\noindent {\bf Claim}: There exists a uniform constant $C$ such
that $|v_l|< C. $\\

Recall that Corollary \ref{local:deltazpotential} implies that
\[
- \triangle_z\; \phi_l(z,\ell_l(z)) = g_{0,\alpha\b\beta}
\eta^{l,\alpha}\;\eta^{l,\b\beta}, \forall \;z\; \in \Sigma_{l}.
\]
Note that $\phi_l(z,\ell_l(z))$ has a uniform bound in $\Sigma_{l}$.
There exists a constant $C$ such that
\[
   -C g_{0,\alpha\b\beta} < R_{0,\alpha\b\beta} < C
g_{0,\alpha\b\beta}.
\]

Consequently, \[
  \triangle_z\; (v_l - C\phi_l) <  0 < \triangle_z\; (v_l +
C\phi_l).
\]
By maximum principle, we have $|v_l|\leq C$ and  our earlier claim
 holds.  Next,
\[
\triangle_z\; \; u_l = S_l(z, \ell(z)) > 0,
\]
and $u_l\mid_{\p \Sigma_{l}}=0.\; $ Obviously $u_l \leq 0$
(maximum principle).  Moreover
\[\begin{array}{lcl}
\int_{\Sigma_{l}}\; e^{- u_l}(z)\cdot \zeta(z) d\,z \wedge d\,\b z
& \leq &  C \int_{\Sigma_{l}}\; e^{-log{\omega_{\phi_l}^n \over
\omega^n }}(z)\cdot
\zeta(z)  d\,z\wedge d\,\b z  \\
& = & \int_{\Sigma_{l}}\;\left({\omega^n\over
\omega_{\phi_l}^n}\right)(z,\ell_l(z))\cdot \zeta(z)\; d\,z\wedge
d\,\b z
\\& = & \int_{\Sigma_{l}} f_l(z,\ell_l(z)) \cdot \zeta(z) \; d\,z\wedge d\,\b z
\\& = & f_l(x_0) \leq C(z_0,\delta,\epsilon),
\end{array}\]
For any small positive number $\delta_1 > 0, $  Theorem
\ref{local:interiorboundofcurvature} implies
\[
0\leq \triangle_z\; u_l = S(z,\ell_l(z)) < {C\over \delta_1^2},
\qquad\forall\;z\; \in\;[-l+\delta_1,l-\delta_1]\times
[\delta_1,1-\delta_1].
\]
These two conditions imply that $u_l$ converges strongly in
$W^{1,\alpha}$ in any fixed compact sub-domain of $
(-l+\delta_1,l-\delta_1)\times (\delta_1,1-\delta_1).\;$
Moreover, \[ \int_{\Sigma_{l\over 2}}\;|\triangle_z\;u_l| =
\int_{\Sigma_{l\over 2}}\; S(z,\ell(z)) = S_l(x_0)\rightarrow 0.
\]

Passing to a subsequence if necessary,  there exists a
non-negative harmonic function  $u_\infty$ in
$(-\infty,\infty)\times [0,1]$ such that for any fixed compact
subset $\Omega$, we have
\[
u_l \rightharpoonup u_\infty
\]
weakly in $L^p(\Omega)$ for any $p > 1;$ and it converges strongly
to $u_\infty$ in $C^{1,\alpha} (0<\alpha < 1)$ in any fixed
compact sub-domain. Moreover, $u_\infty=0\;$
in the boundary of this infinite long strip $\Sigma_\infty = (-\infty, \infty)\times [0,1]$   and
\[
\int_{\Sigma_{\infty}}\; e^{ - u_l}(z)\cdot \zeta(z) d\,z\wedge d\,\b
z\; < C(\delta,\epsilon,x_0).
\]
The only harmonic function with this growth condition in $\Sigma_\infty$
is a constant function.  Thus,  $u_\infty\equiv 0.\;$ In particular, for any
fixed sub-domain $\Omega \subset (-l+\delta_1,l-\delta_1)\times
(\delta_1,1-\delta_1),\;$ we have $ u_l\rightarrow 0$ strongly. In
particular,
\[
\displaystyle \lim_{l\rightarrow \infty} \;u_l(z_0) = 0.
\]
It follows,
\[
\begin{array}{lcl}
\log \epsilon & =& g_l(z_0) \\
&= & u_l(z_0) + v_l(z_0) > -C.
\end{array}.
\]
Consequently,  we have
\[
{\omega_{{\underline{ \phi_0}}}^n\over \omega^n} (z_0,x_0) >
e^{-C} = \epsilon_0.
\]
Note that we chose $x_0$ arbitrarily in $F_\delta, $ thus the
lower bound $e^{-C} = \epsilon_0$ holds for every point in
$F_\delta.\;$ Since $mes (F\setminus F_\delta) < \delta, \;$ our
theorem is then proved.
\end{proof}

\subsection{A regularity theorem  for a $C^{1,1}$ minimizer  of the K energy functional and  the weak ``K\"ahler Ricci" flow}
In this subsection, we want to prove a regularity lemma for
any $C^{1,1}$ minimizer of the K energy in an  arbitrary K\"ahler class.
\begin{theo} \label{appl:volumeinw12}Suppose $\underline {\phi_0} \in C^{1,1}(M)$ is in the closure of $\cal H$ under weak $C^{1,1}$ topology.
If the K energy functional has a uniform lower bound in this K\"ahler class and $E(\underline{\phi_0})$ realizes the infimum of the K energy functional in this K\"ahler class,
then $\left({\omega_{\underline {\phi_0}}^n\over
\omega^n}\right)^{1\over 2}$ is in $W^{1,2} (M,\omega).\;$
\end{theo}
We want to prove this via the ``K\"ahler-Ricci" flow \footnote{This is not K\"ahler Ricci flow in the usual sense since we are not in the Canonical class and the first Chern class might not have a sign.}.
\begin{proof} Let ${\phi_0}(s) (0 \leq s \leq 1) $ be a 1-parameter
K\"ahler potentials such the following holds:
\begin{enumerate}
\item ${\phi_0}(0) = \underline {\phi_0}$ and ${\phi_0}(s)
(0<s\leq 1) \in {\cal H}.\;$ \item ${\phi_0}(s)$ has uniform
$C^{1,1}$ upper bound and ${\phi_0}(s) (s>0) \rightarrow
\underline {\phi_0}$ strongly in $W^{2,p}(M,\omega)$ for $p$ large
enough.

\end{enumerate}
Applying the ``K\"ahler-Ricci flow" to this 1-parameter family of
K\"ahler potentials ${\phi_0}(s) (0< s\leq 1):$
\begin{eqnarray}
  {{\p \phi(s,t)}\over {\p t}} & = & \log {\omega_{\phi}^n\over
\omega^n} \label{eq:krflow1}\\
\phi(s,0) & = & {\phi_0}(s).\;\end{eqnarray}

Clearly, for each $s>0$ fixed, there exists a uniform
$C^{2,\alpha}$ bound for $\phi(s,t) (s>0, 0\leq t\leq \infty).\;$
However, the upper-bound may depends on $s$ and in particular, it may
blows up when $t, s $ are both approaching $0$.  \\

\noindent {\bf Claim 1}: There exists a uniform constant $C$ which
is independent of the parameters $s$ and $t$  such that
\[ n +\triangle_0 \phi \leq C, \qquad \forall\; (s,t) \in \;(0,1]\times [0,\infty).
\]
Here we used $\triangle_\phi, \triangle_0$ to denote the
Lapalacian  operators of the K\"ahler metrics $\omega_\phi,
\omega$
respectively.\\

During this proof, we use $C$  to denote a generic constant which is independent of $s,t.\;$  Taking the derivatives of the flow
equation (\ref{eq:krflow1}), we obtain
\[
{{\p^2 \phi}\over {\p t^2}} = \triangle_\phi \left({{\p
\phi(s,t)}\over {\p t}}\right).
\]
This implies that \[
\begin{array}{lcl} {\omega_{\phi(s,t)}^n\over \omega^n} & = & e^{{\p
\phi(s,t)}\over {\p t}} \leq \displaystyle \max_{x\in M}\; e^{{\p
\phi(s,t,x)}\over
{\p t}}\\
& \leq & \displaystyle \max_{x\in M}\; e^{{\p \phi(s,t,x)}\over
{\p t}}\mid_{t=0} = \displaystyle \max_{x\in M}\;
{\omega_{\phi(s,0,x)}^n\over \omega^n}\\
&=& \displaystyle \max_{x\in M}\; {\omega_{{\phi_0}(s,x)}^n\over
\omega^n} \leq C.
\end{array}
\]
In other words, we have a uniform upper-bound on the evolved
volume
form.\\

Following the calculation in \cite{Yau78}, it is straightforward
to show (for each fixed $s>0$)
\[
\begin{array}{lcl} &&
(\triangle_\phi - {\p \over {\p t}}) \left(exp(-\lambda \phi)
(n+ \triangle_0\phi)\right) \\
& \geq & - exp(-\lambda \phi) (n^2 \displaystyle\inf_{i\neq 1}
(R_{i\b i 1 \b 1})) - \lambda\; exp(-\lambda \phi) (n -
\log {\omega_\phi^n \over \omega^n}) (n +\triangle_0 \phi) \\
&&\qquad+ (\lambda + \displaystyle\inf_{i\neq 1} (R_{i\b i 1 \b
1}) exp(-\lambda \phi) \cdot \left({\omega^n\over
\omega_\phi^n}\right)^{1\over n} (n + \triangle \phi)^{n\over
{n-1}},
\end{array}
\]
where $R_{i\b i 1\b 1}$ is the bisectional curvature of the
K\"ahler metric (corresponds to $\omega$) and $\lambda$ is a
positive number such that
\[
  \lambda + \displaystyle\inf_{i\neq 1} R_{i\b
i 1 \b 1} > 1.
\]
Multiplying $\left({\omega_\phi^n\over \omega^n}\right)^{1\over n}
$ on both sides, we get
\[
\begin{array}{lcl} & & \left({\omega_\phi^n\over \omega^n}\right)^{1\over n}
(\triangle_\phi - {\p \over {\p t}}) \left(exp(-\lambda
\phi) (n+ \triangle_0\phi)\right) \\
& \geq & - exp(-\lambda u) (n^2 \displaystyle\inf_{i\neq 1}
(R_{i\b i 1 \b 1})) \left({\omega_\phi^n\over
\omega^n}\right)^{1\over n}  - \lambda\; exp(-\lambda \phi) (n
-\log {\omega_\phi^n \over \omega^n}) \left({\omega_\phi^n\over \omega^n}\right)^{1\over n} (n +\triangle_0 \phi) \\
&&\qquad+ (C_0 + \displaystyle\inf_{i\neq 1} (R_{i\b i 1 \b 1}))
exp(-C_0 \phi) \cdot (n + \triangle \phi)^{n\over {n-1}},
\end{array}
\]
If $\phi$ is uniformly bound (independent of $s,t$), then (recall
that $\left({\omega_\phi^n\over \omega^n}\right)^{1\over n} < C$)
\[
\left({\omega_\phi^n\over \omega^n}\right)^{1\over n}
(\triangle_\phi - {\p \over {\p t}}) \;v \geq -c_1 -c_2 v + c_0
v^{n\over {n-1}}.
\]
where $c_0,c_1,c_2$ are uniform positive constants and $v =
exp(-\lambda\phi) (n +\triangle_0 \phi).\;$  Therefore,
\[
v(s,t) \leq v(s,0) \leq C.
\]
In other words, there is a uniform constant $C$ such that
\[
0\leq  n +\triangle_0\phi(s,t) \leq C,
\]
provided that $\phi(s,t)$ has  a uniform $C^0$ bound. It is easy
to see that we have uniform upper bound on ${{\p \phi}\over {\p
t}}.\; $  On the other hand, at minimum point of $\phi$, we have
\[
\log {{\omega_\phi^n} \over \omega^n} \geq 0.
\]
Thus,  in the barrier sense, we have
\[
{{\p \displaystyle \min_{M}  \phi}\over {\p t}} \geq 0.
\]

Consequently,
\[
|\phi(s,t)| \leq C .
\]
This concludes  the proof of our first claim above.\\

 For any sequence of number $s_i, t_i \rightarrow 0, $ set
\[
\phi_i =\phi(s_i,t_i),\qquad {\rm and}\qquad {\phi_0}_i =
\phi(s_i,0) ={\phi_0}(s_i,0).\]
 Passing to a subsequence if
necessary, we have that $\phi_i$ converges to some $C^{1,1}$
K\"ahler potential $\tilde{{\phi_0}}$ (strongly in $C^{1,\alpha}
(\forall\; \alpha < 1)$ and weakly in $W^{2,p}$ ($p$ large
enough)). Note that $\tilde{{\phi_0}}$ does not necessary equal to
${\phi_0}$ even though $t_i, s_i \rightarrow 0!$

\noindent {\bf Claim 2}: $\omega_{\tilde{{\phi_0}}}^n \;\equiv\;
\omega_{\phi_0}^n$ and ${\omega_{\phi_i}^n\over \omega^n} $
converge strongly to $ {\omega_{\tilde{{\phi_0}}}^n\over \omega^n}
$
in $L^2(M,\omega).\;$\\

To prove this claim, choose an arbitrary smooth non-negative cut
off function $\chi$ (fixed) and

\[\begin{array}{lcl} {1\over 2} {d\over {d t}} \displaystyle \int_M\; \chi
\left({\omega_\phi^n\over \omega^n }\right)^2\;\omega^n & = &
\displaystyle \int_M\; \chi \triangle_\phi
 \log {\omega_\phi^n\over \omega^n} \cdot \left({\omega_\phi^n\over \omega^n} \right)\;\omega_\phi^n
\\ & = & \displaystyle \int_M\; \chi \left(\triangle_\phi
 {\omega_\phi^n\over \omega^n} -  |\nabla \log {\omega_\phi^n\over \omega^n}|_\phi^2 \; {\omega^n \over \omega_\phi^n}
 \right)\;\omega_\phi^n\\ &\leq & \displaystyle \int_M\; \triangle_\phi \chi
 \left({\omega_\phi^n\over \omega^n} \right) \omega_\phi^n
\leq C. \end{array}\]

The last inequality holds since the evolving K\"ahler potentials
have a uniform $C^{1,1}$ upper-bound and $\chi$ is a fixed smooth
function. Integrating this inequality from $t=0$  to $ t = t_i$,
we have
\[
\begin{array}{lcl}
 \displaystyle \int_M\; \chi
\left({\omega_{\phi_i}^n\over \omega^n} \right)^2\;\omega^n
\mid_{s_i,t_i} & = & \displaystyle \int_M\; \chi
\left({\omega_{{\phi_0}_i}^n\over \omega^n}
\right)^2\;\omega^n\mid_{s_i,0} + \int_0^{t_i} \; {d\over {d t}}
\displaystyle \int_M\; \chi \left({\omega_\phi^n\over \omega^n}
\right)^2\;\omega^n\;d\,t
\\ &\leq & \displaystyle \int_M\; \chi \left({\omega_{{\phi_0}(s_i)}^n\over
\omega^n} \right)^2\;\omega^n  + C \;t_i.
\end{array}
\]
On the other hand,  ${\omega_{\phi_i}^n\over \omega^n}$ converges
weakly to ${\omega_{\tilde{{\phi_0}}}^n\over \omega^n}$ in
$L^2(M,\omega).\;$ Then
\begin{eqnarray}
\displaystyle \int_M\; \chi
\left({\omega_{\tilde{{\phi_0}}}^n\over \omega^n}
\right)^2\;\omega^n &\leq& \displaystyle \lim_{i \rightarrow
\infty}\; \displaystyle \int_M\; \chi
\left({\omega_{\phi_i}^n\over \omega^n} \right)^2\;\omega^n \label{eq:krflow2}\\
& \leq & \displaystyle\lim_{i\rightarrow \infty}
\left(\displaystyle \int_M\; \chi
\left({\omega_{{\phi_0}_i}^n\over
\omega^n} \right)^2\;\omega^n + C \; t_i\right) \label{eq:krflow3}\\
& = & \displaystyle \int_M\; \chi \left({\omega_{{\phi_0}}^n\over
\omega^n} \right)^2\;\omega^n. \label{eq:krflow4}\end{eqnarray}
The last equality holds since ${\phi_0}_i$ converges strongly to
${\phi_0}$ in $W^{2,p}(M,\omega)$ for $p$ large enough (by our
assumption in the beginning).  This holds for any non-negative
smooth cut off function in $M.\;$ Consequently, we have
\[
0\leq {\omega_{\tilde{{\phi_0}}}^n\over \omega^n}  \leq
{\omega_{{\phi_0}}^n\over \omega^n}
\]
a.e. in $M.\;$  However,
\[
\displaystyle \int_M\;{\omega_{\tilde{{\phi_0}}}^n\over \omega^n}
\omega^n =  \displaystyle \int_M\; {\omega_{{\phi_0}}^n\over
\omega^n} \omega^n = vol(M)!
\]
Consequently,
\begin{equation}\omega_{\tilde{{\phi_0}}}^n \equiv
\omega_{\phi_0}^n \label{appl:eq:volumeequal}
\end{equation} in the
sense of $L^2(M,\omega).\;$The uniqueness of $C^{1,1}$ solution to
the Monge-Ampere equation  implies that $\tilde{{\phi_0}} =
{\phi_0}.\;$ In particular, this implies the K energy $E(\phi_i)$
converges to $\bE({\phi_0}).\;$

On the other hand, the equality (\ref{appl:eq:volumeequal}) forces
equality in (\ref{eq:krflow2})-(\ref{eq:krflow4}) to hold. In
particular, we have
\[
\displaystyle \int_M\; \chi
\left({\omega_{\tilde{{\phi_0}}}^n\over \omega^n}
\right)^2\;\omega^n  = \displaystyle \lim_{i\rightarrow
\infty}\displaystyle \int_M\; \chi \left({\omega_{\phi_i}^n\over
\omega^n} \right)^2\;\omega^n.
\]
Thus, ${\omega_{\phi_i}^n\over \omega^n}$ converges strongly to
${\omega_{\tilde{{\phi_0}}}^n\over \omega^n}$ in
$L^2{(M,\omega)}.\;$ Our
second claim is then proved.\\

Now we used these two claims to prove our theorem.   Consider
$\bE(s,t) = \bE(\phi(s,t)) (0< s \leq 1 \;\;{\rm and}\;\; 0\leq t
< \infty).\;$ Set
\[
A = \displaystyle \inf_{\phi\in {\cal H}}\; \bE(\phi)
> -\infty.
\]
By our assumption on ${\phi_0}, $ we have
\[
 \bE({\phi_0}) = A = \lim_{s\rightarrow 0} \; \bE(\phi(s,0)) \leq
\bE(s,t), \qquad \forall\; s>0, t\geq 0.
\]
For any fixed number $c_0>0,$ it is straightforward to show that
there exists a sequence $s_i,t_i \rightarrow 0$ such that
\[
\begin{array} {lcl} - c_0 & \leq &
{{\p \bE(s,t)}\over {\p t}} \mid_{s_i,t_i} \\
& = & \int_M \log {\omega_\phi^n\over \omega^n} \; \triangle_\phi
{ {\p\phi}\over {\p t}} \omega_\phi^n\mid_{s_i,t_i} - \int_M\; {
{\p\phi}\over {\p t}} (Ric(\omega_0)-\underline{R}\;
\omega_\phi)\wedge
\omega_\phi^{n-1}\mid_{s_i,t_i}\\
& = &- \int_M \; |\nabla \log {\omega_\phi^n\over
\omega^n}|_\phi^2 \; \omega_\phi^n \mid_{s_i,t_i} - \int_M\;
\left(\log {\omega_\phi^n\over \omega^n} - C\right) (Ric(\omega_0)
- \underline{R}\; \omega_\phi) \wedge
\omega_\phi^{n-1}\mid_{s_i,t_i} \\& = &- \int_M \; |\nabla \log
{\omega_\phi^n\over \omega^n}|_\phi^2 \; \omega_\phi^n
\mid_{s_i,t_i} - c \int_M\; \left(\log {\omega_\phi^n\over
\omega^n} - C\right) \omega \wedge
\omega_\phi^{n-1}\mid_{s_i,t_i}.
\end{array}
\]
Here $c, C$ are some uniform positive number such that \[ \log
{\omega_\phi^n\over \omega^n} < C + 1 \] and \[
Ric(\omega_0)-\underline{R}\;\omega_\phi \leq c \;\omega.
 \]
Thus
\[\begin{array}{lcl}
\int_M \; |\nabla \log {\omega_{\phi_i}^n\over
\omega^n}|_{\phi_i}^2 \; \omega_{\phi_i}^n & \leq & c_0 - c
\int_M\; \left(\log {\omega_{\phi_i}^n\over \omega^n} -
C\right) \left(\omega + \sqrt{-1} \p \b \p \phi_i - \sqrt{-1} \p \b \p \phi_i  \right) \wedge \omega_{\phi_i}^{n-1}\\
&\leq & c_0 + c \int_M\; \left(\log {\omega_{\phi_i}^n\over
\omega^n} - C\right) \sqrt{-1} \p \b \p \phi_i \wedge
\omega_{\phi_i}^{n-1} - c \int_M\; \left(\log
{\omega_{\phi_i}^n\over \omega^n} - C\right)
\omega_{\phi_i}^n\\
&\leq & c_0 - c\int_M\; \sqrt{-1} \;\p \log
{\omega_{\phi_i}^n\over \omega^n} \wedge  \b \p \phi_i \wedge
\omega_{\phi_i}^{n-1} - c \int_M\; \left(\log
{\omega_{\phi_i}^n\over \omega^n} -
C\right) {\omega_{\phi_i}^n\over \omega^n} \;\omega^n\\
&\leq & c_0 + c \left(\epsilon \int_M \; |\nabla \log
{\omega_{\phi_i}^n\over \omega^n}|_{\phi_i}^2\; \omega_{\phi_i}^n
+ {1\over \epsilon} \int_M\; \sqrt{-1} \p \phi_i \wedge \b \p
\phi_i \wedge \omega_{\phi_i}^{n-1}\right) + C
\\& \leq & C(\epsilon) + c\epsilon \int_M \; |\nabla \log
{\omega_{\phi_i}^n\over \omega^n}|_{\phi_i}^2 \;
\omega_{\phi_i}^n.
\end{array}
\]
Choose $\epsilon $ small enough so that $c\;\epsilon < {1\over
2}.\;$ With this $\epsilon$, we have
\[
\begin{array}{lcl} \int_M\; |\nabla \sqrt{\omega_{\phi_i}^n \over
\omega^n} \mid^2_\omega \;\omega^n & \leq  & C  \int_M\; |\nabla
\sqrt{\omega_{\phi_i}^n \over \omega^n}|^2_{\phi_i} \;\omega^n \\
&= &
  C \int_M |\nabla \log {\omega_{\phi_i}^n\over
\omega^n}|_{\phi_i}^2 \omega_{\phi_i}^n \leq C.
\end{array}
\]
Letting $i\rightarrow \infty$, we see that
\[
\int_M\; |\nabla \sqrt{\omega_{\tilde{{\phi_0}}}^n \over \omega^n}
\mid^2_\omega \;\omega^n \leq C.
\]
Since ${\omega_{\tilde{{\phi_0}}}^n \over \omega^n} \equiv
{\omega_{{{\phi_0}}}^n \over \omega^n}, $ we have
\[
\int_M\; |\nabla \sqrt{\omega_{{{\phi_0}}}^n \over \omega^n}
\mid^2_\omega \;\omega^n \leq C.
\]
The theorem is then proved.
\end{proof}

\section{The Problem of uniqueness of extremal K\"aher metrics}

 Following notations in Subsection \ref{appl:lowerboundkenergy}. Suppose $\phi_1$ is also a constant scalar curvature metric. Then
 \[
 E(\phi_0) = E(\phi_1) = A.
 \] Since
 \[
    \mid \p \b \p \phi^{(l)}\mid_{\Sigma^{(l)} \times M} \leq C
 \]
holds uniformly (independent of $l$), there exists a subsequence of
$\varphi^{(l)}$which converges to $\underline{\phi} \in \b \cH$ in
the weak $C^{1,1}$ norm.
 Following Theorem \ref{appl:kenergyboundbelow}, for any point
 $(s,t) \in
{\Sigma_1}^0, $ we have \[ \displaystyle \lim_{l\rightarrow
\infty} \bE(\phi^{(l)}(s,t)) = A = \displaystyle \inf_{\phi\in
{\cal H}}\; \bE(\phi).\]
In the discussion below, we fix an aribitrary interior point $(s,t) \in {\Sigma_1}^0.\;$
 Theorem
\ref{appl:volumestrongconvergence} implies that
${\omega_{\phi^{(l)}}}^n(s,t,\cdot) $ converges strongly to
$\omega_{\underline{\phi}}^n(s,t,\cdot).\;$ By Theorem
\ref{appl:volumepositive}, we have
${\omega_{\underline{\phi}}^n(s,t,\cdot)\over \omega^n} >
\epsilon_0$ as long as it is positive.  The set of points
where the volume ratio
${\omega_{\underline{\phi}}^n(s,t,\cdot)\over \omega^n} $ vanishes
must have measure $0.\;$ Otherwise,  it contradicts with the fact
that $\sqrt{{\omega_{\underline{\phi}}^n(s,t,\cdot)\over
\omega^n}}$ is in $W^{1,2}(M,\omega) $ (cf. Theorem
\ref{appl:volumeinw12}). Thus,
$${\omega_{\underline{\phi}}^n(s,t,\cdot)\over \omega^n} > \epsilon_0$$
for all points in $M$ except at most a set of measure $0.\;$
Normalizing this volume ratio  in the $L^2$ sense, we obtain
\[
{\omega_{\underline{\phi}}^n(s,t,x)\over \omega^n} >
\epsilon_0, \forall\; x\in M.
\]
 Since ${\underline{\phi}}$
has uniform $C^{1,1}$ bound, this implies that \[
\omega_{\underline{\phi}}(s,t,\cdot)
> c_0\;\omega
\]
for some positive constant $c_0.\;$ In other words, the metric
$g_{\underline{\phi}}$ is equivalent to $g_0.\;$  For any locally
supported test function $\xi, $ we have\footnote{In any fixed open
set, $\omega_{\underline{\phi}}$ can be approximated by a sequence
of smooth K\"ahler metrics such that all K\"ahler metrics in the
sequence has a uniform positive lower bound. Thus one can do small
deformations in the arbitrary direction. Consequently, one can
establish the Euler-Lagrange equation in the weak sense. }
\[
\int_M\; \log {{\omega_{\underline{\phi}}^n}\over {\omega^n}}
\sqrt{-1}\p \b \p \xi \wedge \omega_{\underline{\phi}}^{[n-1]} =
\int_M\; \xi (Ric(\omega_0)-\omega_{\underline{\phi}})\wedge
\omega_{\underline{\phi}}^{[n-1]}.
\]
Write $$\omega_{\underline{\phi}} = \left(g_{\alpha\b \beta} +
{{\p^2 {\underline{\phi}}}\over {\p w_\alpha \p
w_{\b\beta}}}\right) d\,w^\alpha \,d\,w^{\b \beta} =
g_{{\underline{\phi}},\alpha\b \beta}\; d\,w^\alpha \,d\,w^{\b
\beta}, $$ and
$$ f =  {{\omega_{\underline{\phi}}^n}\over {\omega^n}}.\;$$ Since $\log f$ is in
$W^{1,2}, $ we then have
\[
-\int_M\; g^{\alpha\b\beta}_{\underline{\phi}} {{\p\; log f}\over {\p w_\alpha}} {{\p \xi}\over {\p  w^{\b \beta}}}
 \; f= -\int_M\; \xi
(Ric(\omega_0)-\omega_{\underline{\phi}})\wedge
\omega_{\underline{\phi}}^{n-1}
\]
for any locally supported test function $\xi.\;$ Hence $\log f$
satisfies the following 2nd order non-linear equation weakly:
\[
 {1\over f}\; {\p \over {\p w_{\b \beta}}} \left( g^{{\underline{\phi}},\alpha\b\beta} f {\p \over {\p w_\alpha}} \log f \right)
= g^{{\underline{\phi}},\alpha\b\beta}
Ric(\omega_0)_{\alpha\b\beta} - n.
\]
Note that this is a uniformly elliptic second order non-linear
partial differential equation with uniformly bounded coefficients,
while the right hand side is in $L^\infty.\; $ According to the
Holdier estimate (due to De  Giorgi), there exists a small
constant $\alpha \in (0,1)$ such that $\log f \in C^{\alpha}$ for
any interior points. Since $\p M = \emptyset, $ this implies that
$f$ is $C^{\alpha}(M).\;$ Appealing  to the Monge-Ampere equation
\[
  \det\; \left(g_{\alpha\b \beta} + {{\p^2 {\underline{\phi}}}\over
{\p w_\alpha \p w_{\b\beta}}}\right) = f,
\]
it follows that ${\underline{\phi}}\in C^{2,\alpha}.\;$ Returning
to the original equation of divergence form, we have
\[
   g^{\alpha\b\beta}_{\underline{\phi}} {{\p^2}\over {\p w_\alpha \p \b
w_\beta}} \left(\log f \right) = \left(g^{\alpha\b
\beta}_{\underline{\phi}} Ric_{\alpha\b \beta}(\omega) - n\right)
f.
\]
Here the left hand side is a uniformly elliptic operator with
$C^\alpha$ coefficients, the right hand side is also $C^\alpha$
continuous.  The standard elliptic regularity theory  implies that
$\log {{\omega_{\underline{\phi}}^n}\over {\omega^n}} \in
C^{2,\alpha}.\;$ This in turn implies that $\phi \in C^{4,\alpha}
$ or the right hand side is in $C^{2,\alpha}\cdots.\;$ Repeated
boot-strapping between these two equations shows  that
${\underline{\phi}}$ is smooth. Consequently, it must be of
constant scalar curvature.  It is easy to see ${{\p
{\underline{\phi}}}\over {\p s}} = 0,$ and
${\underline{\phi}}(s,t) (0\leq t\leq 1)$ satisfies the geodesic
equation:
\[{{\p^2 \underline{\phi}}\over {\p t^2}} -
g_{\underline{\phi},\alpha\b \beta} \;\eta^\alpha \;\eta^{\b
\beta} = \triangle_z \underline \phi -
g_{\underline{\phi},\alpha\b \beta} \;\eta^\alpha \;\eta^{\b
\beta} = 0,
\]
where \[ \eta^\alpha = - g^{\alpha \b \beta}_{\underline{\phi}}\;
{{\p^2{\underline{\phi}}}\over {\p w_\alpha \p \b w_\beta}}.
\]
The second variation of the K energy must be identically $0$ in
the direction of $t,$ which implies
 \[ \int_{0}^1\;d\,t\int_M\; |{{\p
\eta^\alpha}\over {\p w_{\b \beta}}}|^2_{\underline{\phi}}
\omega_{\underline{\phi}}^{[n]} = 0
\]
or \[ {{\p \eta^\alpha}\over {\p w_{\b \beta}}} \equiv 0,\qquad
\forall\;\alpha,\beta = 0, \forall\; 1,2,\cdots n.
\]
Thus, this path represents a path of holomorphic transformation.
The uniqueness is then proved.

\begin{theo} In any K\"ahler class,  the extremal K\"ahler metric is
unique up to holomorphic transformation.
\end{theo}
This concludes our proof of Theorem \ref{th:uniqueness}.

\section{Appendix:  Loop space of $Gl(n,c)$ and holomorphic discs}

The purpose of this appendix is to give a proof for Lemma
\ref{defo:codimension}.  The lemma is more or less known to
experts in the field, although it is difficult to find exact
statement in literatures.  The proof presented here is shown to us
by Professor E. Lupercio. It uses some standard theory of loop
groups.   We will be very explicit in our presentation here for
the sake of completeness.

To simplify the following explanation, let us suppose that
$G=\GLC{n}$. The \emph{loop group} $\LL G$ of a Lie group $G$ is the
space of maps from the unit circle $S^1$ in
$\complex$ to the corresponding group $G$. 
In this note,  the space $\LL G$ is endowed with the structure of an
infinite dimensional \emph{polarized} manifold. By a
\emph{polarization} of a vector space $H,\;$ we mean a class of
decompositions $H_+\oplus H_-$ that differ only ``by a finite
amount." 
 A manifold is polarized if its tangent bundle is polarized at every
fiber.

There are several important subgroups of the loop group that deserve
consideration. The first of them is the subgroup $\LGp$ of loops
$\gamma \in \LG$ that extend to holomorphic maps of the closed unit
disc $D^2$ on the complex plane $\tilde{\gamma} : D^2 \to \GLC{n}$.

The loop group has very important homogeneous spaces that posses
very nice geometrical interpretations. The most important of them is
the \emph{restricted Grassmannian of a Hilbert space $H$}. The
fundamental idea is that the loop group acts transitively on the
restricted Grassmannian, and the stabilizer of a point is the
subgroup $\LGp$. This action thus realizes the restricted
Grassmannian as a homogeneous space for the loop group of the form
$\LG / \LGp.\;$ To define this Grassmannian, we need the concept of
\emph{Polarization} for the Hilbert space.  Let just say that if we
realize our Hilbert space as the space of functions on the circle
$H^{(n)}=L^2(S^1,\complex^n)\cong H\otimes \complex^n$, then the
natural polarization for $H^{(n)}$ is given by
\begin{equation}
   H^{(n)}=H^{(n)}_+\oplus H^{(n)}_-
\end{equation}
where $H^{(n)}_+$ consists of those elements of $H^{(n)}$ that are boundary values
for a holomorphic map on the unit disc $D^2$, and $H^{(n)}_-$ is the orthogonal
complement of  $H^{(n)}_+$ in $H^{(n)}$. In other words $H^{(n)}_+$ is the space of functions $f(z)$ so that in its Fourier expansion no negative powers of $z$ appear.

We define the \emph{restricted Grassmannian} $\Grrinf$ of $H^{(n)}$
to be the space of all closed subspaces $W$ of $H^{(n)}$ so that the
projections $W\to H^{(n)}_+$ and $W\to H^{(n)}_-$ are respectively a
Fredholm and a Hilbert-Schmidt operator. This definition is crafted
in such a way  that $W$ is then `comparable' in a suitable sense
with $H^{(n)}_ .\;$ That is to say, the decomposition $W\oplus
W^\bot$ is also a
polarization. 
With this definition, the restricted Grassmannian is an infinite
dimensional complex manifold with charts modeled on the Hilbert
space $\II_2 (W; W^\bot)$ of Hilbert-Schmidt operators $W\to
W^\bot$. This shows that the restricted Grassmannian group is a
polarized manifold.

More relevant to this discussion is the Grassmannian $\Grrr{n}  \subseteq \Grrinf$
consisting of those $W\in\Grrinf$ so that $zW\subseteq W$. The index of the
projection $W\to H^{(n)}_+$ is called the virtual dimension of $W$. The \emph{virtual
dimension} index the connected components of the Grassmannian, that is to say that it
can be thought as an isomorphism $\pi_0 \Grrinf \cong \integer$.

If we let the loop group $\LG$ acts on  $H^{(n)}$ by matrix
multiplication (every element $\gamma(z)$ of the loop group is a
matrix valued function on $S^1$), then the action induces a
corresponding action on $\Grrr{n}$ -- this is the purpose of the
condition $z W \subseteq W$ in the definition of this Grassmannian.
The action is transitive and the isotopy group of $H^{(n)}_+$ is
precisely $\LGp$. This produces the identification $\LG / \LGp \cong
\Grrr{n}$. A version of the maximum modulus principle furthermore
implies that
\begin{equation}\label{Birkhoff}
   \Omega U_n \cong \LL U_n / U_n \cong \LG / \LGp \cong \Grrr{n}.
\end{equation}

There is natural stratification of $\Grrr{n}$ whose strata are indexed by homomorphisms $S^1 \to\GLC{n}$. Every such homomorphism can be written in the form
$$ z^{\mathbf{k}}=\left( \begin{array}{ccccc}
               z^{k_1} & & & & \\
               & z^{k_2} & & & \\
               & & z^{k_3} & & \\
               & & & \ddots & \\
               & & & & z^{k_n}
               \end{array} \right)
$$
where $ {\mathbf{k}} = (k_1,k_2,k_3,\ldots, k_n)$ is an integer
partition of the non-negative integer number $k$, namely $k_1 + k_2
+ \cdots k_n = k$.

The \emph{Birkhoff factorization theorem}\footnote{The Birkhoff factorization theorem is equivalent to the theorem of Grothendieck that states that  every holomorphic bundle of rank $n$ over the Riemman sphere can be uniquely written in the form of powers of the Hopf bundle $O(k_1)\oplus\cdots\oplus O(k_n)$.} establishes that any loop $\gamma(z) \in \LG$ can be factored in the form
\begin{equation}\label{factorization}
\gamma(z) = \gamma_-(z) z^{\mathbf{k}} \gamma_+(z)
\end{equation}
where $\gamma_+(z) , \gamma_-(1/z) \in \LGp$ and ${\mathbf{k}}$ is well defined up to the ordering of the $k_i$'s.

We will say that ${\mathbf{k}}$ is the \emph{multi-index} (or Grothendieck index) of $\gamma(z)$ whenever equation \ref{factorization} holds. We will tolerate the ordering ambiguity in this definition.

Let us return to the description of the stratification of $\Grrr{n}$.  Notice that given a loop $\gamma(z) \in \LG$ of index $\kk$, then multiplying on the right by any element $\phi_+(z) \in \LGp$ will not affect the multi-index, that is, $\gamma(z) \phi_+ (z)$ still has multi-index $\kk$. >From this we conclude that the multi-index is constant along orbits of the right-action of $\LGp$ in $\LG$. In other words every point of $\Grrr{n}$ has a well defined multi-index.

Let us define $\LGm$ by declaring that $\phi_-(z) \in \LGm$ if and only if $\phi_-(1/z) \in \LGp$. Again the action of $\LGm$ doesn't affect te multi-index of an element. Then the  orbits of the action of $\LGm$ in $\Grrr{n}$ are precisely the same as the sets of elements in with the same multi-index and all its permutations. This is once more a consequence of the Birkhoff factorization. To avoid the problem of the permutations we will have to consider a smaller subgroup $N^-$ of $\LGm$. The group $N^-$ consists of those elements in $\gamma_-\in\LGm$ so that,
$$ \gamma_-(\infty) =\left( \begin{array}{cccc}
               1 & * & * & * \\
               0 & 1 & * & * \\
               0 & 0 &  \ddots & * \\
               0 & 0 &  0 & 1
               \end{array} \right).
$$

For every partition ${\mathbf{k}}$ of $k$ (and here the order is
important) we define the subspaces\\
 $H_\kk \in \Grrr{n}$ as
\begin{equation}\label{basicspaces}
H_\kk = z^\kk H^{(n)}_+ = \{ f(z)= ( f_1(z), \ldots, f_n(z) ) \colon f_i(z) = \sum_{j=k_i}^{\infty} a^i_j z^j ,\ a_j^i \in \complex \}.
\end{equation}
and we define $\Sigma^\sigma_\kk$ to be the orbit of $H_\kk$ under the action of $\LGm$ in $\Grrr{n}$.

Define now $\Sigma_\kk$ to be $N^- \cdot H_\kk \subset \Grrr{n}$
(again, the order in $\kk$ in important.).  We have

\begin{prop} The set of all elements in $\Grrr{n}$ of multi-index $\kk$ is precisely $\Sigma_\kk^\sigma.\;$ \end{prop}

\begin{proof} Take any element $\gamma(z) \in \LG$ of index $\kk$, and write $W= \gamma(z) H^{(n)}_+ \in \Grrr{n}$. Using the Birkhoff
factorization $\gamma(z) = \gamma_-(z) z^{\mathbf{k}} \gamma_+(z),\;
$ we have $W= \gamma_-(z) z^{\mathbf{k}} \gamma_+(z)  H^{(n)}_+ =
\gamma(z) = \gamma_-(z) z^{\mathbf{k}} H^{(n)}_+ = \gamma_-(z)
H_\kk.\;$ It follows that every $W$ is in some $\Sigma_\kk^\sigma$.
\end{proof}

In fact, a more refined statement is true.

\begin{lem} For every $W \in \Grrr{n},\;$ there is a $\gamma_-(z) \in N^-$ and an $\kk$ so that $W = \gamma_-(z) H_\kk$. \end{lem}

\begin{proof} We use the Pressley-Segal identification of $H^{(n)}= L^2(S^1,\complex^n) \to L^2(S^1, \complex) = H$ given
by $(f_1,\ldots,f_n) \mapsto \widetilde{f}(\zeta) = \sum_{i=1}^n \zeta^{i-1} f_i(\zeta^n).$ We define $\check{W} = \bigcup_m W \cap \zeta^m H_-$..
  Choose an algebraic vector space basis of $\check{W}$ by considering the subset of $\check{W}$ consisting of elements of the
  form  $\zeta^s + \sum_{k=-\infty}^{s-1} a_k \zeta^k$, and choosing one such element for each possible value of $s$.
  Call $S$ the set of all the values of $s$ appearing in this construction. Denote by $w_s = \zeta^s + \sum_{k=-\infty}^{s-1} a_k \zeta^k$
  the chosen element so that our basis is $\BB = \{ w_s \colon s \in S\}$.. Let $H_S = \{ \sum_{s \in S} a_s \zeta^s\} $ be the Hilbert
  space generated by the $\zeta^s$. We may suppose that the orthogonal projection $W\to H_S$ sends $w_s \mapsto \zeta^s$
  (by using reduction of the basis $\BB$ to its reduced echelon form.). This induces an isomorphism $ W\cong H_S$.
   Since $z W \subset W,$ then $s\in S$ implies $s+n \in S$ (cf. p. 98 in \cite{Lu2000}). 
    There are $n$ elements $r_i$ in $S$ so that $r_i-n$ is not in $S$, and $S= \bigcup_{q \in \naturals} \{r_1 + qn, \ldots, r_n + qn\}$.
    Writing $r_i = n k_i + i -1 $ we conclude immediately that $H_S = H_\kk$.

   We can find the element $\gamma_-(z)$ we are seeking for using the isomorphism $W\to H_\kk$ by the following procedure.
   We define the smooth function $v_i(z) \in W$ to be the element in $W$ that projects to $(0,\ldots,z^{k_i},\ldots,0) = \zeta^{i-1+ n k_i}= \zeta^{r_i }$, so that we can write a basis of $W$  as $\BB' = \{ z^k v_i(z) \colon 1\leq n, k\geq0 \}$.  The matrix of smooth functions $v(z) =(v_1(z),\ldots,v_n(z))$ defines an element in $\LG$. Clearly $v(z) \cdot H_+$ = W. Define $\gamma_-(z) = v(z) \cdot z^{-\kk}$. Then $\gamma_-(z) \cdot H_{\kk} = W$. Since $v(z)$ projects to $(0,\ldots,z^{k_i},\ldots,0) \in H_\kk$ then $\gamma_i(z) = v_i(z) z^{-k_i}$ projects to $(0,\ldots,1,\ldots,0) \in H_+$, therefore no positive powers of $z$ appear in th expansion of $\gamma_-(z)$ and moreover the constant term in the Laurent expansion $\gamma_-(\infty)$ is upper triangular  .  We also have that  $\deg \det \gamma_-(z) + v.dim W = v.dim H_\kk$, and therefore $\deg\det\gamma_-(z)=0$. From this we conclude  that $\gamma_-(z) \in N^-$. \end{proof}

\begin{prop} The Grassmannian $\Grrr{n}$ admits a partition $$ \Grrr{n} = \coprod_\kk \Sigma_\kk^\sigma.$$ Moreover, each $\Sigma_\kk^\sigma$ is the union of the $\Sigma_\kk$'s for all permutations in the order of $\kk$, namely $$ \Sigma_\kk^\sigma = \coprod_{\epsilon \in S_n} \Sigma_{\epsilon(\kk)}.$$
\end{prop}

\begin{proof}  It is easy to see that for each permutation $\epsilon$ we have that $H_{\epsilon(\kk)} \in \Sigma_\kk^\sigma$
and hence $ \Sigma_\kk^\sigma \supseteq \bigcup_{\epsilon \in S_n}
\Sigma_{\epsilon(\kk)}.$ Since the whole $\Grrr{n}$ is the disjoint
union of $\Sigma_{\mathbf{j}}$ for all possible ${\mathbf{j}}$ and
since $N^- \subset \LGm$, it is enough to show that
$\Sigma^\sigma_\kk$ does not contain any $H_{\mathbf{j}}$ for a
${\mathbf{j}}$ that is not a permutation of $\kk$. To see this, we
associate a sequence $\omega(W)$ to each $W$ by $\omega_i(W) = \dim
(W \cap \zeta^i H_-)$. It is immediate to see that the sequence
$\omega$ is $\LGm$-invariant, and it nevertheless distinguishes
$H_\kk$ from $H_{\mathbf{j}}$ (because of the proof of the previous
proposition we can recover the ordered multi-index from $\omega$.)
\end{proof}

Finally, we are ready to restate Lemma \ref{defo:codimension} here
for the convenience of readers.
\begin{theo} The set $\Sigma_\kk$ is a contractible submanifold of $\Grrr{n}$ of codimension $$ cd(\kk) = \sum_{i<j} |k_i-k_j| - \varrho(\kk),$$
where $\varrho(\kk)$ is the number of inversions of $\kk$.
\end{theo}

\begin{proof} Let us define an open neighborhood of $H_\kk$. Let $L_0^- \subset \LGm$ be  is the subgroup of elements $\gamma(z)$ so that $\gamma(\infty)$ is the identity matrix. Consider the map $\LG \to \Grrr{n} \colon \gamma(z) \mapsto \gamma(z) \cdot H_\kk$. Let $U_\kk$ be the image under this map of $z^\kk L_0^- z^{-\kk}$. Clearly $H_\kk \in U_\kk$.

   Now $U_\kk$ we prove that is an open set. Here we return to the proof of Lemma 9.0.4 (1?). When we proved that $\gamma_-(z) \in N^-$ we should point out that  the same argument actually proves a little bit more. Indeed we have that $z^{-\kk} v_i(z) z^{\kk} \in H$ is of the form $(0,\dots,1,\ldots,0) + h_-(z) \in H$ where $h_-(z) \in H_-$. Therefore $\gamma_-(z)  \in z^\kk L_0^- z^{-\kk}$. Without loss of generality we assume $\kk=(0,\ldots,0)$, otherwise we shift by the appropiate $z^\kk$. In this case the Birkhoff factorization can be refined to state that every loop $\gamma(z)$ factorizes uniquely as $\gamma_(z)\gamma_+(z)$ where $\gamma_-(z) \in L_0^-$ and $\gamma_+(z) \in \LGp$. This implies that $L_0^- \cong U_{\mathbf{0}}$ is an open chart of $\Grrr{n} = \LG/\LGp$, and therefore $U_\kk$ is open. Also notice that this argument also implies that  $\Sigma_\kk \subseteq U_\kk$.

   We want to show then that the codimension of the inclusion $\Sigma_\kk \subseteq U_\kk$ is given by the formula of the statement of the theorem. Since we know that in the proof of Lemma 1 we actually have $\gamma_-(z) \in (z^\kk L_0^- z^{-\kk}) \cap N^-$. In fact $\LG \to \Grrr{n} \colon \gamma(z) \mapsto \gamma(z) \cdot H_\kk$ induces an identification $(N^- \cap z^\kk L_0^- z^{\kk}) \to \Sigma_\kk$. We claim that the multiplication in $\LG$ indices an identification $(N^- \cap z^\kk L_0^- z^{-\kk} ) \times (N^+ \cap z^\kk L_0^- z^{-\kk} ) \to z^\kk L_0^- z^{-\kk} = U_\kk$, where $N^+$ is just as $N^-$ except that we talk of lower triangular matrices.

   All that remains then is to compute the dimension of $N^+ \cap z^\kk L_0^- z^{-\kk} $. This is done as follows. By taking Laurent expansions of the entries of en element $\gamma(z) \in (N^+ \cap z^\kk L_0^- z^{-\kk})$ we conclude that for $i<j$ then $\gamma_{ii}(z)=1$, $\gamma_{ij}(z) = \sum_{l=1}^{k_i-k_j-1} a_l z^l$ and $\gamma_{ji}=\sum_{l=0}^{k_j-k_i-1} a_l z^l$. By counting coefficients we obtain the desired formula. \end{proof}


{\textit{Remark.}} Notice that the proof of the previous theorem actually shows more. It shows that the codimension of elements $\gamma(z)$ of multi-index $\kk$ inside $\LG$ is given by the same formula. This is done by considering $(N^- \cap z^\kk L_0^- z^{-\kk} ) \times (N^+ \cap z^\kk L_0^- z^{-\kk} )\times \LGp$ instead of $(N^- \cap z^\kk L_0^- z^{-\kk} ) \times (N^+ \cap z^\kk L_0^- z^{-\kk} )$.

X. X. Chen, Department of Mathematics, University of Wisconsin.\\
$~~~~~~~~~~~~~~xiu@math.wisc.edu$\\

G. Tian, Department of Mathematics, Princeton University and
Beijing University. \\
$~~~~~~~~~~~~~~tian@math.princeton.edu$
\end{document}